\documentclass[12pt,a4paper]{article}
\usepackage{amsthm,amssymb,amsfonts,amsmath,amscd,times1}
\usepackage{color}
\usepackage{hyperref}
\usepackage{breakurl}
\usepackage{mysec}
\usepackage[normalem]{ulem}

\makeatletter
\renewcommand{\@listI}%
{\leftmargin=\leftmargini
\partopsep=0pt
\topsep=3pt
\itemsep=3pt
\labelwidth=\leftmargini }
\renewcommand{\@listii}{\setlength{\topsep}{3pt}
}  
 \makeatother

\numberwithin{equation}{section}
\theoremstyle{plain}
\newtheorem{theorem}{Theorem}[section]
\newtheorem{proposition}[theorem]{Proposition}
\newtheorem{lemma}[theorem]{Lemma}
\newtheorem{corollary}[theorem]{Corollary}

\theoremstyle{remark}
\newtheorem{remark}{Remark}[section]

\newtheorem{example}{Example}[section]

\theoremstyle{definition}
\newtheorem{definition}{Definition}[section]
\newtheorem{condition}{Condition}[section]
\newtheorem*{condit}{Condition \ref{condition5}\mypp$'$}

\newcommand{\Ad}{\mathop{\mathrm{Ad}}\nolimits}
\newcommand{\ad}{\mathop{\mathrm{ad}}\nolimits}
\newcommand{\Div}{\mathop{\mathrm{div}}\nolimits}

\newcommand{\Id}{\mathop{\mathrm{id}}\nolimits}
\newcommand{\supp}{\mathop{\mathrm{supp}}\nolimits}

\newcommand{\Diff}{\mathrm{Diff}}
\newcommand{\Vect}{\mathrm{Vect}}
\newcommand{\Euc}{\mathcal{E}}
\newcommand{\SO}{\mathcal{SO}}
\newcommand{\Ort}{\mathcal{O}}
\newcommand{\const}{\mathrm{const}}
\newcommand{\rd}{\mathrm{d}}
\newcommand{\re}{\mathrm{e\myp}}

\newcommand{\CDD}{\mathop{\raisebox{-.11pc}{\Large$\cdot$}}}
\newcommand{\RR}{\mathbb{R}}
\newcommand{\NN}{\mathbb{N}}
\newcommand{\ZZ}{\mathbb{Z}}

\newcommand{\calB}{\mathcal{B}}

\newcommand{\frakX}{\mathfrak{X}}

\newcommand{\mucl}{\mu_{\mathrm{cl}}}
\newcommand{\hatgamma}{{\hat{\gamma\myp}\myn{}}}

\newcommand{\myp}{\mbox{$\:\!$}}
\newcommand{\mypp}{\mbox{$\;\!$}}

\newcommand{\myn}{\mbox{$\;\!\!$}}
\newcommand{\mynn}{\mbox{$\:\!\!$}}

\begin{document}

\title{Cluster point processes on manifolds}
\author{Leonid Bogachev\myp$^{\rm a}$ and Alexei Daletskii\myp$^{\rm b}$}
\date{\small \it $^{\rm a}$Department of Statistics, University of Leeds, Leeds LS2 9JT,
UK.\\
 E-mail: {\tt L.V.Bogachev@leeds.ac.uk}\\[.2pc]
$^{\rm b}$Department of Mathematics, University of York, York YO10
5DD, UK.\\ E-mail: {\tt ad557@york.ac.uk}} \maketitle

\begin{abstract}
\vspace{7pt}
The probability distribution $\mucl $ of a general
cluster point process in a Riemannian manifold $X$ (with independent
random clusters attached to points of a configuration with
distribution $\mu$) is studied via the projection of an auxiliary
measure $\hat{\mu}$ in the space of configurations
$\hatgamma=\{(x,\bar{y})\}\subset X\times\mathfrak{X}$, where $x\in
X$ indicates a cluster ``centre'' and
$\bar{y}\in\mathfrak{X}:=\bigsqcup_{\myp n}\mynn X^n$ represents a
corresponding cluster relative to $x$.
We show that the measure $\mucl $ is quasi-invariant with respect to
the group $\Diff_{0}(X)$ of compactly supported diffeomorphisms of
$X$, and prove an integration-by-parts formula for $\mucl$. The
associated equilibrium stochastic dynamics is then constructed using
the method of Dirichlet forms. General constructions are illustrated
by examples including Euclidean spaces, Lie groups, homogeneous
spaces, Riemannian manifolds of non-positive curvature and metric
spaces. The paper is an extension of our earlier results
for Poisson cluster measures [J.~Funct.\ Analysis 256 (2009)
432--478] and for Gibbs cluster measures
[\url{http://arxiv.org/abs/1007.3148}], where different projection
constructions were utilised.


\medskip
\noindent \textit{MSC\,2010:} Primary 58J65; Secondary 31C25, 46G12,
60G55, 70F45
\end{abstract}
%
%

\tableofcontents

\section{Introduction}\label{sec:1}

The concept of a \textit{configuration space} (over a suitable
Riemannian manifold) is instrumental in the description of various
types of multi-particle structures and naturally appears in many
areas of mathematics and mathematical physics (e.g., theory of
random point processes, statistical mechanics, quantum field theory,
representation theory) and applied sciences (e.g., chemical physics,
image processing, spatial ecology, astronomy, etc.).

Despite not possessing any Banach manifold structure, configuration
spaces have many features of proper manifolds and can indeed be
endowed with ``manifold-like'' structures (see \cite{VGG} and also
\cite{AKR1,AKR2}). As it turns out, the way to do it depends heavily
on the choice of a suitable probability measure $\mu$ on the
configuration space $\varGamma_{X}$ (over the mani\-fold $X$). Such
a choice is often suggested by a physical system under study, but,
in order to furnish a meaningful analytical framework, the measure
$\mu$ must satisfy certain regularity properties, such as the
$\Diff_{0}$-quasi-invariance with respect to the action of certain
diffeomorphism group and/or an integration-by-parts (IBP) formula.
Hence, it is not surprising that the study of the configuration
space as a measure space $(\varGamma_X,\mu)$ requires tools and
techniques at the interface of geometric analysis and measure
theory.
According to this paradigm, it is important
(i) to prove the quasi-invariance  and IBP formulae for a wide class
of measures $\mu$ arising in applications, and (ii) to study the
dependence of the properties of the measure $\mu$ on the topology
and geometry of the underlying manifold $X$ and their interplay with
the multi-particle structure of the space $\varGamma_{X}$.

Such a programme has been implemented for the Poisson and Gibbs
measures on $\varGamma_{X}$ in the case where $X=\RR^d$ (see
\cite{AKR1,AKR2,AKR3,ADKal,AD} and further references therein). The
present paper is a step towards realisation of this programme for
the important class of (in general, non-Gibbsian) measures on
$\varGamma_{X}$ emerging as distributions of \textit{cluster point
processes} in $X$. Intuitively, a cluster point process is obtained
by generating random clusters around points of the background
configuration of cluster ``centres'' \cite{DVJ1}. Cluster models are
well known in the general theory of random point processes
\cite{CI,DVJ1} and are widely used in numerous applications ranging
from neurophysiology (nerve impulses) and ecology (spatial
aggregation of species) to seismology (earthquakes) and cosmology
(constellations and galaxies); see \cite{BD3,CI,DVJ1} for some
references to original papers.

In our earlier papers \cite{BD1,BD2,BD3,BD4}, we have developed a
projection construction of Poisson and Gibbs cluster processes in a
Euclidean space $X=\RR^d$, based on the representation of their
probability distributions (i.e., the corresponding measures on the
configuration space $\varGamma_X$) as the push-forward (image) of
suitable auxiliary measures on a more complex configuration space
$\varGamma_{\mathfrak{X}}$ over the disjoint-union space
$\mathfrak{X}:=\bigsqcup_{\myp n} \mynn X^{n}$, with ``droplet''
points $\bar{y}\in \mathfrak{X}$ representing individual clusters
(of variable size). Such an approach allows one to adapt the ideas
of analysis and geometry on configuration spaces developed earlier
by Albeverio, Kondratiev and R\"ockner \cite{AKR1,AKR2} for plain
(i.e., non-cluster) Poisson and Gibbs measures in $\varGamma_X$, and
to obtain results including the $\Diff_{0}$-quasi-invariance and IBP
formula.

In the present paper, we extend the projection approach to the case
of cluster measures on general Riemannian manifolds $X$ and with
arbitrary centre processes. In so doing, suitable smoothness
properties of the distribution of individual cluster are required,
but it should be stressed that \textit{no smoothness of the centre
process is needed}. That is to say, attaching ``nice'' clusters to
points of a centre configuration acts as smoothing of the entire
cluster process. To an extent, this may be thought of as an
infinite-dimensional analogue of the well-known fact that the
convolution of two measures in $\RR^d$ is absolutely continuous
provided that at least one of those measures is such. However, this
analogy should not be taken too far, since the relationship between
centres and clusters is asymmetric (the latter are attached to the
former, but not vice versa); in particular, as it turns out,
smoothness of the centre process alone is not sufficient for the
smoothness of the resulting cluster process. Let us point out that
the results of the present paper are new even in the case of Poisson
and Gibbs cluster point processes in $\RR^d$, where our new approach
allows one to handle more general models, for instance with the
probability distribution of centres given by a Poisson measure with
a \textit{non-smooth intensity}, or by a Gibbs measure with a
\textit{non-smooth interaction potential}.

The structure of the paper is as follows. In Section \ref{sec:2} we
introduce the general framework of configuration spaces and measures
on them and discuss a ``fibre bundle'' structure of the
configuration space over a product manifold. Section \ref{sec:3} is
devoted to the projection construction of cluster measures $\mu
_{cl}$. Here we derive necessary conditions for the cluster measure
$\mucl$ to be well defined (i.e., with no multiple and accumulation
points) and study the existence of moments. In Section \ref{sec:4}
we prove the $\Diff_{0}$-quasi-invariance and an IBP formula for
$\mucl$. Furthermore, for a general cluster measure we are able to
construct the corresponding Dirichlet form and to prove its
closedness, which implies in a standard way the existence of the
corresponding equilibrium stochastic dynamics. In Section
\ref{sec:5} we discuss examples of cluster distributions that can be
generated in a natural way via certain manifold structures, such as
the group action in the case of homogeneous manifolds and a metric
structure for general Riemannian manifolds. Finally, the Appendix
contains some additional
or supporting
material.

\section{Configuration spaces and measures}\label{sec:2}

\subsection{General setup: probability measures on configuration spaces}
\label{sec:2.1}

Let $X$ be a Polish space
equipped with the Borel $\sigma$-algebra ${\mathcal{B}}(X)$
generated by the open sets. Denote $\overline{\ZZ}_+:=
\ZZ_+\cup\{\infty\}$, where $\ZZ_{+}:=\{0,1,2,\dots \}$, and
consider the space ${\mathfrak{X}}$ built from all Cartesian powers
of $X$, that is, the disjoint union
\begin{equation} \label{eq:calX}
\mathfrak{X}:= {\textstyle\bigsqcup_{\myp n\in\overline{\ZZ}_{+}}}
X^{n},
\end{equation}
including $X^{0}=\{\emptyset \}$ and the space $X^\infty$ of
infinite sequences $(x_1, x_2,\dots)$. That is,
$\bar{x}=(x_{1},x_{2},\dots)\in\mathfrak{X}$ if and only if
$\bar{x}\in X^{n}$ for some $n\in\overline{\ZZ}_{+}$. We take the
liberty to write $x_{i}\in\bar{x}$ if $x_{i}$ is a coordinate of the
``vector'' $\bar{x}$.
The space $\mathfrak{X}$ is endowed with the natural disjoint union topology
induced by the topology in $X$. 

\begin{remark}\label{rm:compact}
Note that a set $B\subset{\mathfrak{X}}$ is compact if and only if
$B=\bigsqcup_{\myp n=0}^{N} B_{n}$, where $N<\infty$ and $B_{n}$ are
compact subsets of $X^{n}$, respectively.
\end{remark}

\begin{remark}
$\mathfrak{X}$ is a Polish space as a disjoint union of Polish spaces.
\end{remark}

Denote by ${\mathcal{N}}(X)$ the space of
$\overline{\ZZ}_{+}$-valued measures on ${\mathcal{B}}(X)$ with
countable (i.e., finite or countably infinite) support.
Consider the natural projection
\begin{equation}  \label{eq:pr0}
{\mathfrak{X}}\ni \bar{x}\mapsto {\mathfrak{p}}(\bar{x}):=\sum_{x_{i}\in
\bar{x}}\delta_{x_{i}}\in {\mathcal{N}}(X),
\end{equation}
where $\delta_{x}$ is the Dirac measure at point $x\in X$.
That is to say, under the map ${\mathfrak{p}}$ each vector from
${\mathfrak{X}}$ is ``unpacked'' into its components to yield a
countable aggregate of (possibly multiple) points in $X$, which can
be interpreted as a generalised configuration $\gamma$,
\begin{equation}  \label{eq:pr}
\mathfrak{p}(\bar{x})\leftrightarrow\gamma :={\textstyle\bigsqcup
\limits_{x_{i}\in \bar{x}}}\{x_{i}\},\qquad
\bar{x}=(x_{1},x_{2},\dots)\in\mathfrak{X}.
\end{equation}

In what follows, we
interpret the notation $\gamma $ either as an aggregate of points in
$X$ or as a $\overline{\ZZ}_{+}$-valued measure or both, depending
on the context. Even though generalised configurations are not,
strictly speaking, subsets of $X$ (because of possible
multiplicities), it is convenient to use set-theoretic notations,
which should not cause any confusion. For instance, we write $\gamma
\cap B$ for the restriction of configuration $\gamma $ to a subset
$B\in {\mathcal{B}}(X)$. For a function $f:X\rightarrow \mathbb{R}$
we denote
\begin{equation}\label{eq:f-gamma}
\langle f,\gamma \rangle :=\sum_{x_{i}\in \gamma}f(x_{i})\equiv
\int_{X}f(x)\,\gamma (\rd{x}).
\end{equation}
In particular, if $\mathbf{1}_{B}(x)$ is the indicator function of a
set $B\in {\mathcal{B}}(X)$ then $\langle \mathbf{1}_{B},\gamma
\rangle=\gamma (B)$ is the total number of points (counted with
their multiplicities) in $ \gamma\cap B$.

\begin{definition}\label{def:gen}
A \textit{configuration space} $\varGamma_{X}^{\sharp}$ is the set
of generalised configurations $\gamma $ in $X$, endowed with the
\textit{cylinder $\sigma $-algebra}
$\mathcal{B}(\varGamma_{X}^{\sharp})$ generated by the class of
cylinder sets $C_{B}^{\myp n}:=\{\gamma
\in\varGamma_{X}^{\sharp}:\gamma (B)=n\}$, \,$B\in \mathcal{B}(X)$,
\,$n\in\ZZ_{+}$\myp.
\end{definition}

\begin{remark}
It is easy to see that the map
$\mathfrak{p}:{\mathfrak{X}}\rightarrow \varGamma_{X}^{\sharp}$
defined by formula (\ref{eq:pr}) is measurable.
\end{remark}

In fact, conventional theory of point processes (and their
distributions as probability measures on configuration spaces)
usually rules out the possibility of accumulation points or multiple
points (see, e.g., \cite{DVJ1}).

\begin{definition}\label{def:proper}
A configuration $\gamma \in \varGamma_{X}^{\sharp}$ is said to be
\emph{locally finite} if $\gamma(B)<\infty$ for any compact set
$B\subset X$. A configuration $\gamma \in \varGamma_{X}^{\sharp}$ is
called \textit{simple} if $\gamma (\{x\})\leq 1$ for each $x\in X$.
A configuration $\gamma \in \varGamma_{X}^{\sharp}$ is called
\emph{proper} if it is both locally finite and simple. The set of
proper configurations will be denoted by $\varGamma_{X}$ and called
the \textit{proper configuration space} over $X$. The corresponding
$\sigma $-algebra ${\mathcal{B}}(\varGamma_{X})$ is generated by the
cylinder sets $\{\gamma \in \varGamma_{X}:\gamma (B)=n\}$ \,($B\in
{\mathcal{B}}(X)$, \,$n\in \ZZ_{+}$).
\end{definition}

Like in the standard theory based on proper configuration spaces
(see, e.g., \cite[\S\myp6.1]{DVJ1}), every probability measure $\mu
$ on the generalised configuration space $\varGamma_{X}^{\sharp}$
can be characterised by its Laplace functional (cf.\ \cite{BD3})
\begin{equation}\label{eq:LAPLACE}
L_{\mu}(f):=\int_{\varGamma_{X}^{\sharp}}\re^{-\langle
f,\myp\gamma\rangle}\,\mu (\rd\gamma),\qquad f\in
{\mathrm{M}}_{+}(X),
\end{equation}
where ${\mathrm{M}}_{+}(X)$ is the class of measurable non-negative
functions on $X$.

\subsection{Cluster point processes}\label{sec:2.2}

Let us recall the notion of a general cluster point process (CPP).
Its realisations are constructed in two steps: (i) a background
random configuration of (invisible) ``centres'' is obtained as a
realisation of some point process $\gamma_{\mathrm{c}}$ governed by
a probability measure $\mu$ on $\varGamma_{X}^{\sharp}$, and (ii)
relative to each centre $x\in\gamma_{\mathrm{c}}$, a set of
observable secondary points (referred to as a \emph{cluster} centred
at~$x$) is generated according to a point process
$\gamma_{x}^{\prime}$ with distribution $\mu_{x}$ on
$\varGamma_{X}^{\sharp}$ ($x\in X$). The resulting (countable)
assembly of random points, called the \emph{cluster point process},
can be expressed symbolically as
\begin{equation*}
\gamma ={\textstyle\bigsqcup\limits_{x\in
\gamma_{\mathrm{c}}}}\,\gamma_{x}^{\prime}\in
\varGamma_{X}^{\sharp},
\end{equation*}
where the disjoint union signifies that multiplicities of points are
taken into account. More precisely, assuming that the family of
secondary processes $\gamma_{x}^{\prime}$ is measurable as a
function of $x\in X$, the integer-valued measure corresponding to a
CPP realisation $\gamma $ is given by
\begin{equation}\label{eq:cluster-gamma}
\gamma(B)=\int_{X}\gamma_{x}^{\myp\prime}(B)\,\gamma_{\mathrm{c}}
(\rd{x})=\sum_{x\in\gamma_{\mathrm{c}}}\gamma_{x}^{\prime}(B),\qquad
B\in {\mathcal{B}}(X).
\end{equation}
We denote by $\mucl$ the probability measure on
$(\varGamma_{X}^{\sharp},\mathcal{B}(\varGamma_{X}^{\sharp}))$ that
    governs the CPP described by formula
\textup{(\ref{eq:cluster-gamma})}.

\begin{remark}
Unlike the standard CPP theory where sample configurations are
\emph{presumed} to be almost surely (a.s.)\ locally finite (see,
e.g., \cite[Definition 6.3.I]{DVJ1}), the description of CPP given
above only implies that its configurations $\gamma$ are countable
aggregates in $X$, but possibly with multiple and/or accumulation
points, even if the background point process $\gamma_{\mathrm{c}}$
is proper. Thus, the CPP measure $\mucl $ is defined on the space
$\varGamma_{X}^{\sharp}$ of \emph{generalised} configurations.
However, developing the differential analysis on configuration
spaces in the spirit of Albeverio, Kondratiev and R\"ockner
\cite{AKR1,AKR2} demands that measures under study
are actually supported on the proper configuration space
$\varGamma_{X}$. We shall address this issue for the general cluster
measure $\mucl$ in Section \ref{sec:3.2} below and give sufficient
conditions in order that $\mucl$-almost all (a.a.)\ configurations
are proper (see \cite{BD3,BD4} for the cases of the Poisson and
Gibbs CPPs, respectively).
\end{remark}

The distribution $\mu_{x}$ of the inner-cluster process
$\gamma_{x}^{\prime}$ determines a probability measure $\eta_{x}$ in
${\mathfrak{X}}$ symmetric with respect to permutations of
coordinates. Conversely, $\mu_{x}$ is a push-forward of the measure
$\eta_{x}$ under the projection map
${\mathfrak{p}}:{\mathfrak{X}}\to\varGamma_{X}^{\sharp}$ defined by
(\ref{eq:pr}), that is,
\begin{equation}\label{eq:p*eta}
\mu_{x}={\mathfrak{p}}^{\ast}\eta_{x}\equiv \eta_{x}\circ
{\mathfrak{p}}^{-1}.
\end{equation}

The following fact is well known for CPPs without accumulation
points (see, e.g., \cite[\S\myp6.3]{DVJ1}); its proof in the general
case is essentially the same (see \cite[Proposition 2.5]{BD3}).

\begin{proposition}
\label{pr:cluster} The Laplace functional of the cluster measure
$\mucl$
is given by
\begin{equation}\label{laplace-G}
L_{\mucl}(f)=\int_{\varGamma_{X}^{\sharp}}\prod_{x\in
\gamma}\biggl(\int_{\mathfrak{X}}\exp \biggl(-\sum_{y_{i}\in
\bar{y}}f(y_{i})\Bigr)\,\eta_{x}(\rd\bar{y})\biggr)\,\mu(\rd\gamma),
\qquad f\in {\mathrm{M}}_{+}(X).
\end{equation}
\end{proposition}
\noindent

\subsection{Measures on configurations in product spaces}\label{sec:2.3}

In this section, we develop a general construction of measures on
configurations in product spaces, which will be useful later on.

Let $Y$ be a Polish space equipped with the Borel $\sigma$-algebra
$\mathcal{B}(Y)$. Consider the product space $Z:=X\times Y$ endowed
with the product $\sigma$-algebra
$\mathcal{B}(Z)=\mathcal{B}(X)\otimes \mathcal{B}(Y)$, and the
corresponding configuration space $\varGamma_Z^{\sharp}$. Let
\begin{equation*}
p_{X}(z):=x,\qquad p_{Y}(z):=y,\qquad  z=(x,y)\in X\times Y,
\end{equation*}
denote the natural projections onto the spaces $X$ and $Y$,
respectively. The maps $p_{X}$ and $p_Y$ can be extended to the
configuration space $\varGamma_{Z}^{\sharp}$:
\begin{align*}
\varGamma_{Z}^{\sharp}\ni\hatgamma&\mapsto
p_X(\hatgamma):=\textstyle{\bigsqcup\limits_{z\in\hatgamma}}
\,p_{X}(z)\in \varGamma_{X}^{\sharp},\\
\varGamma_{Z}^{\sharp}\ni\hatgamma&\mapsto
p_Y(\hatgamma):=\textstyle{\bigsqcup\limits_{z\in\hatgamma}}
\,p_{Y}(z)\in \varGamma_{Y}^{\sharp}.
\end{align*}
For each $\gamma\in \varGamma_{X}^{\sharp}$, consider the space
$Y^{\gamma}:=p_{X}^{-1}(\gamma)$ (i.e., the \textit{fibre} at
$\gamma$), which can be identified with the corresponding Cartesian
product of
identical copies of the space $Y$,
\begin{equation*}
Y^{\gamma}=\prod_{x\in \gamma} Y_{x},\qquad Y_{x}=Y.
\end{equation*}
Therefore, each configuration $\hatgamma\in \varGamma_{Z}^{\sharp}$
can be represented in the form
\begin{equation}\label{eq:hatgamma}
\hatgamma=(\gamma,y^{\gamma})={\textstyle\bigsqcup\limits_{x\in\gamma}\{(x,y_x)\}},
\end{equation}
with
\begin{equation}\label{eq:pX,pY}
\gamma=p_{X}(\hatgamma)\in \varGamma_{X}^{\sharp},\qquad
y^{\gamma}:= (y_{x})_{x\in\gamma}\in Y^{\gamma}.
\end{equation}
More formally, a one-to-one correspondence between $x\in\gamma$ and
$y_x\in y^\gamma$ is described by the relations
\begin{equation}\label{eq:x<->y}
x=p_X\bigl(p_Y^{-1}(y_x)\cap \hatgamma\bigr),\qquad
y_x=p_Y\bigl(p_X^{-1}(x)\cap\hatgamma\bigr).
\end{equation}

Let each space $Y_{x}$, $x\in X$, be equipped with a probability
measure $\eta_{x}$. We assume that the family $\{\eta_{x},\,x\in
X\}$ satisfies the following ``weak measurability'' condition.

\begin{condition}\label{condition1}
For any Borel set $B\in\mathcal{B}(Y)$, the map
\begin{equation*}
X\ni x\mapsto \eta_{x}(B)\in[0,1]
\end{equation*}
is measurable with respect to $\mathcal{B}(X)$.
\end{condition}

\begin{lemma}\label{lm:meas}
Under Condition \textup{\ref{condition1}}, for any function
$f\in\mathrm{M}_+(Z)$ the integral
\begin{equation}\label{eq:int_Y}
\int_{Y}\re^{- f(x,\myp y)}\, \eta_{x}(\rd{y})
\end{equation}
is a measurable function of\/ $x\in X$.
\end{lemma}
\proof It is sufficient to consider functions of the form
$f(x,y)=\mathbf{1}_{A}(x)\cdot\mathbf{1}_{B}(y)$ with Borel sets
$A\in\mathcal{B}(X)$, $B\in\mathcal{B}(Y)$. In that case, the
integral in (\ref{eq:int_Y}) is reduced to
\begin{align*}
\left(\re^{- \mathbf{1}_{A}(x)}-1\right)\int_{B}\,
\eta_{x}(\rd{y})+\int_{Y} \eta_{x}(\rd{y})=\left(\re^{-
\mathbf{1}_{A}(x)}-1\right) \eta_x(B) +1,
\end{align*}
which is measurable in $x$ owing to Condition~\ref{condition1}.
\endproof

For $\gamma\in\varGamma_{X}^{\sharp}$, define the corresponding
product measure on the space $Y^\gamma$,
\begin{equation}\label{eq:eta-gamma}
\eta^{\gamma}(\rd y^\gamma):=
{\textstyle\bigotimes\limits_{x\in\gamma}}\,\eta_{x}(\rd y_x).
\end{equation}
Furthermore, let $\mu$ be a probability measure on
$\varGamma_{X}^{\sharp}$. Let us now define a probability measure
$\hat{\mu}$ on $\varGamma_{Z}^{\sharp}$ as a skew product (see
(\ref{eq:hatgamma}))
\begin{equation}\label{eq:mu-hat}
\hat{\mu}(\rd\hatgamma):=\eta^{\gamma}(\rd{y}^{\gamma})
\,\mu(\rd\gamma),\qquad \hatgamma=(\gamma,y^\gamma)\in
\varGamma_{Z}^{\sharp}.
\end{equation}
More precisely, definition (\ref{eq:mu-hat}) can be rewritten in an
integral form,
\begin{equation}\label{eq:mu-hat-int}
\int_{\varGamma_{Z}^{\sharp}}F(\hatgamma)\,\hat{\mu}(\rd\hatgamma)
=\int_{\varGamma_{X}^{\sharp}}\!\left(\int_{Y^{\gamma}}F(\gamma,y^{\gamma})\,
\eta^{\gamma}(\rd y^{\gamma})\right)\mu(\rd\gamma), \qquad
F\in\mathrm{M}_+(\varGamma_{Z}^{\sharp}).
\end{equation}

\begin{remark}
The measure $\hat{\mu}$ can be viewed as a measure on a
\textit{marked configuration space}, specified by the ground
configurations $\gamma\in \varGamma_{X}^{\sharp}$ with distribution
$\mu(\rd\gamma)$ and the corresponding configurations of marks
$y^\gamma=(y_x)_{x\in\gamma}\in Y^\gamma$ with distribution
$\eta^\gamma(\rd y^\gamma)$ (cf.\ \cite[\S\myp6.4]{DVJ1}).
\end{remark}

\begin{proposition}\label{pr:skew_product}
The measure $\hat{\mu}$ is a probability measure on
$\varGamma_{Z}^{\sharp}$ uniquely determined by its Laplace
functional given by the formula
\begin{equation}\label{eq:L=L}
L_{\hat\mu}(f)=L_{\mu}(\bar{f}\mypp),\qquad f\in\mathrm{M}_+(Z),
\end{equation}
where \begin{equation}\label{eq:f_bar} \bar{f}(x):=-\log
\int_{Y}\re^{- f(x,\myp y)}\,\eta_{x}(\rd{y})\ge 0,\qquad x\in X.
\end{equation}
In particular, the integral formula \textup{(\ref{eq:mu-hat-int})}
is valid, with the internal integral being measurable as a function
of $\gamma$ with respect to the $\sigma$-algebra
$\mathcal{B}(\varGamma_X^\sharp)$.
\end{proposition}

\proof Note that the function $\bar{f}(x)$ defined in
(\ref{eq:f_bar}) is measurable with respect to $\mathcal{B}(X)$
according to Condition \ref{condition1} and Lemma~\ref{lm:meas};
hence the right-hand side of equation (\ref{eq:L=L}) is well
defined. Using the skew-product definition (\ref{eq:mu-hat}) we have
\begin{align*}
L_{\hat\mu}(f)&=\int_{\varGamma_{Z}^{\sharp}}\re^{-\langle
f,\myp\hatgamma\rangle}\,\hat{\mu}
(\rd\hatgamma)\\
&=\int_{\varGamma_{X}^{\sharp}} \left(\int_{Y^\gamma}\re^{-\langle
f,\myp\hatgamma\rangle}\, \eta^\gamma(\rd y^\gamma)\right)
\mu(\rd\gamma)\\
&=\int_{\varGamma_{X}^{\sharp}}\left(\int_{Y^\gamma}\exp
\left(-\sum_{x\in\gamma} f(x,y_x)\right)
{\textstyle\bigotimes\limits_{x\in\gamma}}\eta_{x} (\rd y_x)\right)
\mu(\rd\gamma)\\
&=\int_{\varGamma_{X}^{\sharp}}\prod_{x\in\gamma}\left(\int_{Y_x}
\re^{-f(x,\myp y_x)}\, \eta_{x}(\rd y_x)\right)
\mu(\rd\gamma)\\
&=\int_{\varGamma_{X}^{\sharp}}\exp\left(\sum_{x\in\gamma}\log\int_{Y_x}\re^{-f(x,\myp
y_x)}\,\eta_{x}(\rd y_x)\right)
\mu(\rd\gamma)\\
&=\int_{\varGamma_{X}^{\sharp}}\re^{-\langle
\bar{f},\myp\gamma\rangle}\,\mu(\rd\gamma)=L_{\mu}(\bar{f}\mypp),
\end{align*}
and equation (\ref{eq:L=L}) is established. In particular, taking
$f(x,y)\equiv 0$ we readily obtain $\bar{f}(x)\equiv 0$ and hence
$$
\hat{\mu}(\varGamma_{Z}^{\sharp})=\int_{\varGamma_{Z}^{\sharp}}\hat{\mu}
(\rd\hatgamma)=L_{\hat\mu}(0)=L_{\mu}(0)=1,
$$
as required. Finally, using in (\ref{eq:L=L}) test functions of the
form
$$
f(x,y)=\sum_{i=1}^k s_i \myp\mathbf{1}_{A_i}(x)
\myp\mathbf{1}_{B_i}(y),\quad s_i\ge0,\ \ A_i\in\mathcal{B}(X),\ \
B_i\in\mathcal{B}(Y),
$$
one can show in the usual fashion that equation
(\ref{eq:mu-hat-int}) holds for the indicators of the cylinder sets
$\bigcap_{\myp i=1}^{\myp k} C_{\mynn A_i\times B_i}^{\myp n_i}$,
which in turn implies (\ref{eq:mu-hat-int}) for any measurable
function $F(\gamma,y^\gamma)$. The proof is complete.
\endproof

\begin{remark}
A direct proof of the measurability of the internal integral on the
right-hand side of equation (\ref{eq:mu-hat-int}) appears to be
quite involved. For illustration, we give such a proof in the case
(see Appendix \ref{ap:1}) based on a \textit{measurable} indexation
of the ground configurations $\gamma$, which is known to be
available in the case of proper $\gamma$'s (i.e., with a.s.\ no
multiple or accumulation points). Note that the general theory of
marked point processes (measures) is usually based on the assumption
of properness of the ground process $\gamma$ (see
\cite[\S\myp6.4]{DVJ1}).
\end{remark}

\section{Cluster measures on configuration spaces}\label{sec:3}

From now on, the measure $\mu$ on the space $X$, used in the
construction of the configurational measure $\hat{\mu}$ (see
(\ref{eq:mu-hat})), will represent the distribution of the cluster
centres; furthermore, the role of the space $Y$ will be played by
$\mathfrak{X}\equiv {\textstyle\bigsqcup_{\myp
n\in\overline{\ZZ}_{+}}} X^{n}$ (see (\ref{eq:calX})), with generic
elements $\bar{y}\in\frakX$. To emphasise such a choice, we change
the general notation of the space $Z=X\times Y$ to
$\mathcal{Z}:=X\times \mathfrak{X}$.

\subsection{Projection construction of the cluster measure}\label{sec:3.1}

Recall that the ``unpacking'' map
$\mathfrak{p}:\mathfrak{X}\rightarrow \varGamma_{X}^{\sharp}$ is
defined in (\ref{eq:pr}), and consider a map
$\mathfrak{q}:\mathcal{Z}\rightarrow \varGamma_{X}^{\sharp}$ acting
by the formula
\begin{equation}\label{proj1}
\mathfrak{q}(x,\bar{y}):=\mathfrak{p}(\bar{y})=
{\textstyle\bigsqcup\limits_{y_{i}\in \bar{y}}}\{y_{i}\},\qquad
(x,\bar{y})\in \mathcal{Z}.
\end{equation}
In the usual ``diagonal'' way, the map $\mathfrak{q}$ can be lifted
to the configuration space $\varGamma_{\mathcal{Z}}^{\sharp}$:
\begin{equation}\label{eq:proj}
\varGamma_{\mathcal{Z}}^{\sharp}\ni {\hatgamma}\mapsto
\mathfrak{q}({\hatgamma}):={\textstyle\bigsqcup\limits_{z\in
\hatgamma}}\mathfrak{q}(z)\in \varGamma_{X}^{\sharp}.
\end{equation}

\begin{proposition}\label{pr:q-meas}
The map $\mathfrak{q}:\varGamma_{\mathcal{Z}}^{\sharp}\rightarrow
\varGamma_{X}^{\sharp}$ defined by \textup{(\ref{eq:proj})} is
measurable.
\end{proposition}

\proof Observe that $\mathfrak{q}$ can be represented as a
composition
\begin{equation}\label{eq:circ}
\mathfrak{q}=\mathfrak{p}\circ
p_{\mathfrak{X}}:\varGamma_{\mathcal{Z}}^{\sharp}\stackrel{p_{\mathfrak{X}}}{\longrightarrow
}\varGamma_{\mathfrak{X}
}^{\sharp}\stackrel{\mathfrak{p}}{\longrightarrow
}\varGamma_{X}^{\sharp},
\end{equation}
where the maps $p_{\mathfrak{X}}$ and $\mathfrak{p}$ are defined,
respectively, by
\begin{align}
\label{eq:pY} \varGamma_{\mathcal{Z}}^{\sharp}\ni\hatgamma& \mapsto
p_{\mathfrak{X}}({\hatgamma}):={\textstyle\bigsqcup\limits_{(x,\myp\bar{y})
\in \myp\hatgamma}}\{\bar{y}\}\in
\varGamma_{\mathfrak{X}}^{\sharp}, \\
\label{eq:pr2} \varGamma_{\mathfrak{X}}^{\sharp}\ni \bar{\gamma}&
\mapsto
\mathfrak{p}(\bar{\gamma}):={\textstyle\bigsqcup\limits_{\bar{y}\in
\bar{\gamma}}}\mathfrak{p}(\bar{y})\in \varGamma_{X}^{\sharp}.
\end{align}
For $p_{\mathfrak{X}}:\varGamma_{\mathcal{Z}}^{\sharp}\rightarrow
\varGamma_{\mathfrak{X}}^{\sharp}$ (see (\ref{eq:pY})) we have
\begin{equation*}
p_{\mathfrak{X}}^{-1}(C_{\bar{B}}^{\myp n})=C_{X\times
\bar{B}}^{\myp
n}=\{{\hatgamma}\in\varGamma_{\mathcal{Z}}^{\sharp}:\,{\hatgamma}(X\times
\bar{B})=n\}\in \mathcal{B}(\varGamma_{\mathcal{Z}}^{\sharp}),
\end{equation*}
since $X\times \bar{B}\in \mathcal{B}(\mathcal{Z})$. Furthermore,
the measurability of the projection
$\mathfrak{p}:\varGamma_{\mathfrak{X}}^{\sharp}\rightarrow
\varGamma_{X}^{\sharp}$ (see (\ref{eq:pr2})) was shown in
\cite[\S\myp3.3, p.~455]{BD3}. As a result, the composition of maps
in (\ref{eq:circ}) is measurable, as claimed.
\endproof

Let us define a measure on $\varGamma_{X}^{\sharp}$ as the
push-forward of $\hat{\mu}$ (see (\ref{eq:mu-hat})) under the map
$\mathfrak{q}$ defined in (\ref{proj1}), (\ref{eq:proj}):
\begin{equation}\label{eq:mu*}
\mathfrak{q}^{\ast}\hat{\mu}(A)\equiv
\hat{\mu}\myp(\mathfrak{q}^{-1}(A)),\qquad A\in
\mathcal{B}(\varGamma_{X}^{\sharp}),
\end{equation}
or, equivalently,
\begin{equation}
\int_{\varGamma_{X}^{\sharp}}F(\gamma)\,\mathfrak{q}^{\ast}
\hat{\mu}(\rd\gamma)=\int_{\varGamma_{\mathcal{Z}}^{\sharp}}
F(\mathfrak{q}({\hatgamma}))\,\hat{\mu}({\rd\hatgamma}),\qquad F\in
\mathrm{M}_{+}(\varGamma_{X}^{\sharp}). \label{eq:mu*1}
\end{equation}
The next general result shows that this measure may be identified
with the original cluster measure $\mucl$.

\begin{theorem}\label{th:mucl}
Measure \textup{(\ref{eq:mu*})} coincides with the cluster measure
$\mucl$,
\begin{equation}\label{eq:cl*}
\mucl=\mathfrak{q}^{\ast}\hat{\mu}\equiv \hat{\mu}\circ
\mathfrak{q}^{-1}.
\end{equation}
\end{theorem}

\proof Let us evaluate the Laplace transform of the measure
$\mathfrak{q}^{\ast}\hat{\mu}$. For any function
$f\in\mathrm{M}_{+}(X)$, we obtain, on using (\ref{eq:mu-hat}),
(\ref{eq:proj}) and(\ref{eq:mu*1}),
\begin{align*}
L_{\mathfrak{q}^{\ast}\hat{\mu}}(f)&
=\int_{\varGamma_{X}^{\sharp}}\exp (-\langle f,\gamma \rangle
)\,\mathfrak{q}^{\ast
}\hat{\mu}(\rd\gamma)\\
&=\int_{\varGamma_{\mathcal{Z}}^{\sharp}}
\exp(-\langle f,\mathfrak{q}({\hatgamma})\rangle)\,\hat{\mu}(\rd\hatgamma) \\
&
=\int_{\varGamma_{X}}\biggl(\int_{\varGamma_{\mathfrak{X}}^{\sharp}}
\exp\biggl(-\sum_{x\in \gamma}f(\mathfrak{p}(\bar{y}_{x}+x))\biggr)
{\textstyle\bigotimes\limits_{x\in\gamma}}
\,\eta (\rd\bar{y}_{x})\biggr)\,\mu(\rd\gamma)\\
& =\int_{\varGamma_{X}^{\sharp}}\biggl(\int_{\varGamma
_{\mathfrak{X}}^{\sharp}}\,\prod_{x\in \gamma}\exp
\bigl(-f(\bar{y}_{x})\bigr) {\textstyle\bigotimes\limits_{x\in
\gamma}}
\,\eta_{x}(\rd\bar{y}_{x})\biggr)\,\mu(\rd\gamma) \\
& =\int_{\varGamma_{X}^{\sharp}}
\prod_{x\in\gamma}\biggl(\int_{\mathfrak{X}}\exp \biggl(-\sum_{y\in
\bar{y}}f(y)\biggr)\,\eta_{x}(\rd\bar{y})\biggr) \,\mu(\rd\gamma),
\end{align*}
which coincides with the Laplace transform (\ref{laplace-G}) of the
cluster measure $\mucl$. \endproof

An alternative description of the cluster measure $\mucl$ can be
given as follows. Consider a natural map
$\mathfrak{r}_{\gamma}:\mathfrak{X}^{\gamma}\rightarrow
\varGamma_{\mathfrak{X}}^{\sharp}$ defined by
\begin{equation*}
\mathfrak{X}^{\gamma}\ni(\bar{y}_x)_{x\in\gamma}
\stackrel{\mathfrak{r}_\gamma}{\longmapsto}
\textstyle{\bigsqcup\limits_{x\in\gamma}}\{\bar{y}_x\}
\in\varGamma_{\mathfrak{X}}^{\sharp}
\end{equation*}
(see (\ref{eq:hatgamma})). The map $\mathfrak{r}_\gamma$ is
measurable, which can be shown by repeating the arguments used in
\cite[\S\myp3.3, p.~455]{BD3} in the proof of measurability of
$\mathfrak{p}$. Further, define the map (cf.\ (\ref{eq:pr2}))
\begin{equation}\label{eq:p-gamma}
p_{\gamma}:=\mathfrak{p}\circ \mathfrak{r}_{\gamma}:\
\mathfrak{X}^{\gamma} \stackrel{\mathfrak{r}_{\gamma}}
\longrightarrow\varGamma_{\mathfrak{X}}^{\sharp}
\stackrel{\mathfrak{p}}{\longrightarrow }\varGamma_{X}^{\sharp}.
\end{equation}
Clearly, $p_\gamma$ is measurable as a composition of measurable
maps. Note that the projection $p_{\frakX}$ defined in (\ref{eq:pY})
can be represented as
\begin{equation}\label{eq:p-gamma+}
p_{\frakX}(\hatgamma)=\mathfrak{r}_\gamma(\bar{y}^\gamma), \qquad
\hatgamma=(\gamma,\bar{y}^\gamma)\in\varGamma_{\mathcal{Z}}^\sharp.
\end{equation}
Furthermore, applying $\mathfrak{p}$ to both sides of equality
(\ref{eq:p-gamma+}) and using relations (\ref{eq:circ}) and
(\ref{eq:p-gamma}), we obtain the representation
\begin{equation}\label{eq:q-gamma}
\mathfrak{q}(\hatgamma)=p_\gamma(\bar{y}^\gamma), \qquad
\hatgamma=(\gamma,\bar{y}^\gamma)\in\varGamma_{\mathcal{Z}}^\sharp.
\end{equation}

Consider the probability measures $\varpi^{\gamma}$ and
$\mu^{\gamma}$ on $\varGamma_{\frakX^\gamma}^\sharp$ and
$\varGamma_{X}^{\sharp}$, respectively, defined by
\begin{align}
\label{eq:varpi-gamma} \varpi^{\gamma}:={}
&\mathfrak{r}_{\gamma}^{\ast}\eta^{\gamma},\\[.2pc]
\label{eq:mu-gamma} \mu^{\gamma}:={}&p_{\gamma}^{\ast}\eta^{\gamma}.
\end{align}

\begin{theorem}\label{altern}
The cluster measure $\mucl$ on $\varGamma_X^\sharp$ is represented
in either of the following two forms,
\begin{align}
\label{mu-cluster1}
\mucl(\rd\gamma)&=\int_{\varGamma_{X}^{\sharp}}\mu^{\zeta}(\rd\gamma)
\,\mu(\rd\zeta),\\
\label{mu-cluster2}
\mucl(\rd\gamma)&=\mathfrak{p}^{\ast}\varpi(\rd\gamma),
\end{align}
where $\varpi$ is a measure on $\varGamma_{\frakX}^\sharp$ defined
by
\begin{equation}\label{eq:varpi0}
\varpi(\rd\bar{\gamma}):=\int_{\varGamma_{X}^{\sharp}}
\varpi^{\gamma}(\rd\bar{\gamma})\,\mu(\rd\gamma).
\end{equation}
\end{theorem}

\proof By the change of measure (\ref{eq:mu-gamma}) and relations
(\ref{eq:mu-hat}) and (\ref{eq:q-gamma}), we have, for any Borel
function $F\in\mathrm{M}_+(\varGamma_X^\sharp)$,
\begin{align*}
\int_{\varGamma_{X}^{\sharp}}
\biggl(\int_{\varGamma_{X}^{\sharp}}F(\zeta)\,\mu^{\gamma}(\rd\zeta)
\biggr)\,\mu(\rd\gamma)
&=\int_{\varGamma_{X}^{\sharp}}\left(\int_{\mathfrak{X}}F(p_{\gamma}(\bar{y}^{\gamma}))\,
\eta^{\gamma}(\rd\bar{y}^{\gamma})\right)\mu(\rd\gamma) \\
&=\int_{\varGamma_{\mathcal{Z}}^{\sharp}}
F(\mathfrak{q}(\hatgamma))\,\hat{\mu}(\rd{{\hatgamma}})\\
&=\int_{\varGamma_{\mathcal{Z}}^{\sharp}}
F(\hatgamma)\,\mathfrak{q}^*\hat{\mu}(\rd\gamma),
\end{align*}
according to (\ref{eq:cl*}). Thus, formula (\ref{mu-cluster1}) is
proved.

Similarly,
using relation
(\ref{eq:varpi0}) and the change of measure (\ref{eq:varpi-gamma}),
we get
\begin{align*}
\int_{\varGamma_{\mathfrak{X}}^{\sharp}}
F(\zeta)\:\mathfrak{p}^{\ast}\varpi(\rd\zeta)
&=\int_{\varGamma_{X}^{\sharp}}
\biggl(\int_{\varGamma_{\mathfrak{X}}^{\sharp}}
F(\mathfrak{p}(\bar{\gamma}))\,\varpi^{\gamma}(\rd\bar{\gamma})\biggr)\,\mu(\rd\gamma)\\
&=\int_{\varGamma_{X}^{\sharp}}\biggl(\int_{\mathfrak{X}}
F(\mathfrak{p}\circ \mathfrak{r}_{\gamma}(\bar{y}^{\gamma}))\,
\eta^{\gamma}(\rd\bar{y}^{\gamma})\biggr)\,\mu(\rd\gamma) \\
&=\int_{\varGamma_{X}^{\sharp}}\biggl(\int_{\mathfrak{X}}
F(p_{\gamma}(\bar{y}^{\gamma}))\,\eta^{\gamma}(\rd\bar{y}^{\gamma})\biggr)\,\mu(\rd\gamma) \\
&=\int_{\varGamma_{\mathcal{Z}}^{\sharp}}
F(\mathfrak{q}(\hat\gamma))\hat{\mu}(\rd{{\hatgamma})},
\end{align*}
which proves formula (\ref{mu-cluster2}).
\endproof

\begin{remark}
In the case where $X=\RR^{d}$,
\,$\eta_{x}(\rd\bar{y})=\eta_0(\rd\bar{y}-x)$ and $\mu(\rd\gamma)$
is a Poisson measure $\pi_{\theta}(\rd\gamma)$ with intensity
$\theta$, the measure $\varpi$ coincides with the auxiliary Poisson
measure $\pi_{\sigma}$ considered in \cite{BD3}, with intensity
measure $\sigma(\bar{B})=\int_{X}\eta_{x}(\bar{B})\,\theta
(\rd{x})$, \,$\bar{B}\in\mathcal{B}(\frakX)$.
\end{remark}

The relationships between various measure spaces introduced above
are succinctly illustrated by the following commutative diagrams:

\begin{figure}[h]
\begin{center}
\mbox{$\begin{CD}
(\varGamma^\sharp_{\mathcal{Z}},\hat{\mu})
   @>\mbox{\footnotesize$\Id
   $}>>
     (\varGamma^\sharp_{X}\times \mathfrak{X}^{\gamma},\,\mu\otimes\eta^{\gamma})
     @.\hspace{3.5pc}
       (\varGamma^\sharp_{\mathcal{Z}},\hat{\mu})
       @>\mbox{\footnotesize$\,\Id\otimes \mathfrak{r}_{\gamma}$}>>
        (\varGamma^\sharp_{X}\times \varGamma^\sharp_{\mathfrak{X}},\,\mu\otimes\varpi^{\gamma})\\[.3pc]
@V\mbox{\footnotesize$\mathfrak{q}$}VV
   @VV\mbox{\footnotesize$\Id\otimes p_{\gamma}$}V
     \hspace{3.5pc}@V\mbox{\footnotesize$\mathfrak{q}$}VV
       @VV\mbox{\footnotesize$
       \int \rd\mu
       $}V\\[-.0pc]
(\varGamma^\sharp_{X},\mucl)
   @<\mbox{\footnotesize$\ \int \rd\mu$} <<
     (\varGamma^\sharp_{X}\times \varGamma^\sharp_{\frakX},\,\mu\otimes \mu^{\gamma})
       @.\hspace{3.5pc}(\varGamma^\sharp_{X},\mucl)
           @<\mbox{\footnotesize$\mathfrak{p}$}<
                \phantom{\mbox{\footnotesize$\,id\otimes\mathfrak{r}_{\gamma}$}}<
(\varGamma^\sharp_{\frakX},\varpi)
\end{CD}$}
\end{center}
\end{figure}

\subsection{Conditions for absence of accumulation and multiple
points}\label{sec:3.2}

Let us now give sufficient conditions for the cluster point process
to be proper, so that $\mucl(\varGamma_X)=1$. For any Borel subset
$B\subset X$, consider the set
\begin{align}\label{eq:frakX_B}
\mathfrak{X}_{B}&{}:=\{\bar {y}\in {\mathfrak{X}}:\
{\mathfrak{p}}(\bar {y})\cap B\neq \emptyset
\}\in\mathcal{B}(\mathfrak{X}),
\end{align}
where the projection map $\mathfrak{p}$ is defined in (\ref{eq:pr}).
That is to say, the set $\mathfrak{X}_{B}$ consists of all points
$\bar{y}\in \mathfrak{X}$ with at least one coordinate
$y_i\in\bar{y}$ belonging to $B$.

\begin{condition}\label{condition2}
For any compact set $B\subset X$,
\begin{equation}\label{eq:cond2}
\int_{\varGamma_{X}}\sum_{x\in \gamma
}\eta_{x}(\mathfrak{X}_{B})\,\mu (\rd\gamma)<\infty.
\end{equation}
\end{condition}

\begin{remark}
In view of definition (\ref{eq:frakX_B}), the left-hand side of
(\ref{eq:cond2}) equals the expected number (under the measure
$\mucl$) of clusters that contribute at least one point to the set
$B$.
\end{remark}

We introduce the set
\begin{equation*}
\tilde{\mathfrak{X}}:=\{\bar{y}\in \mathfrak{X}:~\forall\,
y_{i},y_{j}\in \bar{y},\ y_{i}\neq y_{j}\}.
\end{equation*}

\begin{condition}\label{condition3}
For $\mu$-a.a.\ configurations $\gamma\in\varGamma_{X}$, the
probability measure $\eta^\gamma$ on $\mathfrak{X}^{\gamma}$ (see
(\ref{eq:eta-gamma})) is concentrated on the set
\begin{equation}\label{eq:cond3}
\widetilde{\mathfrak{X}^{\gamma}}:=\{\bar{y}^{\gamma}\in
(\tilde{\mathfrak{X}})^{\gamma}:~\forall\, \{x_{1},x_{2}\}\subset
\gamma,\ \,\mathfrak{p}(\bar{y}_{x_{1}})\cap
\mathfrak{p}(\bar{y}_{x_{2}})=\emptyset\},
\end{equation}
that is, $\eta^\gamma(\widetilde{\mathfrak{X}^{\gamma}})=1$.
\end{condition}

\begin{remark}
The set $\widetilde{\mathfrak{X}^{\gamma}}$ ensures that different
clusters attached to the ground configuration $\gamma$ do not have
common points.
\end{remark}

\begin{remark}
A sufficient condition for (\ref{eq:cond3}) is that for any $x\in X$
the measure $\eta_{x}$ is absolutely continuous with respect to the
volume measure in $\mathfrak{X}$.
\end{remark}

From now on, we tacitly assume that the intra-cluster configurations
are a.s.-proper:
\begin{condition}\label{condition3.5}
For each $x\in X$ and for $\eta_x$-a.a.\ $\bar{y}\in\frakX$, the
projection set $\mathfrak{p}(\bar{y})\subset X$ is locally finite
and simple.
\end{condition}

\begin{theorem}\label{proper}
Let $\mucl$ be a cluster measure on the generalised configuration
space $\varGamma_X^\sharp$. Then
\begin{enumerate}
\item[\rm (a)]
under Condition \textup{\ref{condition2}}, \,$\mucl$-a.a\
configurations $\gamma\in\varGamma_X^\sharp$ are locally
finite\textup{;}
\item[\rm (b)]
under Condition \textup{\ref{condition3}}, \,$\mucl$-a.a\
configurations $\gamma\in\varGamma_X^\sharp$ are simple.
\end{enumerate}
Therefore, if\/ both Conditions \textup{\ref{condition2}} and
\textup{\ref{condition3}} are met then the cluster measure $\mucl$
is concentrated on the proper configuration space $\varGamma_{X}$.
\end{theorem}

\proof (a) Let $B\subset X$ be a compact set. From formula
(\ref{mu-cluster2}) and definition (\ref{eq:frakX_B}) of the set
$\mathfrak{X}_B$, it is clear that $\gamma(B)<\infty$ for
$\mucl$-a.a.\ configurations $\gamma\in\varGamma_X^\sharp$ if and
only if
\begin{equation}\label{finiteness1}
\bar\gamma(\mathfrak{X}_{B})<\infty\quad \text{for\ \ }
\varpi\text{-a.a.\ \ }
\bar\gamma\in\varGamma_{\mathfrak{X}}^{\sharp},
\end{equation}
where the measure $\varpi$ is defined in (\ref{eq:varpi0}). Let
$f(\bar{y}):=\mathbf{1}_{\mathfrak{X}_{B}}(\bar{y})$,
\,$\bar{y}\in\frakX$. Recalling definitions (\ref{eq:eta-gamma}),
(\ref{mu-cluster1}), (\ref{mu-cluster2}) and using Condition
\ref{condition2}, we obtain
\begin{align*}
\int_{\varGamma_{\mathfrak{X}}^{\sharp}}\langle
f,\bar\gamma\rangle\,\varpi(\rd\bar\gamma)
&=\int_{\varGamma_{X}}\left(\int_{\varGamma_{\mathfrak{X}}^{\sharp}}\langle
f,\bar\gamma\rangle\,\varpi^{\gamma}(\rd\bar\gamma)\right)\mu(\rd\gamma)\\
&=\int_{\varGamma_{X}}\left(\int_{\mathfrak{X}^{\gamma}}
\sum_{\bar{y}\in
\bar{y}^{\gamma}}f(\bar{y})\,\eta^{\gamma}(\rd\bar{y}^{\gamma})\right)
\mu(\rd\gamma)\\
&=\int_{\varGamma_{X}}\left(\sum_{x\in \gamma}
\int_{\mathfrak{X}_B}\eta_{x}(\rd\bar{y})\right)\mu(\rd\gamma) \\
&=\int_{\varGamma_{X}}\left(\sum_{x\in\gamma}
\eta_{x}(\mathfrak{X}_{B})\right)\mu(\rd\gamma)<\infty,
\end{align*}
which implies (\ref{finiteness1}). Thus, the absence of accumulation
points is proved.

\smallskip
(b) Let $\varGamma_{X}^{\ddagger}$ be the set of all generalised
configurations in $X$ that have multiple points. By definition
(\ref{eq:p-gamma}) of the map $p_{\gamma}$ there is the inclusion
\begin{equation*}
p_{\gamma}^{-1}(\varGamma_{X}^{\ddagger})\subset
\mathfrak{X}^{\gamma}\setminus\widetilde{\mathfrak{X}^{\gamma}},
\end{equation*}
where the set $\widetilde{\mathfrak{X}^{\gamma}}$ is introduced in
(\ref{eq:cond3}). Hence, from (\ref{eq:mu-gamma}) we get, for
$\mu$-a.a.\ $\gamma\in\varGamma_X$,
\begin{equation*}
\mu^{\gamma}(\varGamma_{X}^{\ddagger})
=\eta^{\gamma}(p_{\gamma}^{-1}(\varGamma_{X}^{\ddagger}))\le
1-\eta^{\gamma}(\widetilde{\mathfrak{X}^{\gamma}})=0
\end{equation*}
according to Condition \ref{condition3}, and by formula
(\ref{mu-cluster1}) this implies that
$\mucl(\varGamma_X^\ddagger)=0$.
\endproof

\subsection{Existence of moments}\label{sec:3.3}

\begin{definition}\label{def:M}
We shall say that the measure $\mu$ belongs to the class
$\mathcal{M}_{\theta}^{r}=\mathcal{M}_{\theta}^{r}(\varGamma_{X})$,
for some $r\ge1$ and a Borel measure $\theta$ on $X$, if for any
function $f\in\bigcap_{\myp 1\le \kappa\le r}L^{\kappa}(X,\theta)$
it holds that $\langle f,\gamma\rangle \in
L^{r}(\varGamma_{X},\mu)$, that is,
\begin{equation}\label{eq:Mr}
\int_{\varGamma_X}\left|\langle
f,\gamma\rangle\right|^r\mu(\rd\gamma)<\infty.
\end{equation}
\end{definition}

\begin{lemma}\label{rm:nested} The family of the classes
$\{\mathcal{M}_\theta^r,\ r\ge 1\}$ is nested, that is,
$\mathcal{M}_\theta^{r+\delta}\!\mynn\subset \mathcal{M}_\theta^{r}$
for all $r\ge1$ and any $\delta>0$.
\end{lemma}
\proof Indeed, let $\mu\in\mathcal{M}_\theta^{r+\delta}$, then for
any $f\in \bigcap_{\myp 1\le q\le r+\delta}L^{q}(X,\theta)\subset
\bigcap_{\myp 1\le q\le r}L^{q}(X,\theta)$ we have, by the Lyapunov
inequality,
$$
\int_{\varGamma_X}|\langle f,\gamma\rangle|^r\,\mu(\rd\gamma)\le
\left(\int_{\varGamma_X}|\langle
f,\gamma\rangle|^{r+\delta}\,\mu(\rd\gamma)\right)^{r/(r+\delta)}<\infty,
$$
hence $\mu\in\mathcal{M}_\theta^{r}$, as claimed.
\endproof
\begin{condition}\label{condition4}
There is a locally finite measure $\theta$ on $X$ (i.e.,
$\theta(B)<\infty$ for any compact $B\subset X$), referred to as the
\textit{reference measure}, such that $\mu \in
\mathcal{M}_{\theta}^{1}(\varGamma_{X})$.
\end{condition}

\begin{remark}
The condition $\mu\in \mathcal{M}_{\theta}^{1}(\varGamma_{X})$
implies that $\gamma(B)<\infty$ ($\mu$-a.s.) for any Borel set $B$
such that $\theta (B)<\infty$. Indeed, choose $f(x)={\mathbf
1}_B(x)\in L^1(X,\theta)$, then $\langle f,\gamma\rangle=\gamma(B)$
and Definition \ref{def:M} yields $\int_{\varGamma_X}
\!\gamma(B)\,\mu(\rd\gamma)<\infty$, hence the integrand is
$\mu$-a.s.\ finite.
\end{remark}

\begin{example}
Condition \ref{condition4} holds for a Poisson measure
$\pi_{\theta}$ with a locally finite intensity $\theta$, as well as
for a wide class of Gibbs perturbations of $\pi_{\theta}$ (see,
e.g., \cite{AKR1,AKR2}). More generally, any measure $\mu$ with
bounded correlation functions up to order $n$ (with respect to
$\theta$) belongs to $\mathcal{M}_\theta^n(\varGamma_X)$ (see
Appendix \ref{ap:2}).
\end{example}

\begin{example}
Example of a different type is given by $\mu=\delta_{\gamma_{0}}$,
the Dirac measure on $\varGamma_{X}$ concentrated on a given
configuration $\gamma_{0}\in\varGamma_X$ (e.g., if $X=\RR^d$ then we
can set $\gamma_{0}=\ZZ^{d}$). Here we have
\begin{equation*} \int_{\varGamma_{X}} \vert\langle f,\gamma\rangle
\vert^{n} \,\mu (\rd\gamma)=|\langle f,\gamma_{0}\rangle|^{n},
\end{equation*}
which implies that $\mu \in \mathcal{M}_{\theta}^{n}(\varGamma_{X})$
with $\theta = \sum_{x\in\gamma_{0}}\delta_{x}$.
\end{example}

\begin{definition}
Introduce the measures $\hat{\sigma}$ on $\mathcal{Z}=X\times
\frakX$ and $\bar{\sigma}$ on $\mathfrak{X}$ as follows
\begin{gather}
\label{eq:sigmaZ} \hat{\sigma}(\rd{x}\times\rd\bar{y}) :=
\eta_{x}(\rd\bar{y})\,\theta(\rd{x}),\\
\label{intensity} \bar{\sigma}(\rd\bar{y}) := \int_{X}
\eta_{x}(\rd\bar{y})\,\theta(\rd{x}).
\end{gather}
\end{definition}

\begin{lemma}\label{moments}
Let $\mu \in \mathcal{M}_{\theta}^{n}(\varGamma_{X})$ for some
integer $n\ge 1$. Then $\hat{\mu}\in
\mathcal{M}_{\hat{\sigma}}^{n}(\varGamma_{\mathcal{Z}})$.
\end{lemma}

\proof We have to show that, for any $f\in \bigcap_{\,1\le \kappa\le
n} L^\kappa(\mathcal{Z},\hat{\sigma})$,
\begin{equation}\label{eq:nth}
\int_{\varGamma_{\mathcal{Z}}}|\langle f,\hatgamma\rangle
|^{n}\,\hat{\mu}(\rd\hatgamma)<\infty.
\end{equation}
To this end, let us first observe using definition (\ref{eq:sigmaZ})
that
\begin{equation*}
\int_{X}\left(\int_{\mathfrak{X}}
|f(x,\bar{y})|^\kappa\,\eta_{x}(\rd\bar{y})\right)
\theta(\rd{x})=\int_{\mathcal{Z}}|f(z)|^\kappa
\,\hat{\sigma}(\rd{z})<\infty,
\end{equation*}
which means that the function
\begin{equation}\label{eq:f_m}
\bar{f}_\kappa(x):=\int_{\mathfrak{X}}
|f(x,\bar{y})|^\kappa\,\eta_{x}(\rd\bar{y}),\qquad x\in X,
\end{equation}
belongs to $L^1(X,\theta)$.
Moreover, by the Lyapunov inequality we have, for $q\ge1$,
\begin{align*} \int_X \bar{f}_\kappa(x)^q\,\theta(\rd
x)&=\int_{X}\left(\int_{\mathfrak{X}}
|f(x,\bar{y})|^\kappa\,\eta_{x}(\rd\bar{y})\right)^q \theta(\rd{x})\\
&\le \int_{X}\left(\int_{\mathfrak{X}}
|f(x,\bar{y})|^{\kappa q}\,\eta_{x}(\rd\bar{y})\right) \theta(\rd{x})\\
&=\int_{\mathcal{Z}}|f(z)|^{\kappa q}\,\hat{\sigma}(\rd{z})<\infty,
\end{align*}
as long as $1\le \kappa q\le n$. In other words,
\begin{equation}\label{eq:barf}
\bar{f}_\kappa(x)\in L^{q}(X,\theta), \qquad 1\le q\le n/\kappa.
\end{equation}

Now, using the multinomial expansion we can write for an integer
$n\ge1$
\begin{align}
\notag \int_{\varGamma_{\mathcal{Z}}}|\langle
f,\hatgamma\rangle|^{n}\,\hat{\mu}(\rd\hatgamma)& \leq
\int_{\varGamma_{\mathcal{Z}}}\left( \sum_{z\in\hatgamma}
|f(z)|\right)^{n}\,\hat{\mu}(\rd{\hatgamma})  \\
\label{eq:multi}& =\sum_{m=1}^{n}\int_{\varGamma_{\mathcal{Z}}}
\sum_{\{z_{1},\dots,\myp z_{m}\}\subset
{\hatgamma}}\phi_{n}(z_{1},\dots,z_{m}) \,\hat{\mu}(\rd\hatgamma),
\end{align}
where $\phi_{n}(z_{1},\dots ,z_{m})$ is a symmetric function given
by
\begin{equation}\label{eq:varPsi}
\phi_{n}(z_{1},\dots ,z_{m}):=\sum_{\substack{ i_{1},\dots ,\myp
i_{m}\ge 1 \\[.05pc] i_{1}+\dots +i_{m}=\myp n}}\frac{n!}{i_{1}!\cdots i_{m}!}
\,|f(z_{1})|^{i_{1}}\mynn\cdots |f(z_{m})|^{i_{m}}.
\end{equation}
By definition of the measure $\hat\mu$ (see (\ref{eq:mu-hat}) and
(\ref{eq:mu-hat-int})), the integral on the right-hand side of
(\ref{eq:multi}) is reduced to
\begin{multline}
\label{eq:k=} \int_{\varGamma_{X}}\sum_{\{x_{1},\dots,\myp
x_{m}\}\subset \gamma}\!\left(\int_{\mathfrak{X}^{\gamma}}
\phi_{n}(x_{1},\bar{y}_1;
\dots;x_{m},\bar{y}_m)\,\eta^\gamma(\rd\bar\gamma)\right)\mu(\rd\gamma)\\
=\sum_{\substack{ i_{1},\dots,\myp i_{m}\geq 1 \\[.05pc]
i_{1}+\dots + i_{m}=\myp n}}\!\frac{n!}{i_{1}!\cdots i_{m}!}
\int_{\varGamma_{X}}\sum_{\{x_{1},\dots,\myp x_{m}\}\subset \gamma}
\prod_{j=1}^{m} \bar{f}_{i_j}(x_j)\,\mu(\rd\gamma),
\end{multline}
where we used notation (\ref{eq:f_m}). Furthermore, with the help of
the Jensen inequality the integral on the right-hand side of
(\ref{eq:k=}) may be estimated from above by
\begin{align}
\notag \int_{\varGamma_{X}}\prod_{j=1}^{m}\sum_{x_{j}\in \gamma}
\bar{f}_{i_j}(x_j)\,\mu(\rd\gamma)&=\int_{\varGamma_{X}}\prod_{j=1}^{m}\langle
\bar{f}_{i_j},\gamma\rangle\,\mu(\rd\gamma)\\
\label{eq:ith} &\le\prod_{j=1}^{m}\left(\int_{\varGamma_{X}}\langle
\bar{f}_{i_j},\gamma\rangle^{n/i_j}\,\mu(\rd\gamma)\right)^{i_j/n}.
\end{align}

To summarise, by inspection of relations (\ref{eq:multi}),
(\ref{eq:varPsi}), (\ref{eq:k=}) and (\ref{eq:ith}) we see that in
order to verify (\ref{eq:nth}) it suffices to check that, for any
$k=1,\dots,n$,
\begin{equation}\label{eq:kth}
\int_{\varGamma_{X}}\langle
\bar{f}_{k},\gamma\rangle^{n/k}\,\mu(\rd\gamma) <\infty.
\end{equation}
But we already know (see (\ref{eq:barf})) that $\bar{f}_k\in
L^q(X,\theta)$ for $1\le q\le n/k$. On the other hand, by the
hypothesis of the lemma we have
$\mu\in\mathcal{M}_\theta^n(\varGamma_X)\subset
\mathcal{M}_\theta^{n/k}(\varGamma_X)$ (see Remark \ref{rm:nested}),
and now the required bound (\ref{eq:kth}) follows by condition
(\ref{eq:Mr}) with $r=n/k$.
\endproof

The next condition on the measure $\bar{\sigma}$ will play an
important part in our analysis.

\begin{condition}\label{condition5}
For any compact set $B\subset X$, it holds that
\begin{equation}\label{eq:sigma<infty}
\bar{\sigma}(\mathfrak{X}_{B})<\infty,
\end{equation}
where the set $\frakX_B$ is defined in (\ref{eq:frakX_B}).
\end{condition}

\begin{remark}
Conditions  \ref{condition4} and \ref{condition5} taken together
imply Condition \ref{condition2}. Indeed, from (\ref{intensity}) and
(\ref{eq:sigma<infty}) it follows that $\eta_x(\frakX_B)$ as a
function of $x\in X$ belongs to the space $L^1(X,\theta)$. Thus,
according to Definition \ref{def:M} we can apply Condition
\ref{condition4} to obtain
\begin{equation*}
\sum_{x\in \gamma}\eta_{x}(\mathfrak{X}_{B})\in
L^1(\varGamma_X,\mu),
\end{equation*}
which is nothing else but condition~(\ref{eq:cond2}). Hence, by
Theorem \ref{proper}(a), the measure $\mucl$ is concentrated on
configurations without accumulation points.
\end{remark}

Denote by $N_{X}(\bar{y})$ the ``dimension'' of vector
$\bar{y}\in\frakX$, that is, the total number of its components
\begin{equation}\label{eq:N(y)}
N_{X}(\bar{y}):= \sum_{n=0}^\infty
n\,\mathbf{1}_{X^n}(\bar{y}),\qquad
\bar{y}\in\frakX=\textstyle{\bigsqcup\limits_{n=0}^{\infty}} X^n.
\end{equation}

\begin{lemma}\label{lm:M^n}
Suppose that, in addition to Conditions \textup{\ref{condition4}}
and \textup{\ref{condition5}}, the function $N_{X}(\bar{y})$
satisfies, for any compact set $B\subset X$, an integrability
condition
\begin{equation}\label{eq:N<}
\int_{X}\int_{\mathfrak{X}_{B}}N_{X}(\bar{y})^{n}\,
\eta_{x}(\rd\bar{y})\,\theta (\rd{x})<\infty,
\end{equation}
Then the cluster measure\/ $\mucl$ belongs to the class
$\mathcal{M}^{n}_{\theta}(\varGamma_{X})$.
\end{lemma}

\proof Using the change of measure (\ref{eq:cl*}), for any $\phi \in
C_{0}(X)$ we obtain
\begin{equation*}
\int_{\varGamma_{X}}|\langle \phi ,\gamma \rangle |^{n}\,\mucl
(\rd\gamma)=\int_{\varGamma_{\mathcal{Z}}}|\langle
\phi,\mathfrak{q}(\hatgamma)\rangle|^{n}
\,\hat{\mu}(\rd\hatgamma)=\int_{\varGamma_{\mathcal{Z}}}|\langle
\mathfrak{q}^{\ast}\phi,\hatgamma\rangle|^{n}\,\hat{\mu}(\rd\hatgamma),
\end{equation*}
where
\begin{equation}\label{eq:q*}
\mathfrak{q}^{\ast}\phi (x,\bar{y}):=\sum_{y_{i}\in \bar{y}}\phi
(y_{i}),\qquad (x,\bar{y})\in \mathcal{Z}.
\end{equation}
It suffices to show that $\mathfrak{q}^{\ast}\phi \in
L^{m}(\mathcal{Z},\hat{\sigma})$ for any $m=1,\dots,n$. By the
elementary inequality $(a_{1}+\cdots +a_{k})^{m}\leq
k^{m-1}(a_{1}^{m}+\cdots +a_{k}^{m})$, from (\ref{eq:q*}) we have
\begin{equation}\label{eq:q*1}
\int_{\mathcal{Z}}|\mathfrak{q}^{\ast}\phi
(z)|^{m}\,\hat{\sigma}(\rd{z})\leq
\int_{\mathcal{Z}}N_{X}(\bar{y})^{m-1}\sum_{y_{i}\in \bar{y}}|\phi
(y_{i})|^{m}\,\hat{\sigma} (\rd{x}\times \rd\bar{y}).
\end{equation}
Recalling that $\hat{\sigma} (\rd{x}\times
\rd\bar{y})=\eta_{x}(\rd\bar{y})\,\theta(\rd{x}) $ and denoting
$C_{\phi}:=\sup_{x\in X}|\phi (x)|<\infty $ and $K_{\phi}:=\supp\phi
\subset X$, the right-hand side of (\ref{eq:q*1}) is dominated by
\begin{equation*}
(C_{\phi})^{m}\int_{X}\int_{\mathfrak{X}_{K_{\phi}}}
N_{X}(\bar{y})^{m}\,\eta_{x}(\rd\bar{y})\,\theta(\rd{x}),
\end{equation*}
which proves the result.
\endproof

\subsection{``Translations'' and the droplet cluster}\label{sec:droplet}

Let us describe a general setting that may be used to construct the
family of measures $\{\eta_{x}(\rd\bar{y})\}_{x\in X}$ on $\frakX$
via suitable push-forwards (``translations'') of a pattern measure
$Q$ defined on some auxiliary space. Examples of application of such
an approach will be given in Section \ref{sec:5} below.

More precisely, let $(W,\mypp\calB(W))$ be a measurable space, with
a Borel $\sigma$-algebra $\calB(W)$ generated by the open subsets of
$W$. Consider the corresponding space (cf.\ (\ref{eq:calX}))
\begin{equation}\label{eq:calW}
\mathfrak{W}:={\textstyle\bigsqcup\limits_{n=0}^{\infty}} W^{n},
\end{equation}
and let $Q$ be a probability measure on $\calB(\mathfrak{W})$. For
any map $\varphi: W \to X$, define as usual its diagonal lifting
$\bar{\varphi}:\mathfrak{W}\to \mathfrak{X}$ by
\begin{equation}\label{eq:lift}
\mathfrak{W}\ni \bar{w}\mapsto
\bar{\varphi}(\bar{w}):=(\varphi(w_{i}))_{w_{i}\in \bar{w}}\in
\mathfrak{X}.
\end{equation}
Like in Condition \ref{condition4}, it is assumed that the reference
measure $\theta$ on $X$ is locally finite.

The main assumption in this section is as follows.
\begin{condition}\label{condition7}
Suppose there is a measurable map
$$
W \times X\ni (w,x)\mapsto \varphi_x(w)\in X
$$
such that the measures $\eta_x$ on $\frakX$ are representable as
$\eta_x=\bar{\varphi}_x^{\ast}\myp Q$; that is, for all $x\in X$,
\begin{equation}\label{eq:x0->x}
\eta_x(\bar{B})=Q(\bar{\varphi}_x^{-1}(\bar{B})),\qquad
\bar{B}\in\calB(\frakX).
\end{equation}
\end{condition}

\begin{remark}
In view of formula (\ref{eq:x0->x}), we shall often consider
$\{\varphi_x(\cdot)\}_{x\in X}$ as a \textit{family} of the
maps \,$W\ni w\mapsto \varphi_x(w)\in X$ (indexed by $x\in X$).
\end{remark}

\begin{remark}\label{rm:phix}
Fubini's theorem implies that, for each $x\in X$, the map
$\mathfrak{W}\ni \bar{w}\mapsto\bar{\varphi}_x(\bar{w})\in\frakX$ is
measurable and hence $\bar{\varphi}_x^{-1}(\bar{B})$ is a Borel
subset of $\mathfrak{W}$, so that the right-hand side of formula
(\ref{eq:x0->x}) is well defined. We also have that, for any fixed
Borel set $\bar{B}\subset\frakX$, the function
$\eta_x(\bar{B}):X\to[0,1]$ is measurable.
\end{remark}

\begin{definition}
Given a
map $\varphi_x(w)$ as above, the set
\begin{equation}\label{eq:D}
D_{\myn B}(w):= \{x\in X: \,\varphi_x(w)\in B\}\subset X,\qquad w\in
W,\ \ B\in\mathcal{B}(X),
\end{equation}
is called a \textit{droplet} of shape $B$ anchored at $w$.
Furthermore, the set
\begin{equation}\label{eq:barD}
\bar{D}_{\myn B}(\bar{w}):= \bigl\{x\in X:
\,\bar{\varphi}_x(\bar{w})\in \frakX_B\bigr\}\subset X,\qquad
\bar{w}\in\mathfrak{W},
\end{equation}
is referred to as the \textit{droplet cluster} (of shape $B$)
anchored at $\bar{w}$.
\end{definition}

Note that the droplet $D_{\myn B}(w)$ is a Borel subset of $X$ for
each $w\in W$; moreover, by Remark \ref{rm:phix} the same is true
for the droplet cluster $\bar{D}_{\myn B}(\bar{w})$. On account of
definition (\ref{eq:frakX_B}), formula (\ref{eq:barD}) can be
rewritten in the form
\begin{equation}\label{eq:D-union}
\bar{D}_{\myn B}(\bar{w})
={\textstyle\bigcup\limits_{w_i\in\bar{w}}} D_{\myn B}(w_i), \qquad
\bar{w}\in\mathfrak{W}.
\end{equation}


The following identity enlightens the geometric meaning of Condition
\ref{condition5} stated above.

\begin{lemma}\label{lm:sigma=int}
For any Borel set $B\subset X$, there is the equality
\begin{equation}\label{eq:sigma=int}
\bar{\sigma}(\frakX_B)=\int_{\mathfrak{W}} \theta(\bar{D}_{\myn
B}(\bar{w}))\,Q(\rd\bar{w}).
\end{equation}
In particular, $\bar{\sigma}(\frakX_B)<\infty$ if and only if the
right-hand side of equation \textup{(\ref{eq:sigma=int})} is finite.
\end{lemma}
\proof According to (\ref{intensity}), (\ref{eq:x0->x}) and
(\ref{eq:D}), we have
\begin{align*}
\bar{\sigma}(\frakX_B)
&=\int_{X}\int_{\frakX_B}\eta_{x}(\rd\bar{y})\,\theta (\rd{x})\\
&=\int_{X}\left(\int_{\mathfrak{W}}{\mathbf{1}}_{\frakX_B}(\bar{\varphi}_x(\bar{w}))
\,Q(\rd\bar{w})\right)\theta (\rd{x})\\
&=\int_{\mathfrak{W}}\left(\int_{X}{\mathbf{1}}_{\bar{D}_{\myn
B}(\bar{w})}
(x)\,\theta (\rd{x})\right)Q(\rd\bar{w})\\
&=\int_{\mathfrak{W}}\theta(\bar{D}_{\myn
B}(\bar{w}))\,Q(\rd\bar{w}),
\end{align*}
as claimed.
\endproof

Due to formula (\ref{eq:sigma=int}), Condition \ref{condition5} can
be rewritten as follows.
\begin{condit} 
For any compact set $B\in{\mathcal{B}}(X)$, the mean
$\theta$-measure of the droplet cluster\/ $\bar{D}_{\myn
B}(\bar{w})\subset X$ is finite,
\begin{equation}\label{eq:condA2}
\int_{\mathfrak{W}} \theta(\bar{D}_{\myn B}(\bar{w}))\,
Q(\rd\bar{w})<\infty\myp.
\end{equation}
\end{condit}

Building on Lemma \ref{lm:sigma=int}, let us give two simple
criteria, either of which is sufficient for Condition
\ref{condition5}\mypp$'$ and hence for Condition \ref{condition5}.
The first criterion below (Proposition \ref{pr:3.9}) bounds the
growth of the droplet volume and also assumes a finite mean number
of points in the cluster, while the second criterion (Proposition
\ref{pr:3.10}) requires the continuity and separability of the maps
$\varphi_x(w)$ and puts a restriction on the range of the parent
cluster.

\begin{proposition}\label{pr:3.9}
Suppose that the following two conditions hold\textup{:}
\begin{enumerate}
\item[\rm (i)] \emph{(finite range of ``translations'')}
\,for any compact set $B\subset X$, the map $\varphi_x(w)$ and the
measure $\theta$ satisfy the bound
\begin{equation}\label{eq:C_B}
C_{B}:=\sup_{w\in W}\theta(D_{\myn B}(w))<\infty;
\end{equation}
\item[\rm (ii)]
\emph{(finite mean of the cluster size)} \,the total number of
components in\/ $\bar{w}\in\mathfrak{W}$ \textup{(}cf.\
\textup{(\ref{eq:N(y)})}\textup{)} satisfies the integrability
condition \textup{(}cf.\ \textup{(\ref{eq:N<})}\textup{)}
\begin{equation}\label{eq:EN<}
\int_{\mathfrak{W}}N_{W}(\bar{w})\, Q(\rd\bar{w})<\infty.
\end{equation}
\end{enumerate}
Then Condition~\textup{\ref{condition5}\mypp$'$} is satisfied.
\end{proposition}
\proof By Lemma \ref{lm:sigma=int}, it suffices to show that the
integral on the right-hand side of (\ref{eq:sigma=int}) is finite.
From (\ref{eq:D-union}) and (\ref{eq:C_B})  we readily obtain
\begin{align}
\label{eq:theta<Nxtheta} \theta(\bar{D}_{\myn B}(\bar{w}))
&\le\sum_{w_i\in\bar{w}} \theta(D_{\myn B}(w_i))\le C_B\mypp
N_{W}(\bar{w}),\qquad \bar{w}\in\mathfrak{W},
\end{align}
and by condition (\ref{eq:EN<}) it follows
$$
\int_{\mathfrak{W}}\theta(\bar{D}_{\myn B}(\bar{w})) \,Q(\rd\bar{w})
\le C_B \int_{\mathfrak{W}} N_{W}(\bar{w})\, Q(\rd\bar{w})<\infty,
$$
as required.
\endproof

\begin{remark}\label{rm:3.8}
Bound (\ref{eq:C_B}) holds, for example, if for every $w\in W$ the
map $X\ni x\mapsto \varphi_{x}(w)\in X$ is an isometry and the
measure $\theta$ is absolutely continuous with respect to the volume
measure on $X$, with a bounded Radon--Nikodym density.
\end{remark}

\begin{proposition}\label{pr:3.10}
Suppose that the family of measurable maps $\varphi_x(w)$ described
above satisfies in addition the following two conditions\textup{:}
\begin{enumerate}
\item[\rm(i)] \emph{(continuity in $x$)} \ the map
$\varphi_x(w)$ is continuous in $x\in X$\textup{;} that is, for any
open subset $U\subset X$ and each $w\in W$, the set $\{x\in X:
\varphi_x(w)\in U\}$ is open in $X$\textup{;}

\item[\rm (ii)] \emph{(separability)} \
for any compact set $B\subset X$ and each  $w\in W$, there exists a
compact $B_{w}\subset X$ such that for any $x\notin B_{w}$ we have
$\varphi_{x}(w)\notin B$.
\end{enumerate}
Assume also that there is a compact set $E_0\in\mathcal{B}(W)$ such
that $Q(\mathfrak{E}_0)=1$, where
$\mathfrak{E}_0:=\bigsqcup_{n=0}^\infty E_0^n$ \,\textup{(}cf.\
\textup{(\ref{eq:calW})}\textup{)}\textup{;} that is, all components
of\/ $Q$-a.a.\ vectors $\bar{w}\in\mathfrak{W}$ lie in $E_0\subset
W$. Then \strut{}Condition~\textup{\ref{condition5}}\mypp$'$ is
satisfied.
\end{proposition}
\proof Let $B\subset X$ be an arbitrary compact set.
%
Using formula (\ref{eq:D-union}) and definition (\ref{eq:D}), for
any $\bar{w}\in\mathfrak{E}_0$ we have the inclusion
\begin{align}
\notag \bar{D}_{\myn
B}(\bar{w})={\textstyle\bigcup\limits_{w_i\in\bar{w}}} D_{\myn
B}(w_i) &\subset {\textstyle\bigcup\limits_{w\in E_0}} D_{\myn B}(w)
\\
\label{eq:DB} &\equiv \{x: \varphi_x(E_0)\cap B\ne
\emptyset\}=:\tilde{D}_{\myn B}\subset X.
\end{align}
To complete the proof, it suffices to show that
$\theta(\tilde{D}_{\myn B})<\infty$, since then, by Lemma
\ref{lm:sigma=int}, it will follow
\begin{align*}
\bar{\sigma}(\frakX_B)=\int_{\mathfrak{W}} \theta(\bar{D}_{\myn
B}(\bar{w})) \,Q(\rd\bar{w})&=\int_{\mathfrak{E}_0}
\theta(\bar{D}_{\myn B}(\bar{w}))
\,Q(\rd\bar{w})\\
&\le\theta(\tilde{D}_{\myn B})\int_{\mathfrak{E}_0}Q(\rd\bar{w})=
\theta(\tilde{D}_{\myn B})<\infty.
\end{align*}
To this end, for each $w\in E_0$ consider the set
\begin{equation}\label{eq:T}
A_w :=\{w '\in E_0: \,\varphi_x(w ')\notin B\ \text{\,for all
\,}x\notin B_w \}\subset E_0,
\end{equation}
where $B_w\subset X$ is a compact defined in property (ii); in
particular, it follows that $w \in A_w $ and hence
$E_0\subset\bigcup_{w\in W} A_w$. Furthermore, using property (i) it
can be shown that each $A_w$ is an open set in the topology induced
from $W$ by restriction to $E_0$ (i.e., with open sets in $E_0$
defined as $U\cap E_0$ for all $U$ open in $W$).
Since $E_0$ is compact,
there is a finite subcover; that is, one can choose finitely many
points $w_1,\dots,w_m\in E_0$ such that $E_0\subset \bigcup_{\myp
i=1}^{\myp m} A_{w_i}$. Then, using (\ref{eq:T}), it is easy to see
that for any $x$ outside the set $B_*:=\bigcup_{\myp i=1}^{\myp m}
B_{w_i}$ we have $\varphi_x(w )\notin B$ for all $w \in E_0$.
According to definition (\ref{eq:DB}) of the set $\tilde{D}_{B}$,
this implies that $\tilde{D}_{\myn B}\subset B_*$, hence
$\theta(\tilde{D}_{\myn B})\le \theta(B_*)<\infty$. The proof is
complete.
\endproof

The next statement gives a criterion sufficient for Condition
\ref{condition3} (in turn, implying the simplicity of the cluster
measure $\mucl$, according to Theorem \ref{proper}).

\begin{proposition}\label{pr:simple}
Let $\mucl$ be a cluster measure on the generalised configuration
space $\varGamma_ X^\sharp$. Assume that the background measure
$\mu$ of cluster centres has a locally bounded second-order
correlation function $\kappa_{\myn\mu}^2$
%
%
\textup{(}see
Appendix \textup{\ref{ap:2}}\textup{)}. Assume also that for
$\theta$-a.a.\ $x\in X$ the corresponding ``point'' droplet cluster
$\bar{D}_{\{x\}}(\bar{w})$ has a.s.\ zero $\theta$-measure,
\begin{equation}\label{eq:condB2}
\theta\bigl(\bar{D}_{\{x\}}(\bar{w})\bigr)=0\quad \text{for} \
\,Q\text{-a.a.}\ \bar{w}\in\mathfrak{W}.
\end{equation}
Then Condition \textup{\ref{condition3}} is satisfied and hence
$\mucl$-a.a.\ configurations $\gamma\in \varGamma_{X}^\sharp$ are
simple.
\end{proposition}
\proof It suffices to prove that, for any compact set
$\varLambda\subset X$,
there are $\mucl$-a.s.\ no cross-ties between the clusters whose
centres belong to $\varLambda$. In view of the projection
construction of the cluster measure $\mucl$ (see (\ref{eq:cl*})),
this means that if $A_\varLambda$ is the set of generalised
configurations $\hatgamma\in\varGamma_{\mathcal{Z}}^\sharp$, each
with at least two points $z=(x,\bar{y}_{x})$,
\,$z'=(x',\bar{y}_{x'})$ \,($z,z'\in\hatgamma$, \,$z\ne z'$) such
that $\{x,x'\}\subset \gamma\cap\varLambda$ and
$\mathfrak{p}(\bar{y}_x)\cap\mathfrak{p}(\bar{y}_{x'})\ne\emptyset$,
then we must show that $\hat{\mu}(A_\varLambda)=0$. Note that since
the ground configuration $\gamma\in\varGamma_X^\sharp$ may have
multiple points, the points $x=p_x(z)$, $x'=p_X(z')$ in the pair
$\{x,x'\}\subset\gamma$ are allowed to coincide.

Recalling the skew-product definition (\ref{eq:mu-hat}) of
$\hat{\mu}$,
we see by inspection of all pairs $\{x,x'\}\subset
\gamma:=p_X(\hatgamma)$ that
\begin{align}\label{eq:Psi3}
\hat{\mu}(A_\varLambda)&\le\int_{\varGamma_{X}^\sharp} \mynn
\sum_{\{x,\myp x'\}\subset\gamma} \mathbf{1}_{\varLambda^2}(x,x')\,
f(x,x')\,\mu(\rd\gamma),
\end{align}
where
\begin{align}
f(x,x'):= \eta_{x}\otimes\eta_{x'}(\mathcal{D}_{\mathfrak{X}})=
\int_{\frakX^2}\mathbf{1}_{\mathcal{D}_{\mathfrak{X}}}(\bar{y},\myp\bar{y}')\,
\eta_{x}(\rd\bar{y})\,\eta_{x'}(\rd\bar{y}') \label{eq:F}
\end{align}
and the set $\mathcal{D}_{\mathfrak{X}}\subset \mathfrak{X}^2$ is
defined by
\begin{equation}\label{eq:DX}
\mathcal{D}_{\mathfrak{X}}:= \{(\bar{y},\bar{y}')\in \mathfrak{X}^2:
\,\mathfrak{p}(\bar{y})\cap \mathfrak{p}(\bar{y}')\ne\emptyset\}.
\end{equation}
By definition (\ref{corr-funct}) of correlation functions, the
right-hand side of (\ref{eq:Psi3}) is reduced to
\begin{equation}\label{eq:=0}
\frac{1}{2!}\int_{\varLambda^2} f(x,x')\mypp \kappa_\mu^2(x,x')\,
\theta(\rd x)\,\theta(\rd x')\le \const \int_{\varLambda^2}
f(x,x')\, \theta(\rd x)\,\theta(\rd x'),
\end{equation}
since, by assumption, $\kappa_\mu^2$ is
bounded on $\varLambda^2$. Furthermore, substituting (\ref{eq:F})
and changing the variables $\bar{y}=\varphi_x(\bar{w})$,
$\bar{y}'=\varphi_{x'}(\bar{w}')$ (see (\ref{eq:x0->x})), the
integral on the right-hand side of (\ref{eq:=0}) is rewritten as
\begin{align*}
\int_{\varLambda^2}&
\left(\int_{\frakX^2}\mathbf{1}_{\mathcal{D}_{\mathfrak{X}}}
(\varphi_x(\bar{w}),\myp\varphi_{x'}(\bar{w})')\,
Q^{\otimes\myp2}(\rd\bar{w}\times\rd\bar{w}')\right)\theta^{\otimes\myp2}(\rd
x\times \rd x')\\
&\qquad=\int_{\frakX^2}
\left(\int_{\varLambda^2}\mathbf{1}_{\mathcal{D}_{\mathfrak{X}}}
(\varphi_x(\bar{w}),\myp\varphi_{x'}(\bar{w})')\,
\theta^{\otimes\myp2}(\rd
x\times \rd
x')\right)Q^{\otimes\myp2}(\rd\bar{w}\times\rd\bar{w}')\\[.3pc]
&\qquad\qquad =\int_{\frakX^2}
\theta^{\otimes\myp2}\bigl(B_\varLambda(\bar{w},\bar{w}')\bigr)
\,Q^{\otimes\myp2}(\rd\bar{w}\times\rd\bar{w}'),
\end{align*}
where the set $B_\varLambda(\bar{w},\bar{w}')\subset \varLambda^2$
is given by (cf.\ (\ref{eq:DX}))
$$
B_\varLambda(\bar{w},\bar{w}'):= \{(x,x')\in\varLambda^2:
\,\varphi_x(w)=\varphi_{x'}(w')\ \text{\,for some }\, w\in\bar{w},\
w'\in\bar{w}'\}.
$$
It remains to note that
\begin{align*}
\theta^{\otimes\myp2}\bigl(B_\varLambda(\bar{w},\bar{w}')\bigr)
&=\int_\varLambda
\theta\Bigl(\textstyle\bigcup\nolimits_{w_i\in\bar{w}}
\!\bigcup\nolimits_{w_j'\in\bar{w}'}
\{x':\varphi_{x'}(w_j')=\varphi_{x}(w_i)\}\Bigr)\,\theta(\rd x)\\[.2pc]
&\le \sum_{w_i\in\bar{w}}\int_\varLambda \theta
\Bigl(\textstyle\bigcup\nolimits_{w_j'\in\bar{w}'}\mynn
\{x':\varphi_{x'}(w_j')=\varphi_{x}(w_i)\}\Bigr)\,\theta(\rd x)\\
&=\sum_{w_i\in\bar{w}}\int_\varLambda \theta
\bigl(\bar{D}_{\{\varphi_{x}(w_i)\}}(\bar{w}')\bigr)\,\theta(\rd x)
=0\qquad (Q^{\otimes\myp2}\text{-a.s.}),
\end{align*}
since, by assumption (\ref{eq:condB2}), $\theta
\bigl(\bar{D}_{\{\varphi_{x}(w)\}}(\bar{w}')\bigr)=0$ for
$\theta$-a.a.\ $x\in\varLambda$, \,$Q$-a.a.\ $\bar{w}\in
\mathfrak{W}$ and each $w_i\in\bar{w}$, and, moreover, $\bar{w}\in
\mathfrak{W}$ \strut{}contains at most countably many coordinates.
Hence, the right-hand side of (\ref{eq:=0}) vanishes and due to
estimates (\ref{eq:Psi3}) and (\ref{eq:=0}) the claim of the
proposition follows.
\endproof

It is easy to give simple sufficient criteria for condition
(\ref{eq:condB2}) of Proposition \ref{pr:simple}. The first
criterion below is set out in terms of the reference measure
$\theta$, whereas the second one exploits the in-cluster parent
distribution $Q$.

\begin{proposition}\label{pr:a2-1}
Assume that for each $x\in X$, the equation $\varphi_{y}(w)=x$ has
at most one solution $y=y(x;w)$ for every $w\in \bar{w}$ and
$Q$-a.a.\ $\bar{w}\in \mathfrak{W}$. Furthermore, let the measure
$\theta$ be non-atomic, that is, $\theta\{y\}=0$\myp{} for each
$y\in X$. Then condition \textup{(\ref{eq:condB2})} is satisfied.
\end{proposition}
\proof Using formula (\ref{eq:D-union}) and definition (\ref{eq:D}),
we obtain
\begin{align*}
0\le \theta\bigl(\bar{D}_{\{x\}}(\bar{w})\bigr)&\le
\sum_{w_i\in\bar{w}}\theta\bigl(D_{\{x\}}(w_i)\bigr)\\
&=\sum_{w_i\in\bar{w}}\theta\{y\in X: \,\varphi_{y}(w_i)=x\}\\
&=\sum_{w_i\in\bar{w}}\theta\{y(x;w_i)\}=0,
\end{align*}
since the measure $\theta$ is non-atomic.
\endproof

\begin{proposition}\label{pr:a2-2}
Suppose that the in-cluster configurations a.s.\ have no fixed
points, that is, for any $x\in X$ and $\theta$-a.a.\ $y\in X$,
\begin{equation}\label{eq:Q=0}
Q\bigl\{\bar{w}\in\mathfrak{W}: \,\exists\,w_i\in\bar{w}\ \,
\text{such that}\ \,\varphi_y(w_i)=x\bigr\}=0.
\end{equation}
Then condition \textup{(\ref{eq:condB2})} follows.
\end{proposition}
\proof Observe that identity (\ref{eq:sigma=int}) together with the
change of measure (\ref{eq:x0->x}) yields, for each $x\in X$,
\begin{align*}
\int_{\mathfrak{W}}
\theta\bigl(\bar{D}_{\{x\}}(\bar{w})\bigr)\,Q(\rd\bar{w})
&=\int_{X}\biggl(\int_{\frakX_{\{x\}}}\!\eta_{x'}(\rd\bar{y})\biggr)\,\theta(\rd{x}')\\
&=\int_{X} Q\bigl(\bar{\varphi}^{-1}_{x'}(\frakX_{\{x\}})\bigr)\,\theta(\rd{x}')\\
&=\int_{ X} Q\{\bar{w}\in\mathfrak{W}: \,\varphi_{x'}(w_i)=x\ \
\text{for some}\ \ w_i\in\bar{w}\} \,\theta(\rd x')=0,
\end{align*}
according to (\ref{eq:Q=0}). Thus, the proof is complete.
\endproof

\section{Quasi-invariance and integration by parts}
\label{sec:4}

From now on, we assume that the probability distribution of the
centre process (see the beginning of Section \ref{sec:3}) satisfies
the following natural condition.
\begin{condition}\label{condPROP}
The measure $\mu$ is supported on the proper configuration space,
$\mu(\varGamma_{X})=1$; that is, $\mu$-a.a.\ configurations $\gamma$
are locally finite and simple.
\end{condition}
Note that in this case the measure $\hat{\mu}$ is concentrated on
the marked configuration space
\begin{equation*}
\varGamma_{X}(\mathfrak{X}):=\{\hatgamma\in
\varGamma_{\mathcal{Z}}:\,p_{X}(\hatgamma)\in \varGamma_{X}\}.
\end{equation*}

Let us also assume that $X$ is a Riemannian manifold (with a fixed
Riemannian structure). Our aim in this section is to prove the
quasi-invariance of the measure $\mucl$ with respect to the action
of compactly supported diffeomorphisms of $X$
(Section~\ref{sec:4.2}), and to establish an IBP formula
(Section~\ref{sec:4.3}). We begin in Section~\ref{sec:4.1} with a
brief description of some convenient ``manifold-like'' concepts and
notations first introduced in \cite{AKR1} (see also
\cite[\S\myp4.1]{BD3}), which furnish a suitable framework for
analysis on configuration spaces.

\subsection{Differentiable functions on configuration spaces}
\label{sec:4.1}

Let $T_{x}X$ be the tangent space of $X$ at point $x\in X$, with the
corresponding (canonical) inner product denoted by a ``fat''
dot~$\CDD$\,. The gradient on $X$ is denoted by $\nabla $. Following
\cite{AKR1}, we define the ``tangent space'' of the configuration
space $\varGamma_{X}$ at $\gamma \in \varGamma_{X}$ as the Hilbert
space $T_{\gamma}\varGamma_{X}:=L^{2}(X\rightarrow TX;\,\rd\gamma)$,
or equivalently $T_{\gamma}\varGamma_{X}=\bigoplus_{x\in
\gamma}T_{x}X$. The scalar product in $T_{\gamma}\varGamma_{X}$ is
denoted by $\langle \cdot,\cdot \rangle_{\gamma}$, with the
corresponding norm $|{\cdot}|_{\gamma}$. A vector field $V$ over
$\varGamma_{X}$ is a map $\varGamma_{X}\ni \gamma \mapsto V(\gamma
)=(V(\gamma)_{x})_{x\in \gamma}\in T_{\gamma}\varGamma_{X}$ . Thus,
for vector fields $V_{1},V_{2}$ over $\varGamma_{X}$ we have
\begin{equation*}
\left\langle V_{1}(\gamma),V_{2}(\gamma)\right\rangle_{\gamma}
=\sum_{x\in \gamma} V_{1}(\gamma)_{x}\CDD V_{2}(\gamma)_{x},\qquad
\gamma \in \varGamma_{X}.
\end{equation*}

For $\gamma \in \varGamma_{X}$ and $x\in \gamma $, denote by
${\mathcal{O}}_{\gamma,\myp x}$ an arbitrary open neighborhood of
$x$ in $X$ such that ${\mathcal{O}}_{\gamma,\myp x}\cap \gamma
=\{x\}$. For any measurable function $F:\varGamma_{X}\rightarrow
{{\mathbb{R}}}$, define the function $F_{x}(\gamma ,\cdot
):{\mathcal{O}}_{\gamma,\myp x}\rightarrow \RR$ by
$F_{x}(\gamma,y):=F((\gamma \setminus\{x\})\cup\{y\})$, and set
\begin{equation*}
\nabla_{\myn x}F(\gamma):=\left. \nabla
F_{x}(\gamma,y)\right\vert_{y=x},\qquad x\in X,
\end{equation*}
provided that $F_{x}(\gamma ,\cdot)$ is differentiable at $x$.

Recall that for a function $\phi :X\rightarrow \mathbb{R}$ its
support $\supp\phi $ is defined as the closure of the set $\{x\in
X\!:\phi (x)\neq 0\}$. Denote by ${\mathcal{FC}}(\varGamma_{X})$ the
class of functions on $\varGamma_{X}$ of the form
\begin{equation}\label{local-funct}
F(\gamma)=f(\langle \phi_{1},\gamma \rangle ,\dots ,\langle
\phi_{k},\gamma \rangle),\qquad \gamma \in \varGamma_{X},
\end{equation}
where $k\in \mathbb{N}$, \thinspace $f\in C_{b}^{\infty
}(\mathbb{R}^{k})$ ($:=$ the set of $C^{\infty}$-functions on
${\mathbb{R}}^{k}$ globally bounded together with all their
derivatives), and $\phi_{1},\dots ,\phi_{k}\in C_{0}^{\infty}(X)$
($:=$ the set of $C^{\infty}$-functions on $X$ with compact
support). Each $F\in {\mathcal{FC}}(\varGamma_{X})$ is local, that
is, there is a compact $B\subset X$ (e.g., $B=\cup_{j=1}^k
\supp\phi_j$) such that $F(\gamma)=F(\gamma \cap B)$ for all $\gamma
\in \varGamma_{X}$. Thus, for a fixed $\gamma $ there are finitely
many non-zero derivatives $\nabla_{\myn x}F(\gamma)$.

For a function $F\in {\mathcal{FC}}(\varGamma_{X})$ its
$\varGamma$-gradient $\nabla^{\varGamma}F$ is defined as
\begin{equation}
\nabla^{\varGamma}\mynn F(\gamma):=(\nabla_{\myn x} F(\gamma))_{x\in
\gamma}\in T_{\gamma}\varGamma_{X},\qquad \gamma \in \varGamma_{X},
\label{eq:G-gradient}
\end{equation}
so the directional derivative of $F$ along a vector field $V$ is given by
\begin{equation*}
\nabla_{\myn V}^{\varGamma}\myp F(\gamma):=\langle
\nabla^{\varGamma} \mynn
F(\gamma),V(\gamma)\rangle_{\gamma}=\sum_{x\in \gamma} \nabla_{\myn
x}F(\gamma) \CDD V(\gamma)_{x},\qquad \gamma \in \varGamma_{X}.
\end{equation*}%
Note that the sum here contains only finitely many non-zero terms.

Further, let ${\mathcal{FV}}(\varGamma_{X})$ be the class of
cylinder vector fields $V$ on $\varGamma_{X}$ of the form
\begin{equation}\label{vf}
V(\gamma)_{x}=\sum_{i=1}^{k}G_{i}(\gamma)\mypp v_{i}(x)\in
T_{x}X,\qquad x\in X,
\end{equation}
where $G_{i}\in {\mathcal{FC}}(\varGamma_{X})$ and $v_{i}\in
\Vect_{0}(X)$ ($:=$ the space of compactly supported $C^{\infty
}$-smooth vector fields on $X$), \thinspace $i=1,\dots,k$
\,($k\in\NN$). Any vector filed $v\in \Vect_{0}(X)$ generates a
constant vector field $V$ on $\varGamma_{X}$ defined by
$V(\gamma)_{x}:=v(x)$. We shall preserve the notation $v$ for it.
Thus,
\begin{equation}
\nabla_{\myn v}^{\varGamma}F(\gamma)=\sum_{x\in \gamma} \nabla_{\myn
x}F(\gamma) \CDD v(x),\qquad \gamma \in \varGamma_{X}.
\label{eq:grad-new}
\end{equation}

The approach based on ``lifting'' the differential structure from
the underlying space $X$ to the configuration space $\varGamma_{X}$
as described above can also be applied to the spaces
$\mathfrak{X}=\bigsqcup_{\myp n=0}^{\infty}X^{n}$,
$\mathcal{Z}=X\times \mathfrak{X}$ and $\varGamma_{\mathfrak{X}}$,
$\varGamma_{\mathcal{Z}}$, respectively. In such cases, we will use
the analogous notations as above without further explanation.

\subsection{Quasi-invariance}\label{sec:4.2}

In this section, we discuss the property of quasi-invariance of the
measure $\mucl$ with respect to diffeomorphisms of $X$. Let us start
by describing how diffeomorphisms of $X$ act on configuration
spaces. For a measurable map $\varphi :X\rightarrow X$, its
\textit{support} $\supp\varphi $ is defined as the closure of the
set $\{x\in X:\,\varphi(x)\ne x\}$. Let $\Diff_{0}(X)$ be the group
of diffeomorphisms of $X$ with \textit{compact support}. For any
$\varphi \in {\mathrm{Diff}}_{0}(X)$, consider the corresponding
``diagonal'' diffeomorphism $\bar{\varphi}:\mathfrak{X}\rightarrow
\mathfrak{X}$ acting on each constituent space $X^{n}$ ($n\in
\ZZ_{+}$) as
\begin{equation}\label{eq:phi-hat0}
X^{n}\ni \bar{y}=(y_{1},\dots ,y_{n})\mapsto
\bar{\varphi}(\bar{y}):=(\varphi (y_{1}),\dots ,\varphi (y_{n}))\in
X^{n}.
\end{equation}
Finally, we introduce a special class of diffeomorphisms
$\hat{\varphi}$ on $\mathcal{Z}$ acting only in the
$\bar{y}$-coordinate,
\begin{equation}\label{eq:hat-phi}
\hat{\varphi}(z):=(x,\bar{\varphi}(\bar{y})), \qquad
z=(x,\bar{y})\in \mathcal{Z}.
\end{equation}

\begin{remark}\label{supp}
Despite $K_\varphi:= \supp\varphi $ is compact in $X$, the support
of the diffeomorphism $\hat{\varphi}$ (again defined as the closure
of the set $\{z\in \mathcal{Z}: \hat{\varphi}(z)\ne z\}$) is given
by $\supp\hat{\varphi}=X\times \frakX_{K_\varphi}$ (see
(\ref{eq:frakX_B})), where $\frakX_{K_\varphi}$ is \textit{not}
compact in the topology of $\frakX$ (cf.\ Remark~\ref{rm:compact}).
\end{remark}

In a standard fashion, the maps $\varphi $, $\bar\varphi$ and
$\hat{\varphi}$ can be lifted to measurable ``diagonal''
transformations (denoted by the same letters) of the configuration
spaces $\varGamma_{X}$, $\varGamma_{\frakX}$ and
$\varGamma_{\mathcal{Z}}$, respectively:
\vspace{-.8pc}
\begin{align}
\label{di*1} \varGamma_{X}\ni \gamma \mapsto \varphi (\gamma):=
&\{\varphi (x),\ x\in \gamma \}\in \varGamma_{X},\\[.2pc]
\label{di*2} \varGamma_{\frakX}\ni \bar{\gamma} \mapsto
\bar{\varphi}(\bar{\gamma}):=
&\{\bar{\varphi}(\bar{y}),\ \bar{y}\in \bar{\gamma} \}\in \varGamma_{\frakX},\\[.2pc]
\label{di*3} \varGamma_{\mathcal{Z}}\ni \hatgamma\mapsto
\hat{\varphi}(\hatgamma):=&\{\hat{\varphi}(z),\ z\in\hatgamma\}\in
\varGamma_{\mathcal{Z}}.
\end{align}

The following lemma shows that the operator $\mathfrak{q}$ commutes
with the action of diffeomorphisms (\ref{di*1}) and (\ref{di*3}).

\begin{lemma}
For any diffeomorphism $\varphi \in {\mathrm{Diff}}_0(X)$ and the
corresponding diffeomorphism $\hat{\varphi}$, it holds
\begin{equation}  \label{comm1}
\varphi\circ\mathfrak{q}=\mathfrak{q}\circ \hat{\varphi}.
\end{equation}
\end{lemma}

\proof The statement follows from definition (\ref{eq:proj}) of the
map $\mathfrak{q}$ in view of the structure of diffeomorphisms
$\varphi $ and $\hat{\varphi}$ (see (\ref{eq:hat-phi}), (\ref{di*1})
and (\ref{di*3})).
\endproof

Assume that, for all $x\in X$, the measure $\eta_x$ is absolutely
continuous with respect to the Riemannian volume $\rd\bar{y}$ on
$\mathfrak{X}$ and, moreover,
\begin{equation}\label{QI}
h_{x}(\bar{y}):=\frac{\eta_{x}(\rd\bar{y})}{\rd \bar{y}}>0 \qquad
{\text{for \,a.a.}}\ \,\bar{y}\in {\mathfrak{X}}.
\end{equation}
This implies that the measure $\eta_{x}$ is quasi-invariant with
respect to the action of transformations
$\bar{\varphi}:{\mathfrak{X}}\rightarrow\mathfrak{X}$ \,($\varphi
\in\Diff_{0}(X)$), that is, the measure
$\bar{\varphi}^{\ast}\eta_{x}$ is absolutely continuous with respect
to $\eta_{x}$ with the Radon--Nikodym density
\begin{equation}\label{density'}
\rho_{\eta}^{\bar{\varphi}}(x,\bar{y}):=\frac{\rd(\bar{\varphi}
^{\ast}\eta_{x})}{\rd\eta_{x}}(\bar{y})=\frac{h_{x}(\bar{\varphi}
^{-1}(\bar{y}))}{h_{x}(\bar{y})}\,J_{\bar{\varphi}}(\bar{y})^{-1}
\end{equation}
(we set $\rho_{\eta}^{\bar{\varphi}}(x,\bar{y})=1$ if
$h_{x}(\bar{y})=0$ or $h_{x}(\bar{\varphi}^{-1}(\bar{y}))=0$). Here
$J_{\bar{\varphi}}(\bar{y})$ is the Jacobian determinant of the
diffeomorphism $\bar{\varphi}$; due to the diagonal structure of
$\bar{\varphi}$ (see (\ref{eq:phi-hat0})) we have
$J_{\bar{\varphi}}(\bar{y})=\prod_{y_{i}\in
\bar{y}}J_{\varphi}(y_{i})$, where $J_{\varphi}(y)$ is the Jacobian
determinant of $\varphi $.

\begin{theorem}\label{th:inv}
The measure $\hat{\mu}$ is quasi-invariant with respect to the
action of $\hat{\varphi}$ on $\varGamma_{\mathcal{Z}}$ defined by
formula \textup{(\ref{eq:hat-phi})}, with the Radon--Nikodym density
$R_{\hat{\mu}}^{\hat{\varphi}}=\rd(\hat{\varphi}^{\ast}
\hat{\mu})/\rd\hat{\mu}$ given by
\begin{equation}\label{RND}
R_{\hat{\mu}}^{\hat{\varphi}}({\hatgamma})=\prod_{z\in
\hatgamma}\myp\rho_{\eta}^{\bar{\varphi}}(z),\qquad {\hatgamma}\in
\varGamma_{\mathcal{Z}}.
\end{equation}
Moreover, $R_{\hat{\mu}}^{\hat{\varphi}}\in
L^{1}(\varGamma_{\mathcal{Z}},\hat{\mu})$.
\end{theorem}

\proof First of all, note that $\rho_{\varphi}(z)=1$ for any
$z=(x,\bar{y})\notin\supp\hat{\varphi}=X\times\frakX_B=:
\mathcal{Z}_{K_{\varphi}}$, where $K_{\varphi}=\supp\varphi$ (see
Remark \ref{supp}), and $\hat{\sigma}(\mathcal{Z}_{K_{\varphi}})=
\bar{\sigma}(\frakX_{K_{\varphi}})<\infty$ by Condition
\ref{condition5} (see (\ref{eq:sigmaZ})). Therefore,
$\hatgamma(\mathcal{Z}_{K_{\varphi}})<\infty $ for $\hat{\mu}$-a.a.\
configurations $\hatgamma\in \varGamma_{\mathcal{Z}}$, hence the
product in (\ref{RND}) contains only finitely many terms different
from $1$ and so the function
$R_{\hat{\mu}}^{\hat{\varphi}}(\hatgamma)$ is well defined.
Moreover, it satisfies the ``localisation'' equality
\begin{equation}
R_{\hat{\mu}}^{\hat{\varphi}}(\hatgamma)=
R_{\hat{\mu}}^{\hat{\varphi}}(\hatgamma\cap
\mathcal{Z}_{K_{\varphi}})\qquad \text{for}\ \
\hat{\mu}\text{-a.a.}\ \,\hatgamma\in \varGamma_{\mathcal{Z}}.
\label{eq:R=R}
\end{equation}

Now, using definitions (\ref{di*1}), (\ref{di*2}) and
(\ref{eq:eta-gamma}), we obtain
\begin{align}
\notag
\int_{\varGamma_{\mathcal{Z}}}F(\hatgamma)\,\hat{\varphi}^{\ast}\hat{\mu}(\rd\hatgamma)
&=\int_{\varGamma_{\mathcal{Z}}}F(\hat{\varphi}(\hatgamma))
\,\hat{\mu}(\rd\hatgamma) \\
\notag
&=\int_{\varGamma_{X}}\left(\int_{\frakX^{\gamma}}
F(\gamma,\bar{\varphi}(\bar{y}^{\gamma}))
\,\eta^{\gamma}(\rd\bar{y}^{\gamma})\right) \mu(\rd\gamma) \\
\notag
&=\int_{\varGamma_{X}}\left(\int_{\frakX^{\gamma}}F(\gamma,\bar{\varphi}(\bar{y}^{\gamma}))
\,{\textstyle\bigotimes\limits_{x\in\gamma}}\,\eta_x(\rd\bar{y}_x)\right) \mu(\rd\gamma) \\
&=\int_{\varGamma_{X}}\left(\int_{\frakX^{\gamma}}F(\gamma,
\bar{y}^{\gamma}) \,{\textstyle
\bigotimes\limits_{x\in\gamma}}\,\bar{\varphi}^*\eta_x(\rd\bar{y}_x)\right)
\mu(\rd\gamma). \label{eq:wait}
\end{align}
Furthermore, by the quasi-invariance property of the measure
$\eta_x$ (see formula (\ref{density'}) for the density), the
right-hand side of (\ref{eq:wait}) is represented in the form
\begin{align*} \int_{\varGamma_{X}}\left(
\int_{\frakX^{\gamma}}F(\gamma,\bar{y}^{\gamma})
\prod_{x\in\gamma}\rho_{\eta}^{\bar{\varphi}}(x,\bar{y}_x) \
\eta^{\gamma}(\rd\bar{y}^{\gamma})\right) \mu(\rd\gamma)
=\int_{\varGamma_{\mathcal{Z}}}F(\hatgamma)
R_{\hat{\mu}}^{\hat{\varphi}}(\hatgamma)\,\hat{\mu}(\rd\hatgamma),
\end{align*}
which proves the quasi-invariance of $\hat{\mu}{}$. In particular,
for $F\equiv 1$  this yields
$\int_{\varGamma_{\mathcal{Z}}}\!R_{\hat{\mu}}^{\hat{\varphi}}(\hatgamma)
\,\hat{\mu}{}(\rd\hatgamma)\allowbreak=1$, and hence
$R_{\hat{\mu}}^{\hat{\varphi}}\in
L^{1}(\varGamma_{\mathcal{Z}},\hat{\mu})$, as claimed.
\endproof

Let $\mathcal{I}_{\mathfrak{q}}:
L^{\infty}(\varGamma_{X},\mucl)\rightarrow
L^{\infty}(\varGamma_{\mathcal{Z}},\hat{\mu})$ be the isometry
defined by the map ${\mathfrak{q}}$ (see (\ref{eq:proj})),
\begin{equation}\label{eq:I}
({\mathcal{I}_{\mathfrak{q}}}F)({\hatgamma}):=F\circ
\mathfrak{q}({\hatgamma}),\qquad
{\hatgamma}\in\varGamma_{\mathcal{Z}}.
\end{equation}
The adjoint operator $\mathcal{I}_{\mathfrak{q}}^{\ast}$ is a
bounded operator on the corresponding dual spaces,
\begin{equation}\label{eq:I*}
\mathcal{I}_{\mathfrak{q}}^{\ast }: \myp L^\infty
(\varGamma_{\mathcal{Z}},\hat{\mu})^{\prime}\rightarrow
L^{\infty}(\varGamma_{X},\mucl)^{\prime}.
\end{equation}

\begin{lemma}\label{lm:L1}
The operator $\mathcal{I}_{\mathfrak{q}}^{\ast }$ defined by
\textup{(\ref{eq:I*})} can be restricted to the operator
\begin{equation}\label{eq:I*L1}
\mathcal{I}_{\mathfrak{q}}^{\ast}:\myp
L^{1}(\varGamma_{\mathcal{Z}},\hat{\mu})\rightarrow
L^{1}(\varGamma_{X},\mucl).
\end{equation}
\end{lemma}

\proof It is known (see \cite{GL}) that, for any $\sigma$-finite
measure space $(M,\mu)$, the corresponding space $L^{1}(M,\mu)$ can
be identified with the subspace $V$ of the dual space
$L^{\infty}(M,\mu)^{\prime}$ consisting of all linear functionals on
$L^{\infty}(M,\mu)$ continuous with respect to the bounded
convergence in $L^{\infty}(M,\mu)$. That is, $\ell\in V$ if and only
if $\ell(\psi_{n})\rightarrow 0$ for any $\psi_{n}\in
L^{\infty}(M,\mu)$ such that $|\psi_{n}|\leq 1$ and
$\psi_{n}(x)\rightarrow 0$ as $n\to\infty$ for $\mu$-a.a.\ $x\in M$.
Hence, to prove the lemma it suffices to show that, for any $F\in
L^{1}(\varGamma_{\mathcal{Z}},\hat{\mu})$, the functional
$\mathcal{I}_{\mathfrak{q}}^{\ast }F\in L^{\infty
}(\varGamma_{\mathcal{Z}},\hat{\mu})^{\prime}$ is continuous with
respect to bounded convergence in
$L^{\infty}(\varGamma_{\mathcal{Z}},\hat{\mu})$. To this end, for
any sequence $(\psi_{n})$ in $L^{\infty }(\varGamma_{X},\mucl)$ such
that $|\psi_{n}|\leq 1$ and $\psi_{n}(\gamma )\rightarrow 0$ for
$\mucl$-a.a.\ $\gamma \in \varGamma_{X}$, we have to prove that
$\mathcal{I}_{\mathfrak{q}}^{\ast }F(\psi_{n})\rightarrow 0$.

Let us first show that
$\mathcal{I}_{\mathfrak{q}}\psi_{n}(\hatgamma) \equiv
\psi_{n}(\mathfrak{q}(\hatgamma))\rightarrow 0$ for
$\hat{\mu}$-a.a.\ $\hatgamma\in \varGamma_{\mathcal{Z}}$. Set
\begin{align*}
A_{\psi }:=& \{\gamma \in \varGamma_{X}:\psi_{n}(\gamma)\rightarrow
0\}\in \mathcal{B}(\varGamma_{X}), \\[0.3044pc]
\hat{A}_{\psi }:=& \{{\hatgamma}\in
\varGamma_{\mathcal{Z}}:\psi_{n}(\mathfrak{q}(\hatgamma))\rightarrow
0\}\in \mathcal{B}(\varGamma_{\mathcal{Z}}),
\end{align*}
and note that $\hat{A}_{\psi }=\mathfrak{q}^{-1}(A_{\psi })$; then,
recalling relation (\ref{eq:cl*}), we get
\begin{equation*}
\hat{\mu}(\hat{A}_{\psi
})=\hat{\mu}\bigl(\mathfrak{q}^{-1}(A_{\psi})\bigr)
=\mucl(A_{\psi})=1,
\end{equation*}
as claimed. Now, by the dominated convergence theorem this implies
\begin{equation*}
\mathcal{I}_{\mathfrak{q}}^{\ast}
F(\psi_{n})=\int_{\varGamma_{\mathcal{Z}}}
F({\hatgamma})\,\mathcal{I}_{\mathfrak{q}}\psi_{n}(\hatgamma)\,
\hat{\mu}(\rd\hatgamma)\to 0,
\end{equation*}
and the proof is complete.
\endproof

\begin{corollary}
For any measurable functions $F\in L^{\infty}(\varGamma_{X},\mucl)$
and $G\in L^{1}(\varGamma_{\mathcal{Z}},\hat{\mu})$, we have the
identity
\begin{equation}\label{eq:F*G}
\int_{\varGamma_{\mathcal{Z}}}
\!G(\hatgamma)\,\mathcal{I}_{\mathfrak{q}}F(\hatgamma)\,\hat{\mu}(\rd\hatgamma)
=\int_{\varGamma_{X}} \!F(\gamma)\,\mathcal{I}_{\mathfrak{q}}^{\ast}
G(\gamma)\,\mucl(\rd\gamma).
\end{equation}
\end{corollary}
Taking advantage of Theorem \ref{th:inv} and applying the projection
construction, we obtain our main result in this section.

\begin{theorem}\label{q-inv}
The cluster measure $\mucl$ is quasi-invariant with respect to the
action of the diffeomorphism group\/ $\Diff_{0}(X)$ on
$\varGamma_{X}$. The corresponding Radon--Nikodym density is given
by $R_{\mucl}^{\varphi}=\mathcal{I}_{\mathfrak{q}}^{\ast
}R_{\hat{\mu}}^{\hat{\varphi}}\in L^1(\varGamma_X,\mucl)$.
\end{theorem}

\proof Note that, due to (\ref{eq:cl*}) and (\ref{comm1}),
\begin{equation*}
\mucl\circ \varphi^{-1}=\hat{\mu}\circ \mathfrak{q}^{-1}\circ
\varphi^{-1}=\hat{\mu}\circ \hat{\varphi}^{-1}\circ
\mathfrak{q}^{-1}.
\end{equation*}
That is, $\varphi^{\ast}\mucl=\mucl\circ \varphi^{-1}$ is a
push-forward of the measure $\hat{\varphi}^{\ast
}\hat{\mu}=\hat{\mu}\circ \hat{\varphi}^{-1}$ under the map $
\mathfrak{q}$, that is, $\varphi^{\ast}\mucl= \mathfrak{q}^{\ast
}\hat{\varphi}^{\ast}\hat{\mu}$. In particular, if
$\hat{\varphi}^{\ast}\hat{\mu}$ is absolutely continuous with
respect to $\hat{\mu}$ then so is $\varphi^{\ast}\mucl$ with respect
to $\mucl$. Moreover, by formula (\ref{eq:cl*}) and Theorem
\ref{th:inv}, for any $F\in L^{\infty }(\varGamma_{X},\mucl)$ we
have
\begin{equation}\label{eq:L1}
\int_{\varGamma_{X}}F(\gamma)\,\varphi^{\ast}
\mucl(\rd\gamma)=\int_{\varGamma_{\mathcal{Z}}}
{\mathcal{I}_{\mathfrak{q}}}F({\hatgamma})\,\hat{\varphi}^{\ast}
\hat{\mu}(\rd\hatgamma)=\int_{\varGamma_{\mathcal{Z}}}
\mynn{\mathcal{I}_{\mathfrak{q}}}F({\hatgamma})\myp
R_{\hat{\mu}}^{\hat{\varphi}}(\hatgamma)\,\hat{\mu}(\rd\hatgamma).
\end{equation}
By Lemma \ref{lm:L1}, the operator $\mathcal{I}_{\mathfrak{q}}^{\ast
}$ acts from $L^{1}(\varGamma_{\mathcal{Z}},\hat{\mu})$ to
$L^{1}(\varGamma_{X},\mucl)$. Therefore, again using (\ref{eq:cl*})
the right-hand side of (\ref{eq:L1}) can be rewritten as
\begin{equation*}
\int_{\varGamma_{X}}\mynn
F(\gamma)\mypp(\mathcal{I}_{\mathfrak{q}}^{\ast
}R_{\hat{\mu}}^{\hat{\varphi}})(\gamma)\,\mucl(\rd\gamma),
\end{equation*}
which completes the proof.
\endproof

\begin{remark}
 Cluster measure
$\mucl$ on the configuration space $\varGamma_{X}$ can be used to
construct a unitary representation $U$ of the diffeomorphism group
$\Diff_{0}(X)$ by operators in $L^{2}(\varGamma_{X},\mucl)$, given
by the formula
\begin{equation}\label{eq:U}
U_{\varphi}F(\gamma)
=\sqrt{R_{\mucl}^{\varphi}(\gamma)}\,F(\varphi^{-1}(\gamma)),\qquad
F\in L^{2}(\varGamma_{X},\mucl).
\end{equation}
Such representations, which can be defined for arbitrary
quasi-invariant measures on $\varGamma_{X}$, play a significant role
in the representation theory of the group $\Diff_{0}(X)$
\cite{Ism,VGG} and quantum field theory \cite{GGPS,Goldin}. An
important question is whether the representation (\ref{eq:U}) is
irreducible. According to \cite{VGG}, this is equivalent to the
$\Diff_{0}(X)$-ergodicity of the measure $\mucl$, which in our case
is equivalent to the ergodicity of the measure $\hat{\mu}$ with
respect to the group of transformations $\hat{\varphi}$ \,($\varphi
\in \Diff_{0}(X)$).
\end{remark}

\subsection{Integration-by-parts (IBP) formulae }
\label{sec:4.3}

In this section, we assume that the conditions of Lemma \ref{lm:M^n}
are satisfied with $n=1$. Thus, the measures $\mu$, $\hat{\mu}$ and
$\mucl$ belong to the corresponding $\mathcal{M}^{1}$-classes. It is
also assumed, as before, that for each $x\in X$ the measure $\eta_x$
is absolutely continuous with respect to the Riemannian volume
$\rd\bar{y}$ on $\frakX$, with the Radon--Nykodym density
$h_x(\bar{y})$.

\subsubsection{Integration by parts for the cluster distributions $\eta_x$.}
Let $v\in \mathrm{Vect}_{0}(X)$ (:= the space of compactly supported
smooth vector fields on $X$), and define a ``vertical'' vector field
$\hat{v}$ on $\mathcal{Z}$ by the formula
\begin{equation}\label{eq:v-hat}
\hat{v}(x,\bar{y}):=(v(y_{i}))_{y_i\in\bar{y}},\qquad \bar{y}
=(y_{i})\in\mathfrak{X}.
\end{equation}
Observe that if the density $h_x(\bar{y})$ is differentiable
($\rd\bar{y}$-a.e.) then the measure $\eta_{x}$ satisfies the IBP
formula (see, e.g., \cite[\S1.3, \S\myp2.4]{BelDal}; cf.\
\cite[\S\myp5.1.3, p.~207]{Bo})
\begin{equation}\label{eq:etaIBP}
\int_{{\mathfrak{X}}}\nabla^{\hat{v}}\mynn
f(\bar{y})\,\eta_{x}(\rd\bar{y})=-\int_{\mathfrak{X}}f(\bar{y})
\,\beta_{\eta}^{\hat{v}}(x,\bar{y})\,\eta_x(\rd\bar{y}),\qquad f\in
C_{0}^{\infty}(\mathfrak{X}),
\end{equation}
where $\nabla^{\hat{v}}$ is the derivative along the vector field
$\hat{v}$ and
\begin{equation}\label{eq:log-der}
\beta_{\eta}^{\hat{v}}(x,\bar{y}):=(\beta_{\eta}(x,\bar{y}),
\hat{v}(x,\bar{y}))_{T_{\bar{y}}{\mathfrak{X}}}+\Div\hat{v}(x,\bar{y})
\end{equation}
is the logarithmic derivative of
$\eta_{x}(\rd\bar{y})=h_{x}(\bar{y})\mypp\rd\bar{y}$ along
$\hat{v}$, expressed in terms of the vector logarithmic derivative
\begin{equation} \label{eq:nabla-h}
\beta_{\eta}(x,\bar{y}):=\frac{\nabla
h_{x}(\bar{y})}{h_{x}(\bar{y})}\in T_{\bar{y}}{\mathfrak{X}},\qquad
(x,\bar{y})\in X\times {\mathfrak{X}}.
\end{equation}

Denote for brevity
$$
\|\bar{y}\|_1:=\sum_{y_i\in\bar{y}}|y_i|,\qquad \bar{y}\in\frakX.
$$
\begin{lemma}\label{lm:4.6}
Suppose that
$\int_{\mathcal{Z}_B}\|\beta_\eta(z)\|_1^{n}\,\hat{\sigma}(\rd{z})<\infty$
for any compact $B\subset X$, and assume that condition\/
\textup{(\ref{eq:N<})} is satisfied. Let $\hat{v}$ be a vector filed
on $\mathcal{Z}$ defined by \textup{(\ref{eq:v-hat})} with
$v\in\Vect_0(X)$. Then $\beta_{\eta}^{\hat{v}}\in
L^{n}(\mathcal{Z},\hat{\sigma})$.
\end{lemma}

\proof To show that $\beta_{\eta}^{\hat{v}}\in
L^{n}(\mathcal{Z},\hat{\sigma})$, it suffices to check that each of
the two terms on the right-hand side of (\ref{eq:log-der}) belongs
to $L^{n}(\mathcal{Z},\hat{\sigma})$. Setting $b_{v}:=\sup_{x\in
X}|v(x)|<\infty$ and noting that $K_{v}:=\supp v$ is a compact in
$X$, we have
\begin{align}
\int_{\mathcal{Z}}|
(\beta_\eta(z),\hat{v}(z))|^{n}\,\hat{\sigma}(\rd{z})&\leq
\int_{X}\int_{\mathfrak{X}_{K_v}} \!\left(\myn\sum_{y_{i}\in
\bar{y}}|\beta_\eta(x,\bar{y})_{i}|
\cdot |v(y_{i})|\right)^{\!n} \eta_{x}(\rd\bar{y})\,\theta(\rd{x})  \notag \\
& \leq (b_v)^n
\int_{X}\int_{\mathfrak{X}_{K_v}}\left(\myn\sum_{y_{i}\in
\bar{y}}|\beta_{\eta}(x,\bar{y})_{i}|\right)
^{\!n}\eta_{x}(\rd\bar{y})\,\theta (\rd x)
\notag\\
\label{eq:beta<infty} & = (b_v)^n
\int_{\mathcal{Z}_{K_v}}\|\beta_\eta(z)\|_1^{n}\,\hat{\sigma}(\rd
z)<\infty,
\end{align}
by the first hypothesis of the theorem. Similarly, denoting
$d_v:=\sup_{x\in X}|\Div v(x)|<\infty$, we obtain
\begin{align}
\int_{\mathcal{Z}}|\mynn\Div\hat{v}(x,\bar{y})|^{n}\,\hat{\sigma}
(\rd{x}\times \rd\bar{y})&=\int_{\mathcal{Z}}
\left(\myn\sum_{y_{i}\in \bar{y}}|\mynn\Div v
(y_{i})|\right)^{\!n}\eta_{x}(\rd\bar{y})\,\theta(\rd{x})\,\notag\\
\label{eq:div<infty} & \le (d_v)^n
\int_{X}\int_{\frakX_{K_{v}}}N_{X}(\bar{y})^{n}
\,\eta_{x}(\rd\bar{y})\,\theta(\rd{x})<\infty,
\end{align}
according to assumption (\ref{eq:N<}). As a result, combining bounds
(\ref{eq:beta<infty}) and (\ref{eq:div<infty}), we see that
$\beta_{\eta}^{\hat{v}}\in L^{n}(\mathcal{Z},\hat{\sigma})$, as
claimed.
\endproof

Let us define the space $H_{\rm loc}^{1,\myp n}(\mathfrak{X)}$
($n\geq 1$) as the set of functions $f\in
L^{n}(\mathfrak{X},\rd\bar{y})$ satisfying, for any compact
$B\subset X$, the condition
\begin{equation}\label{eq:Sobolev}
\kappa^B_n(f):=\int_{\mathfrak{X}}\|\nabla\myn
f(\bar{y})\|_1^{n}\,\rd\bar{y}\equiv \int_{\mathfrak{X}_B}\left(\myn
\sum_{y_{i}\in \bar{y}}\,|\nabla_{\myn
y_{i}}f(\bar{y})|\right)^{\!n}\rd\bar{y}<\infty.
\end{equation}
Due to the elementary inequality $(|a|+|b|)^{n}\leq
2^{n-1}\bigl(|a|^{n}+|b|^{n}\bigr)$, $H^{1,\myp n}(\mathfrak{X)}$ is
a linear space.

The integrability condition in Lemma \ref{lm:4.6} on the vector
logarithmic derivative $\beta_\eta(z)$ can be characterised as
follows.
\begin{lemma}\label{lm:B^n}
Assume that, for some integer\/ $n\geq 1$,\/
$h_{x}^{1/n}\mynn\in H_{\rm loc}^{1,\myp n}(\mathfrak{X)}$ for
$\theta$-a.a.\ $x\in X$. Then\/
$\int_{\mathcal{Z}_B}\|\beta_\eta(z)\|_1^{n}\,
\sigma_{\mathcal{Z}}(\rd z)<\infty$\/ if and only if\/ for any
compact $B\subset X$
\begin{equation}\label{eq:kappan}
\int_{X} \kappa^B_n(h_{x}^{1/n})\,\theta(\rd x)<\infty.
\end{equation}
\end{lemma}

\proof Substituting formulae (\ref{QI}) and (\ref{eq:nabla-h}), it
is easy to see that
\begin{align*}
\int_{\mathcal{Z}_B}\|\beta_\eta(z)\|_1^{n}\,\hat{\sigma}(\rd z)&=
\int_{X}\int_{\mathfrak{X}_B}\left(\myn \sum_{y_{i}\in\bar{y}}
\frac{|\nabla_{\myn y_{i}} h_{x}(\bar{y})|}
{h_{x}(\bar{y})}\right)^{\!n} h_x(\bar{y})\,\rd\bar{y}\;\theta (\rd
x)
\\&=\int_{X}\int_{\mathfrak{X}_B}\left(\myn
\sum_{y_{i}\in\bar{y}} \frac{|\nabla_{\myn y_{i}}h_{x}(\bar{y})|}
{h_{x}(\bar{y})^{1-1/n}}\right)^{\!n}\rd\bar{y}\;\theta (\rd x) \\
&=n^n\int_{X}\int_{\mathfrak{X}_B}\left(\myn\sum_{y_{i}\in
\bar{y}}\,\bigl|\nabla_{\myn
y_{i}}\bigl(h_{x}(\bar{y})^{1/n}\bigr)\bigr|\right)^{n}
\rd\bar{y}\:\theta(\rd x)\\
&=n^n\int_{X}\kappa^B_n (h_{x}^{1/n})\, \theta(\rd x)<\infty,
\end{align*}
according to (\ref{eq:Sobolev}) and (\ref{eq:kappan}).
\endproof

From now on, we assume the following
\begin{condition}\label{cond4.1}
For any compact $B\subset X$, the vector logarithmic derivative
$\beta_\eta$ defined in \textup{(\ref{eq:nabla-h})} satisfies the
integral bound
$$
\int_{\mathcal{Z}_B}\|\beta_\eta(z)\|_1\, \hat{\sigma}(\rd
z)<\infty.
$$
\end{condition}

\subsubsection{Integration by parts for $\eta_x$ as a push-forward measure.}\label{sec:4.3.2}

Using the general IBP framework outlined in Appendix \ref{ap:2}, and
in particular picking up on Remark \ref{rm:sqcup}, let us consider
the special case with
$\mathcal{W}:=\mathfrak{W}\equiv\bigsqcup_{n=1}^\infty W^n$,
\,$\mathcal{Y}:=\mathfrak{X}\equiv\bigsqcup_{n=1}^\infty X^n$ and
$\phi:=\bar{\varphi}_{x}$, where the maps
$\bar{\varphi}_{x}:\mathcal{W}\to \frakX$ \,($x\in X$) are described
in Section \ref{sec:droplet}. We assume that $\varphi _{x}\in
C_{b}^{2}(W,X)$ uniformly in $x\in X$ (i.e., with global constants
bounding the first two derivatives, $\rd\varphi_x(\bar{w})$ and
$\rd^2\varphi_x(\bar{w})$). Furthermore, given a probability measure
$Q$ on $\mathfrak{W}$, consider the family of measures
$\{\eta_{x}\}_{x\in X}$ on $\mathfrak{X}$ defined by (cf.\
(\ref{eq:x0->x}))
\begin{equation}\label{nu-x}
\eta_{x}:=\bar{\varphi}_{x}^{\ast}\myp Q,\qquad x\in X.
\end{equation}

We need the following two integrability conditions on the vector
logarithmic derivative $\beta_{Q}(\bar{w})$ and the number of
components $N_{W}(\bar{w})$ in a (random) vector
$\bar{w}\in\mathfrak{W}$, both involving the $\theta$-measure of the
droplet cluster $\bar{D}_{B}(\bar{w})$ for any compact $B\subset X$
(see (\ref{eq:barD})):
\begin{gather}
\label{beta} \int_{\mathfrak{W}}\left\Vert \beta_{Q
}(\bar{w})\right\Vert_{1}^{n}\,\theta (\bar{D}_{B}(\bar{w}))\,Q
(\rd\bar{w})<\infty,\\
\label{N}
\int_{\mathfrak{W}}N_{W}(\bar{w})^{n}\,\theta(\bar{D}_{B}(\bar{w}))
\,Q(\rd\bar{w})<\infty.
\end{gather}

We can now prove the following result.

\begin{theorem}\label{th:C2}
Suppose that conditions \textup{(\ref{beta})} and \textup{(\ref{N})}
hold for some $n\ge1$. Then the following statements are
true\textup{:}

\textup{(a)} The function $N_{X}(\bar{y})$ satisfies the
integrability condition (\textup{\ref{eq:N<})}, that is, for any
compact set $B\subset X$
\begin{equation}\label{C-eq:N<}
\int_{X}\int_{\mathfrak{X}_{B}}N_{X}(\bar{y})^{n}\,
\eta_{x}(\rd\bar{y})\,\theta(\rd{x})<\infty.
\end{equation}

\textup{(b)} For any $v\in \Vect_{0}(X)$, the measure $\eta_{x}$
satisfies the IBP formula \textup{(\ref{eq:log-der})} with the
corresponding logarithmic derivative $\beta_{\eta }^{\hat{v}}\in
L^{n}(\mathcal{Z},\hat{\sigma})$.
\end{theorem}

\proof (a) By the change of measure (\ref{nu-x}), we obtain (cf.\
the proof of Lemma~\ref{lm:sigma=int})
\begin{align*}
\int_{X}\int_{\mathfrak{X}_{B}}N_{X}(\bar{y})^{n}\,\eta
_{x}(\rd\bar{y})\,\theta(\rd{x})&=\int_{X}\int_{\mathfrak{W}}
\mathbf{1}_{\mathfrak{X}_{B}}(\bar{\varphi}_{x}(\bar{w}))N_{W}(\bar{w})^{n}
\,Q(\rd\bar{w})\,\theta(\rd x)\\
&=\int_{\mathfrak{W}}N_{W}(\bar{w})^{n}\left(\int_{X}
\mathbf{1}_{\mathfrak{X}_{B}}(\bar{\varphi}_{x}(\bar{w}))\,\theta(\rd x)\right)Q(\rd\bar{w})\\
&=\int_{\mathfrak{W}}N_{W}(\bar{w})^{n}
\,\theta(\bar{D}_{B}(\bar{w}))\,Q(\rd\bar{w})<\infty,
\end{align*}
according to condition (\ref{N}), and so the first part of the
theorem is proved.

\smallskip
(b) Recall that the vector field $\hat{v}$ on $\mathcal{Z}=X\times
\mathfrak{X}$ is defined by
\begin{equation}
\hat{v}(x,\bar{y}):=\left( v(y_{i})\right)_{y_i\in\bar{y}} ,\qquad
\bar{y}=(y_{i})\in\mathfrak{X}.
\end{equation}
Then, owing to the component-wise structure of the map
$\bar{\varphi}_{x}$ (cf.\ (\ref{eq:lift})), we have
\begin{equation*}
(\mathcal{I}_{\bar{\varphi}_{x}}\hat{v})(x,\bar{w})
=\bigl((\mathcal{I}_{\varphi_{x}}v)(w_i)\bigr)_{w_i\in\bar{w}}
,\qquad \bar{w}=(w_{i})\in \mathfrak{W}.
\end{equation*}
It is clear that $\hat{v}\in \Vect_{b}^{1}(\mathfrak{X})$. Moreover,
$\mathcal{I}_{\varphi_{x}}v\in \Vect_{b}^{1}(W)$,
\,$\mathcal{I}_{\bar{\varphi}_{x}}\hat{v}\in
\Vect_{b}^{1}(\mathfrak{W})$ uniformly in $x\in X$, which implies
that
\begin{gather}\label{est-1}
C_{1}:=\sup_{x\in X,\,w\in W}
|\myp\mathcal{I}_{\varphi_{x}}v(w)|<\infty,\\
\label{est-2} C_{2}:=\sup_{x\in X,\,w\in W} |\myn\Div
\mathcal{I}_{\varphi_{x}}v(w)|<\infty.
\end{gather}
By Theorem \ref{IBP-proj} and Remark \ref{rm:sqcup} in Appendix
\ref{ap:3}, the measure $\eta_x=\bar{\varphi}_{x}^{\ast}\myp Q$
satisfies the IBP formula (\ref{eq:log-der}) with the logarithmic
derivative
\begin{equation}\label{Log-der}
\beta_{\eta}^{\hat{v}}(x,\bar{y})=\bigl(\mathcal{I}_{\varphi_{x}}^{\ast}
\beta_{Q}^{\mathcal{I}_{\bar{\varphi}_{x}}\hat{v}}\bigr)(\bar{y}),\qquad
x\in X,\ \ \bar{y}\in\frakX,
\end{equation}
where
\begin{align*}
\beta_{Q}^{\mathcal{I}_{\bar{\varphi}_{x}}\hat{v}}(x,\bar{w})
&=\bigl(\beta_{Q}(\bar{w}),\mathcal{I}_{\bar{\varphi}_{x}}\hat{v}(\bar{w})\bigr)
_{T_{\bar{w}}\mathfrak{W}}+\Div \mathcal{I}_{\bar{\varphi}_{x}}
\hat{v}(\bar{w})
\\[.3pc]
&=\sum_{w_{i}\in \bar{w}}\bigl(\beta_{Q}(\bar{w})_{i},
\mathcal{I}_{\varphi_{x}}v(w_{i})\bigr)_{T_{w_{i}}W}+\sum_{w_{i}\in
\bar{w}}\Div\mathcal{I}_{\varphi_{x}}v(w_{i}).
\end{align*}

Let us show that $\beta_{\eta }^{\hat{v}}\in
L^{n}(\mathcal{Z},\hat{\sigma})$. Recall that the map
$\mathcal{I}_{\varphi_{x}}^{\ast}:L^{n}(\mathfrak{W},Q)\rightarrow
L^{n}(\mathfrak{X},\eta_{x})$ is an isometry. Thus, according to
(\ref{Log-der}) and after the change of measure (\ref{nu-x}), we
have
\begin{align*}
\int_{\mathcal{Z}}\bigl|\beta_{\eta}^{\hat{v}}(z)\bigr|^{n}
\,\hat{\sigma}(\rd{z})&\le
\int_{X}\int_{\mathfrak{X}}\bigl\|\beta_{\eta}^{\hat{v}}(x,\bar{y})\bigr\|_1^{n}
\:\hat{\sigma}(\rd{x}\times\rd\bar{y})\\
&=\int_{X}\int_{\mathfrak{W}} \bigl\|
\beta_{Q}^{{\mathcal{I}_{\bar{\varphi}_{x}}\hat{v}}}(x,\bar{w})\bigr\|_1^{n}
\:Q(\rd\bar{w})\,\theta(\rd x).
\end{align*}
Observe that $\supp
{\mathcal{I}_{\bar{\varphi}_{x}}\hat{v}}=\bar{\varphi}_{x}^{-1}(\mathfrak{X}_{K_v})$,
where $K_v:=\supp v$. Then, similarly to the proof of Lemma
\ref{lm:B^n}, we obtain
\begin{align*}
\int_{X}\int_{\mathfrak{W}}&\bigl|\bigl(\beta_{Q}(\bar{w}),
{\mathcal{I}_{\bar{\varphi}_{x}}\hat{v}}(\bar{w})\bigr)_{T_{\bar{w}}\mathfrak{W}}
\bigr|^{n} \,Q(\rd\bar{w})\,\theta (\rd x)\\
&\leq
\int_{X}\int_{\bar{\varphi}_{x}^{-1}(\mathfrak{X}_{K_v})}\left(\sum_{w_{i}\in
\bar{w}}| \beta_{Q}(\bar{w})_{i}|\cdot|
\mathcal{I}_{\varphi_{x}}v(w_{i})|
\right) ^{n}\, Q(\rd\bar{w})\,\theta (\rd x) \\
&\leq \sup_{x\in X,\,w\in W}
|\mathcal{I}_{\varphi_{x}}v(w)|^{n}\int_{\mathfrak{W}}\left(
\sum_{w_{i}\in \bar{w}}\left\vert \beta_{Q} (\bar{w})_{i}\right\vert
\right)^{n}\!\left(\int_{X}
\mathbf{1}_{\mathfrak{X}_{K_v}}(\bar{\varphi}_{x}(\bar{w}))\,\theta
(\rd x)\right)Q(\rd\bar{w})\\
&=C_{1}^{n} \int_{\mathfrak{W}}\left\Vert
\beta_{Q}(\bar{w})\right\Vert_{1}^{n}\,\theta
(\bar{D}_{K_v}(\bar{w}))\,Q(\rd\bar{w})<\infty,
\end{align*}
according to condition (\ref{beta}). Similarly, using bound
(\ref{est-2}) and making the change of measure (\ref{nu-x}), we get
\begin{align*}
\int_{X}\int_{\mathfrak{W}}&|\Div {\mathcal{I}_{\bar{\varphi}_{x}}
\hat{v}}(\bar{w})|^{n} \,Q(\rd\bar{w})\,\theta (\rd{x})\\
&\leq \int_{X}\int_{\mathfrak{W}}\left( \sum_{w{i}\in \bar{w}}\vert
\Div \mathcal{I}_{\varphi_{x}}
v(w_{i})\vert \right) ^{n}\,Q(\rd\bar{w})\,\theta (\rd x) \\
&\leq \sup_{x\in X,\,w\in W}\vert
\Div\mathcal{I}_{\varphi_{x}}v(w)\vert^{n}\int_{X}
\int_{\bar{\varphi}_{x}^{-1}(\mathfrak{X}_{K_v})}
N_{W}(\bar{w})^{n}\,Q(\rd\bar{w})\,\theta(\rd x) \\
&=C_{2}^{n}\int_{X}\int_{\mathfrak{X}_{K_v}}N_{X}(\bar{y})^{n}\,
\eta_{x}(\rd\bar{y})\,\theta(\rd{x})<\infty,
\end{align*}
according to part (a). Thus, part (b) of the theorem is proved.
\endproof

\begin{remark}
Recalling a simple bound (\ref{eq:theta<Nxtheta}) for the
$\theta$-measure of the droplet cluster $\bar{D}_{B}(\bar{w})$, we
observe that, under condition (\ref{eq:C_B}) (see Proposition
\ref{pr:3.9}), conditions (\ref{beta}) and (\ref{N}) of Theorem
\ref{th:C2} specialise, respectively, as follows:
\begin{gather*}
\int_{\mathfrak{W}}\left\Vert \beta_{Q}(\bar{w})\right\Vert
_{1}^{n}N_{W}(\bar{w})\,Q(\rd\bar{w})<\infty, \qquad
\int_{\mathfrak{W}}N_{W}(\bar{w})^{n+1}\,Q(\rd\bar{w})<\infty.
\end{gather*}
Similarly, the assumptions of Proposition \ref{pr:3.10} imply that
$\sup_{\bar{w}\in\mathfrak{W}} \theta(\bar{D}_{B}(\bar{w}))<\infty$
(see the proof), so that conditions (\ref{beta}) and (\ref{N})
transcribe, respectively, as
\begin{gather*}
\int_{\mathfrak{W}}\left\Vert \beta_{Q }(\bar{w})\right\Vert_{1}^{n}
\,Q(\rd\bar{w})<\infty, \qquad
\int_{\mathfrak{W}}N_{W}(\bar{w})^{n}\,Q(\rd\bar{w})<\infty.
\end{gather*}
\end{remark}

\subsubsection{Integration by parts for the cluster measure $\mucl$.} Denote by
${\mathcal{FC}}_{\hat{\sigma}}(\varGamma_{\mathcal{Z}})$ the class
of functions on $\varGamma_{\mathcal{Z}}$ of the form
\begin{equation}\label{sigma-local}
F(\hatgamma)=f(\langle \phi_{1},\hatgamma\rangle,\dots ,\langle
\phi_{k},\hatgamma\rangle),\qquad \hatgamma\in
\varGamma_{\mathcal{Z}},
\end{equation}
where $k\in\NN$, \,$f\in C_{b}^{\infty}(\RR^{k})$ and
$\phi_{1},\dots ,\phi_{k}\in C_{\hat{\sigma}}^{\infty}(\mathcal{Z})$
$:=$ the set of $C^{\infty}$-functions on $\mathcal{Z}$ with
$\hat{\sigma}$-finite support (cf.\ (\ref{local-funct})).

For any $F\in {\mathcal{FC}}(\varGamma_{X})$ we introduce the
function $\hat{F}=\mathcal{I}_{\mathfrak{q}}
F:\varGamma_{\mathcal{Z}}\rightarrow \RR$. It follows from condition
(\ref{eq:sigma<infty}) that
\begin{equation}\label{local-functions}
\hat{F}\in {\mathcal{FC}}_{\hat{\sigma}}(\varGamma_{\mathcal{Z}}).
\end{equation}

\begin{theorem}\label{th:4.7}
The measure $\hat{\mu}$ satisfies the following IBP formula
\begin{equation}\label{eq:4.25}
\int_{\varGamma_{\mathcal{Z}}}\nabla_{\myn \hat{v}}^{\varGamma}
\mynn F(\hatgamma)\,\hat{\mu}(\rd\hatgamma)
=-\int_{\varGamma_{\mathcal{Z}}}
F(\hatgamma)B^{\hat{v}}_{\hat{\mu}}(\hatgamma)\,\hat{\mu}(\rd\hatgamma),
\end{equation}
where
\begin{equation}\label{eq:4.26}
B^{\hat{v}}_{\hat{\mu}}(\hatgamma):=\sum_{(x,\myp\bar{y})\in
\hatgamma}\beta^{\hat{v}}_{\eta}(x,\bar{y})\in L^{1}(\varGamma
_{\mathcal{Z}},\hat{\mu}).
\end{equation}
\end{theorem}

\proof Let us first observe that the integral on the left-hand side
of (\ref{eq:4.25}) is well defined because $\hat{\mu}\in
\mathcal{M}_{\theta }^{1}(\varGamma_{\mathcal{Z}})$. Indeed, the
inclusion (\ref{local-functions}) implies that the function
\begin{equation*}
G(\hatgamma):=\sum_{z\in
\hatgamma}\nabla_{\myn\bar{y}}\hat{F}(\hatgamma)\CDD\hat{v}(\bar{y}),\qquad
z=(x,\bar{y})\in \mathcal{Z},\ \ \ \hatgamma\in
\varGamma_{\mathcal{Z}},
\end{equation*}
is bounded and has $\hat{\sigma}$-finite support, which implies that
$G\in L^{1}(\mathcal{Z},\hat{\sigma})$. Thus the function
\begin{equation*}
\varGamma_{\mathcal{Z}}\ni \hatgamma\mapsto \left\langle G,\gamma
\right\rangle \equiv
\nabla_{\myn\hat{v}}^{\varGamma}\hat{F}(\hatgamma)
\end{equation*}
belongs to $L^{2}(\varGamma_{\mathcal{Z}},\hat{\mu})$ by the
definition of the class $\mathcal{M}_{\theta
}^{1}(\varGamma_{\mathcal{Z}})$.

Using decomposition (\ref{eq:mu-hat-int}) of the measure $\hat{\mu}$
and taking the notational advantage of the one-to-one association
$x\leftrightarrow \bar{y}_x$ for
$(x,\bar{y}_x)\in\hatgamma=(\gamma,\bar{y}^\gamma)$ (see
(\ref{eq:x<->y})), we obtain
\begin{align}
\notag
\int_{\varGamma_{\mathcal{Z}}}\nabla_{\myn\hat{v}}^{\varGamma}
F(\hatgamma)\,\hat{\mu}(\rd\hatgamma)
&=\int_{\varGamma_{X}}\left(\int_{{\frakX}^{\gamma}}\sum_{x\in\gamma}
\nabla_{\myn \bar{y}_x}^{\hat{v}}F(\gamma,\bar{y}^{\gamma})
\,\eta^{\gamma}(\rd\bar{y}^{\gamma})\right) \mu(\rd\gamma)\\
\notag
&=\int_{\varGamma_{X}}\sum_{x\in\gamma}\left(\int_{{\frakX}^{\gamma}}
\nabla_{\myn\bar{y}_x}^{\hat{v}}F(\gamma,\bar{y}^{\gamma})
\,\eta^{\gamma}(\rd\bar{y}^{\gamma})\right) \mu(\rd\gamma)\\
\label{eq:IBP-take1}
&=\int_{\varGamma_{X}}\sum_{x\in\gamma}\left(\int_{{\frakX}^{\gamma}}
\nabla_{\myn \bar{y}_x}^{\hat{v}}F(\gamma,\bar{y}^{\gamma})
\,\textstyle{\bigotimes\limits_{x'\in\gamma}}\eta_{x'}(\rd\bar{y}_{x'})\right)
\mu(\rd\gamma),
\end{align}
by a product structure of $\eta^\gamma$ (see (\ref{eq:eta-gamma})).
Furthermore, on applying the IBP formula (\ref{eq:etaIBP}) the
right-hand side of (\ref{eq:IBP-take1}) is represented in the form
\begin{align*}
-\int_{\varGamma_{X}}&\sum_{x\in\gamma}\left(\int_{{\frakX}^{\gamma}}
F(\gamma,\bar{y}^{\gamma})\,\beta^{\hat{v}}_{\eta}(x,\bar{y}_x)
\,\textstyle{\bigotimes\limits_{x'\in\gamma}}\eta_{x'}(\rd\bar{y}_{x'})\right) \mu(\rd\gamma) \\
&=-\int_{\varGamma_{X}}\left(\int_{{\frakX}^{\gamma}}
\sum_{x\in\gamma}
F(\gamma,\bar{y}^{\gamma})\,\beta^{\hat{v}}_{\eta}(x,\bar{y}_x)
\,\eta^{\gamma}(\rd\bar{y}^{\gamma})\right) \mu(\rd\gamma) \\
&=-\int_{\varGamma_{\mathcal{Z}}}F(\hatgamma)B^{\hat{v}}_{\hat{\mu}}(\hatgamma)
\,\hat{\mu}(\rd\hatgamma).
\end{align*}
which proves formula (\ref{eq:4.25}).

Finally, in view of Condition \ref{condition4}, Lemma \ref{moments}
implies that $\hat{\mu}\in
\mathcal{M}_{\hat{\sigma}}^{n}(\varGamma_{\mathcal{Z}})$, and by
Definition \ref{def:M} and Condition \ref{cond4.1} it follows that
$B_{\hat{\mu}}^{\hat{v}}\in
L^{1}(\varGamma_{\mathcal{Z}},\hat{\mu})$.
\endproof

The next two theorems are our main results in this section.

\begin{theorem}\label{IBP-}
For any function $F\in {\mathcal{FC}}(\varGamma_{X})$, the cluster
measure $\mucl$ satisfies the following IBP formula
\begin{equation}\label{IBP0-}
\int_{\varGamma_{X}}\sum_{x\in \gamma}\nabla_{\myn x}F(\gamma) \CDD
v(x)\,\mucl(\rd\gamma)=-\int_{\varGamma_{X}}F(\gamma)\mypp
B_{\mucl}^{v}(\gamma)\,\mucl(\rd\gamma),
\end{equation}
where $B_{\mucl}^{v}(\gamma):=\mathcal{I}_{\mathfrak{q}}^{\ast}
B^{\hat{v}}_{\hat{\mu}} \in L^{1}(\varGamma_{X},\mucl)$
\textup{(}see \textup{(\ref{eq:I*L1})}\textup{)} and the logarithmic
derivative $B^{\hat{v}}_{\hat{\mu}}(\hatgamma)$ is defined in
\textup{(\ref{eq:4.26})}.
\end{theorem}

\proof For any function $F\in {\mathcal{FC}}(\varGamma_{X})$ and
vector field $v\in \mathrm{Vect}_{0}(X)$, let us denote for brevity
\begin{equation}\label{eq:H}
H(x,\gamma):= \nabla_{\myn x}F(\gamma)\CDD v(x),\qquad x\in X,\ \
\gamma \in \varGamma_{X}.
\end{equation}
Furthermore, setting
$\hat{F}=\mathcal{I}_{\mathfrak{q}}F:\varGamma_{\mathcal{Z}}\rightarrow
\mathbb{R}$ we introduce the notation
\begin{equation}\label{eq:H-hat}
\hat{H}(z,\hatgamma):=\nabla_{\myn \bar{y}}\hat{F}(\hatgamma) \CDD
\hat{v}(\bar{y}),\qquad z=(x,\bar{y})\in \mathcal{Z},\ \
{\hatgamma}\in \varGamma_{\mathcal{Z}}.
\end{equation}
From these definitions, it is clear that
\begin{equation}\label{eq:QHH}
{\mathcal{I}_{\mathfrak{q}}}\Biggl(\mynn\sum_{x\in \gamma}
H(x,\gamma)\Biggr)(\hatgamma)=\sum_{z\in
\hatgamma}\hat{H}(z,\hatgamma),\qquad {\hatgamma}\in
\varGamma_{\mathcal{Z}}.
\end{equation}

By Theorem \ref{th:4.7}, the measure $\hat{\mu}$ satisfies the IBP
formula
\begin{equation}\label{IBPp}
\int_{\varGamma_{\mathcal{Z}}}\sum_{z\in
\hatgamma}\hat{H}(z,\hatgamma)\,\hat{\mu}(\rd\hatgamma)=-\int_{\varGamma_{\mathcal{Z}}}
\!\hat{F}(\hatgamma)B_{\hat{\mu}}^{\hat{v}}(\hatgamma)\,\hat{\mu}(\rd\hatgamma),
\end{equation}
where the logarithmic derivative
$B_{\hat{\mu}}^{\hat{v}}(\hatgamma)=\langle
\beta_{\eta}^{\hat{v}},\hatgamma\rangle$ belongs to
$L^{1}(\varGamma_{\mathcal{Z}},\hat{\mu})$ by Theorem \ref{th:4.7}.

Now, using formulae (\ref{eq:H-hat}), (\ref{eq:QHH}) and
(\ref{IBPp}), we obtain
\begin{align*}
\int_{\varGamma_{X}}\sum_{x\in \gamma}H(x,\gamma)\,\mucl(\rd\gamma)
&
=\int_{\varGamma_{\mathcal{Z}}}\!\Biggl(\mynn\sum_{(x,\myp\bar{y})\in
\hatgamma}\!\nabla_{\myn
\bar{y}}{\mathcal{I}_{\mathfrak{q}}}F(\hatgamma) \CDD
\hat{v}(\bar{y})\mynn\Biggr)\,\hat{\mu}(\rd\hatgamma) \\
&
=-\int_{\varGamma_{\mathcal{Z}}}{\mathcal{I}_{\mathfrak{q}}}F(\hatgamma)
\mypp B_{\hat{\mu}}^{\hat{v}}(\hatgamma)\,\hat{\mu}(\rd\hatgamma) \\
& =-\int_{\varGamma_{X}}F(\gamma)\,\mathcal{I}_{\mathfrak{q}}^{\ast}
B_{\hat{\mu}}^{\hat{v}}(\gamma)\,\mucl(\rd\gamma),
\end{align*}
where $\mathcal{I}_{\mathfrak{q}}^{\ast}B_{\hat{\mu}}^{\hat{v}}\in
L^{1}(\varGamma_{X},\mucl)$ by Lemma~\ref{lm:L1}. Thus, formula
(\ref{IBP0-}) is proved.
\endproof

\begin{remark}
Observe that the logarithmic derivative
$B_{\hat{\mu}}^{\hat{v}}=\langle
\beta_{\eta}^{\hat{v}},\hatgamma\rangle$ (see (\ref{IBPp})) does not
depend on the underlying measure $\mu$, and so it is the same as,
say, in the Poisson case with $\mu=\pi_{\theta}$. Nevertheless, the
logarithmic derivative $B_{\mucl}^{v}$ does depend on $\mu$ via the
mapping $\mathcal{I}_{\mathfrak{q}}^{\ast}$.
\end{remark}

According to Theorem \ref{IBP-}, $B_{\mucl}^{v}\!\in
L^{1}(\varGamma_{\mathcal{Z}},\mucl)$. However, under the conditions
of Lemma \ref{lm:4.6} with $n\geq 2$, this statement can be
enhanced.

\begin{lemma}\label{lm:B^n*}
Assume that\/ $\int_{\mathcal{Z}}\vert
\beta_\eta(z)\vert^{m}\,\hat{\sigma}(\rd z)<\infty $ for
$m=1,2,\dots,n$ and some integer\/ $n\geq 2$, and let condition\/
\textup{(\ref{eq:N<})} hold. Then\/ $B_{\mucl}^{v}\myn\in
L^{n}(\varGamma_{\mathcal{Z}},\mucl)$.
\end{lemma}

\proof By Lemmata \ref{lm:M^n} and \ref{lm:4.6}, it follows that
$\langle \beta_{\eta}^{\hat{v}},\hatgamma\rangle \in
L^{n}(\varGamma_{\mathcal{Z}},\hat{\mu})$. Let $r:=n/(n-1)$, so that
$n^{-1}+r^{-1}=1$. Note that $\mathcal{I}_{\mathfrak{q}}$ can be
treated as a bounded operator acting from
$L^{r}(\varGamma_{X},\mucl)$ to
$L^{r}(\varGamma_{\mathcal{Z}},\hat{\mu})$. Hence,
$\mathcal{I}_{\mathfrak{q}}^{\ast}$ is a bounded operator from
$L^{r}(\varGamma_{\mathcal{Z}},\hat{\mu})^{\prime
}=L^{n}(\varGamma_{\mathcal{Z}},\hat{\mu})$ to $L^{r}(\varGamma
_{X},\mucl)^{\prime}=L^{n}(\varGamma_{X},\mucl)$, which implies that
$B_{\mucl}^{v}=\mathcal{I}_{\mathfrak{q}}^{\ast}\langle \beta_{\eta
}^{\hat{v}},\hatgamma\rangle \in
L^{n}(\varGamma_{\mathcal{Z}},\mucl)$.
\endproof

Formula (\ref{IBP0-}) can be extended to more general vector fields
on $\varGamma_{X}$. Let ${\mathcal{FV}}(\varGamma_{X})$ be the class
of vector fields $V$ of the form $V(\gamma)=(V(\gamma)_{x})_{x\in
\gamma}$,
\begin{equation*}
V(\gamma)_{x}=\sum_{j=1}^{k}G_{j}(\gamma)\,v_{j}(x)\in T_{x}X,
\end{equation*}
where $G_{j}\in {\mathcal{FC}}(\varGamma_{X})$ and
$v_{j}\in\Vect_{0}(X)$, \thinspace $j=1,\dots,k$. For any such $V$
we set
\begin{equation*}
B_{\mucl}^{V}(\gamma):=(\mathcal{I}_{\mathfrak{q}}^{\ast
}B_{\hat{\mu}}^{{\mathcal{I}_{\mathfrak{q}}}V})(\gamma),
\end{equation*}
where $B_{\hat{\mu}}^{\mathcal{I}_{\mathfrak{q}} V}(\hatgamma)$ is
the logarithmic derivative of $\hat{\mu}$ along
$\mathcal{I}_{\mathfrak{q}}
V(\hatgamma):=V(\mathfrak{q}(\hatgamma))$ (see \cite{AKR1}). Note
that ${\mathcal{I}_{\mathfrak{q}}}V$ is a vector field on
$\varGamma_{\mathcal{Z}}$ owing to the obvious equality
\begin{equation*}
T_{\hatgamma}\varGamma_{\mathcal{Z}}
=T_{{\mathfrak{q}}(\hatgamma)}\varGamma_{X}.
\end{equation*}
Clearly,
\begin{equation*}
B_{\mucl}^{V}(\gamma)=\sum_{j=1}^{k}\biggl(G_{j}(\gamma)
B_{\mucl}^{v_{j}}(\gamma)+\sum_{x\in \gamma}\nabla_{\myn
x}G_{j}(\gamma) \CDD v_{j}(x)\biggr).
\end{equation*}

\begin{theorem}
\label{IBP1} For any $F_{1},F_{2}\in {\mathcal{FC}}(\varGamma_{X})$
and\/ $V\in {\mathcal{FV}}(\varGamma_{X})$, we have
\begin{equation*}
\begin{aligned} \int_{\varGamma_{X}}&\sum_{x\in \gamma}\nabla_{\myn
x}F_{1}(\gamma)\CDD V(\gamma)_{x}\,F_{2}(\gamma)\;\mucl(\rd\gamma) \\
&=-\int_{\varGamma_{X}}F_{1}(\gamma)\,\sum_{x\in \gamma}\nabla_{\myn
x}F_{2}(\gamma)\CDD  V(\gamma)_{x}\ \mucl(\rd\gamma)
-\int_{\varGamma_{X}}F_{1}(\gamma)
F_{2}(\gamma)B_{\myn\mucl}^{V}(\gamma)\,\mucl(\rd\gamma).
\end{aligned}
\end{equation*}
\end{theorem}

\proof The proof can be obtained by a straightforward generalisation
of the arguments used in the proof of Theorem \ref{IBP-}. \endproof

We define the \emph{vector logarithmic derivative} of $\mucl$ as a
linear operator
\begin{equation*}
B_{\mucl}\!:\,{\mathcal{FV}}(\varGamma_{X})\to
L^{1}(\varGamma_{X},\mucl)
\end{equation*}
via the formula
\begin{equation*}
B_{\mucl}V(\gamma):=B_{\mucl}^{V}(\gamma).
\end{equation*}
This notation will be used in the next section.

\subsection{Dirichlet forms and equilibrium stochastic dynamics}
\label{sec:4.4}

Throughout this section, we assume that the conditions of Lemma
\ref{lm:M^n} are satisfied with $n=2$. Thus, the measures $\mu$,
$\hat{\mu}$ and $\mucl$ belong to the corresponding
$\mathcal{M}^{2}$-classes. Our considerations will involve the
$\varGamma $-gradients (see Section~\ref{sec:4.1}) on different
configuration spaces, such as $\varGamma_{X}$,
$\varGamma_{\mathfrak{X}}$ and $\varGamma_{\mathcal{Z}}$; to avoid
confusion, we shall denote them by $\nabla_{\myn X}^{\varGamma}$,
$\nabla_{\myn \mathfrak{X}}^{\varGamma}$ and $\nabla_{\mynn
\mathcal{Z}}^{\varGamma}$, respectively.

Let us introduce a pre-Dirichlet form ${\mathcal{E}}_{\mucl}$
associated with the Gibbs cluster measure $\mucl$, defined on
functions $F_{1},F_{2}\in \mathcal{FC}(\varGamma_{X})\subset
L^{2}(\varGamma_{X},\mucl)$ by
\begin{equation}\label{eq:E-mu}
{\mathcal{E}}_{\mucl}(F_{1},F_{2}):=\int_{\varGamma_{X}}\langle
\nabla_{\myn X}^{\varGamma}F_{1}(\gamma),\nabla_{\myn X}^{\varGamma}
F_{2}(\gamma)\rangle_{\gamma}\,\mucl(\rd\gamma).
\end{equation}
Let us also consider the operator $H_{\mucl}$ defined by
\begin{equation}
H_{\mucl}F:=-\Delta^{\varGamma}F+B_{\mucl}\!\nabla_{\myn
X}^{\varGamma }F,\qquad F\in {\mathcal{FC}}(\varGamma_{X}),
\label{eq:Hcl}
\end{equation}
where $\Delta^{\varGamma}F(\gamma):=\sum_{x\in
\gamma}\Delta_{x}F(\gamma)$.

The next theorem readily follows from the general theory of
(pre-)Dirichlet forms associated with measures from the class
$\mathcal{M}^{2}(\varGamma_{X})$ (see \cite{AKR2,MR}).

\begin{theorem}
The pre-Dirichlet form \textup{(\ref{eq:E-mu})} is well defined,
i.e., ${\mathcal{E}}_{\mucl}(F_{1},F_{2})<\infty $ for all
$F_{1},F_{2}\in {\mathcal{FC}}(\varGamma_{X})$. Furthermore,
expression \textup{(\ref{eq:Hcl})} defines a symmetric operator
$H_{\mucl}$ in $L^{2}(\varGamma_{X},\mucl)$, which is the generator
of\/ ${\mathcal{E}}_{\mucl}$, that is,
\begin{equation}\label{generator}
{\mathcal{E}}_{\mucl}(F_{1},F_{2})=\int_{\varGamma_{X}}F_{1}(\gamma)
\,H_{\mucl}F_{2}(\gamma)\,\mucl(\rd\gamma),\ \quad F_{1},F_{2}\in
{\mathcal{FC}}(\varGamma_{X}).
\end{equation}
\end{theorem}

Formula (\ref{generator}) implies that the form
${\mathcal{E}}_{\mucl}$ is closable. It follows from the properties
of the \textit{carr\'{e} du champ} \,$\sum_{x\in\gamma}
\mynn\nabla_{\myn x}F_{1}(\gamma) \CDD \nabla_{\myn x}F_{2}(\gamma)$
that the closure of ${\mathcal{E}}_{\mucl}$ (for which we shall keep
the same notation) is a quasi-regular local Dirichlet form on a
bigger state space $\overset{\,..}{\varGamma}_{X}$ consisting of all
integer-valued Radon measures on $X$ (see \cite{MR}). By the general
theory of Dirichlet forms (see \cite{MR0}), this implies the
following result (cf.\ \cite{AKR1,AKR2,BD3}).

\begin{theorem}
\label{th:7.2} There exists a conservative diffusion process
$\mathbf{X}=(\mathbf{X}_{t},\,t\geq 0)$ on
$\overset{\,\myp..}{\varGamma}_{X}$, properly associated with the
Dirichlet form $\mathcal{E}_{\mucl}$, that is, for any function
$F\in L^{2}(\overset{\,\myp..}{\varGamma}_{X},\mucl)$ and all\/
$t\geq 0$, the map
\begin{equation*}
\overset{\,\myp..}{\varGamma}_{X}\ni \gamma \mapsto p_{t}F(\gamma
):=\int_{\varOmega}F(\mathbf{X}_{t})\,\rd P_{\gamma}
\end{equation*}
is an\/ $\mathcal{E}_{\mucl}$-quasi-continuous version of\/ $\exp
(-tH_{\mucl})F$. Here $\varOmega$ is the canonical sample space
\textup{(}of $\overset{\,\myp..}{\varGamma}_{X}$-valued continuous
functions on $\RR_{+}$\textup{)} and $(P_{\gamma},\,\gamma \in
\overset{\,\myp..}{\varGamma}_{X})$ is the family of probability
distributions of the process $\mathbf{X}$ conditioned on the initial
value $\gamma =\mathbf{X}_{0}$. The process $\mathbf{X}$ is unique
up to $\mucl$-equivalence. In particular, $\mathbf{X}$ is
$\mucl$-symmetric \textup{(}i.e., $\int F_{1}\mypp
p_{t}F_{2}\,\rd\mucl=\int F_{2}\,p_{t}F_{1}\,\rd\mucl$ for all
measurable functions $ F_{1},F_{2}:\overset{\,\myp..}{\varGamma
}_{X}\rightarrow \mathbb{R}_{+}$\textup{)} and $\mucl$ is its
invariant measure.
\end{theorem}

\subsection{On the irreducibility of the Dirichlet form}

Let ${\mathcal{E}}_{\hat{\mu}}$ be the pre-Dirichlet form associated
with the Gibbs measure $\hat{\mu}$, defined on functions
$F_{1},F_{2}\in
\mathcal{FC}_{\hat{\sigma}}(\varGamma_{\mathcal{Z}})\subset
L^{2}(\varGamma_{\mathcal{Z}},\hat{\mu})$ by
\begin{equation}\label{generator1}
{\mathcal{E}}_{\hat{\mu}}(F_{1},F_{2}):=\int_{\varGamma_{\mathcal{Z}
}}\langle \nabla_{\mynn\mathcal{Z}}^{\varGamma}F_{1}(\hatgamma),
\nabla_{\mynn\mathcal{Z}}^{\varGamma}F_{2}(\hatgamma)\rangle
_{\hatgamma}\,\hat{\mu}(\rd\hatgamma).
\end{equation}
The integral on the right-hand side of (\ref{generator1}) is well
defined because $\hat{\mu}\in \mathcal{M}_{\theta}^{2}
(\varGamma_{\mathcal{Z}})\subset \mathcal{M}_{\theta}^{1}
(\varGamma_{\mathcal{Z}})$. Indeed, the function
\begin{equation*}
G(z):=\left( \nabla_{\mynn z} F_{1}(\hatgamma),\nabla_{\mynn z}
F_{2}(\hatgamma)\right)
\end{equation*}
is bounded and has a $\hat{\sigma}$-finite support, which implies
that $G\in L^{1}(\mathcal{Z},\hat{\sigma})$. Thus the function
\begin{equation*}
\varGamma_{\mathcal{Z}}\ni \hatgamma\mapsto \left\langle G,\gamma
\right\rangle \equiv \langle
\nabla_{\mynn\mathcal{Z}}^{\varGamma}F_{1}(\hatgamma),
\nabla_{\mynn\mathcal{Z}}^{\varGamma} F_{2}
(\hatgamma)\rangle_{\hatgamma}
\end{equation*}
belongs to $L^{2}(\varGamma_{\mathcal{Z}},\hat{\mu})$ by the
definition of the class $\mathcal{M}_{\theta}^{1}
(\varGamma_{\mathcal{Z}})$. It can be shown by a direct computation
that
\begin{equation}\label{dirform-ident}
\mathcal{E}_{\mucl}(F,F)=\mathcal{E}_{\hat{\mu}}
(\mathcal{I}_{\mathfrak{q}}F,\mathcal{I}_{\mathfrak{q}}F),\qquad
F\in {\mathcal{FC}}(\varGamma_{X}).
\end{equation}

Note that the pre-Dirichlet form
$(\mathcal{E}_{\hat{\mu}},{\mathcal{FC}}_{\hat{\sigma}}(\varGamma_{\mathcal{Z}}))$
is not necessarily closable.
A sufficient condition of its closability is an IBP formula for the
measure $\hat{\mu}$ with respect to \textit{all} directions in
$\varGamma_{Z}$ rather then only in $\mathfrak{X}^{\gamma}$ (cf.\
Theorem \ref{th:4.7}), which requires in turn some smoothness
conditions on the measure $\mu$ and also on the measure $\eta_{x}$
as a function of $x\in X$. Such conditions are satisfied, for
instance, if $X=\RR^{d}$, $\mu $ is a Poisson measure or, more
generally, Gibbs measure with a smooth interaction potential, and
the family $\{\eta_{x}\}$ is defined by translations of a parent
measure $\eta_{0}$ (i.e., $\eta_{x}(B):=\eta_{0}(B-x)$).
This case has been studied in great detail in \cite{BD3,BD4}, where
formula (\ref{dirform-ident}) was extended to all functions $F$ from
the domain of $\mathcal{E}_{\mucl}$ (with the closure
$\bar{\mathcal{E}}_{\hat{\mu}}$ of the pre-Dirichlet form
$(\mathcal{E}_{\hat{\mu}},{\mathcal{FC}}_{\hat{\sigma}}(\varGamma_{\mathcal{Z}}))$
on the right-hand side). In turn, this makes it possible to
characterise the kernel of the Dirichlet form $\mathcal{E}_{\mucl}$
via the kernels of the forms $\bar{\mathcal{E}}_{\hat{\mu}}$ and
$\mathcal{E}_{\mu }$; in particular, it has been proved in
\cite{BD3,BD4} that $\mathcal{E}_{\mucl}$ is \textit{irreducible}
(that is, its kernel consists of constants) whenever
$\mathcal{E}_{\mu}$ is such.

Let us remark that irreducibility is an important property closely
related to the ergodicity of stochastic dynamics and extremality of
invariant measures. It seems plausible that in our situation the
irreducibility of $\mathcal{E}_{\mucl}$ is controlled by the
properties of the distribution of centres $\mu$ rather then the
cluster distributions $\{\eta_{x}\}$, but this remains an open
question.

\section{Examples}\label{sec:5}

In order to make tractable the general cluster model discussed
above, one needs an efficient method to construct the family
$\{\eta_{x}\}_{x\in X}$ of cluster distributions attached to centres
$x$ lying on a ground configuration $\gamma$. In the situation where
$X$ is a linear space, this is straightforward by translations of a
parent distribution $\eta_{0}$ specified at the origin (see Section
\ref{sec:5.1}). For other classes of spaces, the linear action has
to be replaced by another suitable transformation (see Sections
\ref{sec:5.2}, \ref{sec:5.3} and \ref{sec:5.4}). Alternative, more
direct methods may also be applicable based on specific properties
of the space structure, for instance by confining oneself to a class
of distributions with a suitable invariance property (cf.\ Section
\ref{sec:5.3}) or by exploiting the space metric, leading to
``radially symmetric'' distributions (see Section \ref{sec:5.5}).

We discuss below a number of selected examples where this programme
can be realised. In so doing, we will mostly be using the
push-forward method of Section \ref{sec:droplet}. Specifically, the
discussion of the resulting cluster measure $\mucl$ in each example
will be essentially confined to the following two important aspects:
\begin{itemize}
\item[(i)]
verification of general sufficient conditions for the cluster
process configurations to be proper, such as Condition
\ref{condition5}\mypp$'$ in Proposition \ref{pr:3.9} specialised to
the conditions of Propositions \ref{pr:3.9} and \ref{pr:3.10} (local
finiteness), and the conditions of Propositions \ref{pr:simple} and
their particular cases in Propositions \ref{pr:a2-1} and
\ref{pr:a2-2} (simplicity); and
\item[(ii)]
verification of appropriate smoothness conditions on the mapping
$\varphi_x$ that we imposed as a prerequisite of an IBP formula for
the cluster measure $\mucl$ (see the beginning of Section
\ref{sec:4.3.2}).
\end{itemize}

\subsection{Euclidean spaces}\label{sec:5.1}

In the situation where $X=\RR^{d}$, the family $\{\eta _{x}\}_{x\in
X}$ of cluster distributions can be constructed by translations of a
parent distribution $\eta _{0}$ specified at the origin
\cite{BD3,BD4}. This can be formulated in terms of the construction
of Section \ref{sec:droplet}. Take $W:=X$ and define the family of
maps $\varphi_{x}:X\rightarrow X$ ($x\in X$) as translations
\begin{equation}\label{eq:phi-vector}
\varphi _{x}(y):=y+x,\qquad y\in X.
\end{equation}
Then definition (\ref{eq:D}) of the droplet $D_{\myn B}(y)$
specialises to
\begin{equation*}
D_{\myn B}(y)=B-y,\qquad y\in X, \ \ B\in \mathcal{B}(X).
\end{equation*}
Furthermore, formula (\ref{eq:D-union}) for the droplet cluster now
reads
$$
\bar{D}_B(\bar{y})=
{\textstyle\bigcup\limits_{y_i\in\bar{y}}}(B-y_i),\qquad
\bar{y}\in\mathfrak{X},
$$
which makes the notion of the droplet cluster particularly
transparent as a set-theoretic union of ``droplets'' of shape $B$
shifted to the centrally reflected coordinates of the vector
$\bar{y}=(y_i)$. The parent measure $Q$ on $\mathfrak{X}$ (see
(\ref{eq:x0->x})) can then be interpreted as a pattern distribution
$\eta_0$, and the measures $\eta_x$ are obtained by translations of
$\eta_0$ to points $x\in X$:
\begin{equation}\label{eq:eta0-vector}
\eta_x(\bar{B}):=\bar{\varphi}_x^*\eta_0(\bar{B})\equiv
\eta_0(\bar{B}-x),\qquad \bar{B}\in\mathcal{B}(\mathfrak{X}).
\end{equation}

Let us discuss in this context criteria of properness of the
corresponding cluster measure $\mucl$ laid out in Section
\ref{sec:droplet}. First of all, conditions (\ref{eq:C_B}) and
(\ref{eq:EN<}) of Proposition \ref{pr:3.9} (which guarantee
Condition \ref{condition5}\mypp$'$) are reduced, respectively, to
\begin{gather}\label{co1-R}
\sup_{y\in X}\theta (B-y)<\infty,\qquad
\int_{\mathfrak{X}}N_{X}(\bar{y})\,\eta_0(\rd\bar{y})<\infty .
\end{gather}
In turn, the first condition in (\ref{co1-R}) is satisfied, for
instance, if the measure $\theta(\rd x)$ is absolutely continuous
with respect to Lebesgue measure $\rd{x}$ on $X$ and the
corresponding Radon--Nikodym density is bounded (cf.\ Remark
\ref{rm:3.8}). Next, condition (i) of Proposition \ref{pr:3.10}
(i.e., continuity of $\varphi_x$ in $x$ is obviously satisfied for
(\ref{eq:phi-vector}), while condition (ii) holds with a compact
$B_y=B-y$ \,($y\in X$). Finally, let us point out that the use of
Propositions \ref{pr:a2-1} and \ref{pr:a2-2} is greatly facilitated
by the fact that the equation $\varphi_y(w)=x$, reducing for
(\ref{eq:phi-vector}) to equation $w+y=x$, has the unique solution
$y=x-w$.

Regarding conditions for IBP formulae, note that map
(\ref{eq:phi-vector}) is of course smooth, with $\rd\varphi_x=\Id$
(the identity operator) and $\rd^2\varphi_x=0$. Finally, if the
probability measure $\eta_0(\rd\bar{x})$ is absolutely continuous
with respect to Lebesgue measure $\rd\bar{x}$ on $\mathfrak{X}$,
then conditions (\ref{beta}) and (\ref{N}) in Theorem \ref{th:C2}
can be easily rewritten in terms of the corresponding Radon--Nikodym
density.

\subsection{Lie groups}\label{sec:5.2}

Let $X=G$ be a (non-compact) Lie group, and $\mathfrak{g}$ the
corresponding Lie algebra endowed with a scalar product
$(\cdot,\cdot)_{\mathfrak{g}}$ (see, e.g., \cite{Helgason}). This
scalar product generates in a standard way a right-invariant
Riemannian structure on $G$. The group product of elements $g_{1},
g_{2}\in G$ is denoted by $g_{1}g_{2}\in G$, and $e\in G$ stands for
the identity of the group $G$.

Let us show how a family of measures $\{\eta_{x}\}_{x\in G}$ on
$\mathfrak{G}:=\bigsqcup_{n=0}^\infty G^n$ can be set out using the
push-forward construction of Section \ref{sec:droplet}. Take $W:=G$
and define the map $\varphi_{x}(g):G\times G\rightarrow G$ as a
translation
\begin{equation}\label{eq:gx}
\varphi_{x}(g):=gx,\qquad g,x\in G.
\end{equation}
By the properties of the Lie group multiplication, the map
$\varphi_{x}(g)$ is continuous in $(g,x)\in G\times G$ and therefore
automatically measurable. In view of (\ref{eq:gx}), definition
(\ref{eq:D}) of the droplet $D_{\myn B}(g)$ specialises to
\begin{equation*}
D_{\myn B}(g)=g^{-1}B,\qquad B\in\mathcal{B}(G),\ \ g\in G.
\end{equation*}
Accordingly, by formula (\ref{eq:D-union}) the corresponding droplet
cluster is represented as
\begin{equation*}
\bar{D}_{\myn B}(\bar{g})={\textstyle\bigcup\limits_{g_i\in
\bar{g}}} \,g_i^{-1}B,\qquad
\bar{g}\in\mathfrak{G}:={\textstyle\bigsqcup\limits_{n=0}^\infty}
G^{n}.
\end{equation*}
If $Q$ is a probability measure on $\mathfrak{G}$, then on
substituting (\ref{eq:gx}) into definition (\ref{eq:x0->x}) we get
\begin{equation}\label{eq:eta_x-group}
\eta_{x}(\bar{B}):=(\bar{\varphi}_{x}^{\ast}\myp Q)(\bar{B})\equiv
Q(\bar{B}\myp x^{-1}),\qquad
\bar{B}\in\mathcal{B}(\mathfrak{G})\quad (x\in G).
\end{equation}
Observe from (\ref{eq:eta_x-group}) that in fact the measure $Q$
coincides with $\eta_e$; hence definition (\ref{eq:eta_x-group}) can
be rewritten in a ``translation'' form naturally generalising
formula (\ref{eq:eta0-vector}) in the Euclidean case, namely
\begin{equation}\label{eq:eta_0-group'}
\eta_{x}(\bar{B})=\eta_e(\bar{B}\myp x^{-1}),\qquad
\bar{B}\in\mathcal{B}(\mathfrak{G})\quad (x\in G).
\end{equation}

Specialising the general criteria of properness of $\mucl$ described
in Section \ref{sec:droplet}, we have that conditions (\ref{eq:C_B})
and (\ref{eq:EN<}) of Proposition \ref{pr:3.9} are reduced,
respectively, to
\begin{gather}\label{co1-gr}
\sup_{g\in G}\theta (g^{-1}B)<\infty,\qquad
\int_{\mathfrak{G}}N_{G}(\bar{g})\,Q(\rd\bar{g})<\infty.
\end{gather}
Similarly to the previous section, the first condition in
(\ref{co1-gr}) is satisfied proviso the reference measure $\theta$
is absolutely continuous with respect to a left Haar measure on $G$
and the corresponding Radon--Nikodym density is bounded (cf.\ Remark
\ref{rm:3.8}). As mentioned above, maps (\ref{eq:gx}) automatically
satisfy the continuity condition (i) of Proposition \ref{pr:3.10},
whereas condition (ii) holds with a compact $B_g=B\myp g^{-1}$
($g\in G$). Moreover, as a natural extension of the Euclidean case,
the equation $\varphi_y(g)=x$ with (\ref{eq:gx}) takes the form
$g\myp y=x$, which has the unique solution $y=g^{-1}x$. Hence,
Propositions \ref{pr:a2-1} and \ref{pr:a2-2} can be easily applied.

In a standard fashion, the Lie algebra $\mathfrak{g}$ of the group
$G$ can be identified with the space of right-invariant vector
fields on $G$; moreover, all tangent spaces $T_{g}G $ are identified
with $T_{e}G$ (and therefore with $\mathfrak{g}$) via right
translations. Under this identification, for the map $\varphi_x(w)$
defined in (\ref{eq:gx}) we have $\rd\varphi_{x}(g)=\Id$ for any
$x,g\in G$, where $\Id:\mathfrak{g}\to\mathfrak{g}$ is the identity
operator. It follows that $\|\rd\varphi_{x}(g)\|=1$ and
$\rd^{2}\varphi _{x}(g)=0$ for all $x,g\in G$, which automatically
implies that $\varphi _{x}\in C_{b}^{2}(G,G)$ uniformly in $x\in G$.
Thus, one can apply Theorem \ref{th:C2} 
provided that conditions (\ref{beta}) and (\ref{N}) are satisfied.
Finally, if the probability measure $Q$ is absolutely continuous
with respect to a left Haar measure on $\mathfrak{G}$, then
conditions (\ref{beta}) and (\ref{N}) can be easily specified in
terms of the corresponding Radon--Nikodym density.

\subsection{Homogeneous manifolds}\label{sec:5.3}

\subsubsection{Construction of cluster distributions $\eta_x$.}
Let $G$ be a (non-compact) Lie group and $X$ a $G$-homogeneous
Riemannian manifold (see, e.g., \cite{Besse,Helgason}). More
precisely, $G$ is a closed subgroup of the group of isometries of
$X$ acting on $X$ transitively, that is, for any $x,y\in X$ there
exists an element $g\in G$ such that $g\cdot x = y$ (equivalently,
$G\cdot x=X$ for some, and hence for all $x\in X$), and the mapping
\begin{equation}\label{eq:g(x)}
G\times X\ni (g,x)\mapsto g\cdot x\in X
\end{equation}
is differentiable. Given a fixed point $x_{0}\in X$, the manifold
$X$ is diffeomorphic to the quotient manifold $G/H_{x_{0}}$, where
$H_{x_0}:=\{g\in G: g\cdot x_0=x_0\}$ is the isotropy subgroup of
$G$ at $x_0$.

\begin{example}\label{ex:5.1}
Take $X=\RR^{d}$ and the group $G=\mathbb{R}^{d}$ with the natural
additive structure acting on $X$ by translations. In this case,
$H_{x_{0}}=\left\{ 0\right\} $ for every $x_{0}\in X$.
\end{example}

\begin{example}\label{ex:5.2}
Let $X=\RR^{d}$ and consider $G=\Euc^+(d)$, the Euclidean group of
isometries of $\RR^{d}$ preserving orientation. In this case,
$H_{0}=\SO(d)$ and $X\cong \Euc^+(d)/\SO(d)$.
\end{example}

\begin{example}\label{ex:5.3}
Let $X=\mathbb{H}^{d}$ be a $d$-dimensional hyperbolic space. In
this situation, $G=\SO_{0}(d,1)$ is the connected component of the
identity in the orthogonal group $\Ort(d,1)$ of the canonical
quadratic form with signature $(n,1)$, and
$X\cong\SO_{0}(d,1)/\SO(d)$.
\end{example}

\begin{example}\label{ex:5.4}
If $G$ is a Lie group and $H$ is its compact subgroup, then one can
use the quotient manifold $X=G/H$ with the natural $G$-action on it.
\end{example}

Define a family of maps $\varphi_x:G\to X$ as the group action (see
(\ref{eq:g(x)}))
\begin{equation}\label{eq:phi_x}
\varphi_{x}(g):=g\cdot x,\qquad g\in G,\ \ x\in X. \end{equation}
Then definition (\ref{eq:D}) of the droplet $D_{\myn B}(g)$ takes
the form
\begin{equation*}
D_{\myn B}(g)= g^{-1}\myn\cdot B,\qquad B\in \mathcal{B}(X), \ \
g\in G,
\end{equation*}
and the droplet cluster is given by
\begin{equation*}
\bar{D}_{\myn B}(\bar{g})={\textstyle\bigcup\limits_{g_i\in\bar{g}}}
(g_i^{-1}\myn\cdot B),\qquad
\bar{g}\in\mathfrak{G}:=\textstyle{\bigsqcup\limits_{n=0}^\infty}G^{n}.
\end{equation*}
According to Section \ref{sec:droplet}, we can now use
(\ref{eq:phi_x}) to define the family of distributions
\begin{equation}\label{eta-ex1}
\eta_{x}:=\bar{\varphi}_{x}^{\ast}\myp Q\equiv Q\circ
\bar{\varphi}_{x}^{-1},
\end{equation}
where $Q$ is a given probability measure on $\mathfrak{G}$.

Conditions (\ref{eq:C_B}) and (\ref{eq:EN<}) take the form,
respectively,
\begin{equation}\label{co1}
\sup_{g\in G}\theta (g^{-1}\myn\cdot B)<\infty,\qquad
\int_{\mathfrak{G}}N_{G}(\bar{g})\,Q(\rd\bar{g})<\infty.
\end{equation}
The first of conditions (\ref{co1}) is satisfied, for instance, if
$\theta$ is absolutely continuous with respect to the volume measure
on $X$ and the corresponding Radon--Nikodym density is bounded (cf.\
Remark \ref{rm:3.8}).

Let us point out that a special class of measures
$\{\eta_{x}\}_{x\in X}$ on $\mathfrak{G}$ can be constructed
somewhat more naturally by essentially reproducing the group
translations method for Lie groups (cf.\ (\ref{eq:eta_x-group})).
More precisely, fix an arbitrary point $x_{0}\in X$ and an
$H_{x_{0}}$-invariant measure $\eta_{x_0}$ on $\mathfrak{X}$ (i.e.,
$\eta_{x_0}(h\bar{B})=\eta_{x_0}(\bar{B})$ for any
$\bar{B}\in\mathcal{B}(\frakX)$ and all $h\in H_{x_0}$); such a
measure always exists due to the compactness of $H_{x_{0}}$. Since
the group action is transitive, the group orbit of $x_0$ coincides
with $X$, hence each $x\in X$ can be represented in the form
$x=g\cdot x_0$ with some $g=g_x\in G$. Let us now define the measure
$\eta_{x}$ on $\calB(\frakX)$ by the formula
\begin{equation}\label{eta-ex1-direct}
\eta_{x}:=g_x^*\mypp\eta_{x_0}\equiv \eta_{x_0}\circ g_x^{-1},\qquad
x=g_x\cdot x_{0}.
\end{equation}
It follows that $\eta_x$ is $H_x$-invariant for each $x\in X$.
Definition (\ref{eta-ex1-direct}) does not depend on the choice of a
solution $g_x$ of the equation $g\cdot x_0=x$; indeed, if there is
another solution $\tilde{g}_x$ then
$$
\eta_{x_0}\circ \tilde{g}_x^{-1}=\eta_{x_0}\circ
(\tilde{g}_x^{-1}g_x)\,g_x^{-1}=\eta_{x_0}\circ g_x^{-1},
$$
since $\tilde{g}_x^{-1}g_x\in H_{x_0}$ and $\eta_{x_0}$ is
$H_{x_0}$-invariant.

\begin{remark}
Choosing various subgroups $G$ of the general group of isometries of
$X$ may lead to different representations of $X$ as a homogeneous
space. Consequently, formula (\ref{eta-ex1}) will define different
cluster measures. This is illustrated in the next simple example for
the Euclidean space.
\end{remark}

\begin{example}\label{ex:5.5}
Let $X=\RR^d$ \,($d\ge2$). If $G$ is the group of
\textit{translations} $x\mapsto x-g$ $(x,g\in \RR^d$), then the
corresponding homogeneous space is isomorphic to $\RR^d$ and, as
described in Section \ref{sec:5.1}, the measures $\eta_x$ are
obtained by translations, $\eta_x(\cdot)=\eta_0(\cdot-x)$ (see
Example \ref{ex:5.1}). Let now $G=\Euc^+(d)$ (see Example
\ref{ex:5.2}), that is, the group of \textit{rotations} $g=(\xi,A)$
with the action $\varphi_x(g):=g\cdot x = A(x-\xi) + \xi$
\,($x\in\RR^d$), where $\xi\in\RR^d$
and
$A\in\SO(d)$.
It is easy to check that, for a given Borel set $B\subset\RR^d$,
\begin{align*}
\varphi_x^{-1}(B\setminus\{x\})&=\bigl\{g\in
G:\,A\ne I\ \
\text{and}\ \ \xi\in (I-A)^{-1}(B-Ax)\bigr\},\\
\varphi_x^{-1}(\{x\})&=\bigl\{g\in G:\,A\ne I,\ \xi=x\ \
\text{or}\ \ A=I,\ \xi\in\RR^d\bigr\},
\end{align*}
where $I$ is the identity matrix. Consider the simplest case where
each cluster contains only one point; in other words, the measures
$\eta_x$ are supported on $X$ (i.e., $\eta_x(X^n)=0$ for $n\ne 1$).
Let $Q(\rd{}\xi\times\rd{A})$ be a probability measure on $G$;
assume for simplicity that $Q\{A\ne I\}=1$.
Then definition (\ref{eta-ex1}) specialises to
\begin{equation}\label{eq:g0|alpha}
\eta_x(B)=Q(\varphi_x^{-1}(B))=\int_{\SO(d)}
Q(\RR^d\times\rd A) \int_{(I-A)^{-1}(B-Ax)}\,Q(\rd{}\xi\myp|\myp A),
\end{equation}
where $Q(\RR^d\times\rd A)$ is the marginal distribution of $A$ and
$Q(\rd{}\xi\myp|\myp A)$ is the conditional distribution of $\xi$
given $A$.
Conditionally on $A$,
$\eta_x$ is obtained from $\eta_0$ via a translation by the vector
$-(I-A)^{-1}Ax$, which is different from $x$. If $A$ is truly
random, then averaging with respect to its distribution will further
mix up the random shifts $-(I-A)^{-1}Ax$.
\end{example}

\subsubsection{Verification of smoothness.}
Our next goal is to show that $\varphi_{x}(\cdot)\in C_{b}^{2}(G,X)$
uniformly in $x\in G$ for a special Riemannian metric on $G$.
Following \cite[Ch.~7, pp.\ 181--186]{Besse}, fix any $x\in X$ and
let, as before, $H_{x}$ be the isotropy subgroup at $x$. Then the
manifold $X$ can be identified with the quotient manifold $G/H_{x}$
in such a way that the map $\varphi_{x}:G\to X$ coincides with the
natural projection $G\rightarrow G/H_{x}$. Let $\mathfrak{h}_{x}$ be
the Lie algebra of $H_{x}$. It is known that the Lie algebra
$\mathfrak{g}$ of the Lie group $G$ admits a decomposition
\begin{equation}\label{decomp}
\mathfrak{g}=\mathfrak{h}_{x}\oplus\mathfrak{x}_{x},
\end{equation}
where $\mathfrak{x}_{x}$ is a subspace of $\mathfrak{g}$ invariant
with respect to the adjoint representation $H_x\ni h\mapsto \Ad_{h}$
of $H_x$ in $\mathfrak{g}$. Then the tangent space $T_{x}X$ can be
identified with the space $\mathfrak{x}_{x}$. The Riemannian metric
of $X$ induces an $\Ad_{h}$-invariant scalar product
$(\cdot,\cdot)_{\mathfrak{x}_{x}}$ on $\mathfrak{x}_{x}$.

Let us choose an auxiliary $\Ad_{h}$-invariant scalar product
$(\cdot,\cdot)_{\mathfrak{h}_{x}}$ on $\mathfrak{h}_{x}$. Such a
product always exists thanks to the compactness of $H_{x}$; for
instance, we can set
$(\cdot,\cdot)_{\mathfrak{h}_{x}}:=-B(\cdot,\cdot)$, where $B$ is
the Killing--Cartan form (see, e.g., \cite[Ch.~7, pp.\
184--185]{Besse} or \cite[Ch.~II, \S\myp6, p.~131]{Helgason}).
Observe that the isotropy subgroup at $g\cdot x$ has the form
$H_{g\cdot x}=g\cdot H_{x}\myp g^{-1}$, therefore the corresponding
Lie algebra is given by $\mathfrak{h}_{g\cdot
x}=\Ad_{g}(\mathfrak{h}_{x})$. We equip it with the scalar product
\begin{equation*}
(\cdot,\cdot)_{\mathfrak{h}_{g\cdot x}}:=(\Ad_{g^{-1}}\cdot,
\Ad_{g^{-1}}\cdot)_{\mathfrak{h}_{x}}.
\end{equation*}
Moreover, we can set $\mathfrak{x}_{g\cdot
x}=\Ad_{g}(\mathfrak{x}_{x})$, so that decomposition (\ref{decomp})
at $g\cdot x$ takes the form
\begin{equation}\label{decomp1}
\mathfrak{g}=\Ad_{g}(\mathfrak{h}_{x})\oplus
\Ad_{g}(\mathfrak{x}_{x}).
\end{equation}
Now we can define a scalar product $(\cdot,\cdot)_{\mathfrak{g,}g}$
on $\mathfrak{g}$ by setting for all $h\in\mathfrak{h}_{g\cdot x}$,
\,$r\in \mathfrak{x}_{g\cdot x}$
\begin{equation*}
(h + r,h + r)_{\mathfrak{g},g}=(h,h)_{\mathfrak{h}_{g\cdot x}}
+(r,r)_{\mathfrak{x}_{g\cdot x}}.
\end{equation*}
The $G$-invariance of the Riemannian metric on $X$ implies that
\begin{equation}\label{eq:5.15}
(\cdot,\cdot)
_{\mathfrak{g},g}=(\Ad_{g^{-1}}\cdot,\Ad_{g^{-1}}\cdot)
_{\mathfrak{g}, e}.
\end{equation}
The family of scalar products $(\cdot,\cdot)_{\mathfrak{g},g}$
($g\in G$), defines a Riemannian metric on $G$. Note that this
metric is neither left nor right invariant.

For a fixed $x\in X$, let us compute the derivative
$\rd\varphi_{x}(g):T_{g}G\rightarrow T_{g\cdot x}X$ of the map $G\ni
g\mapsto \varphi_{x}(g)=g\cdot x\in X$. As in the previous section,
we identify the tangent space $T_{g}G$ with the Lie algebra
$\mathfrak{g}$ by right translations; under this identification,
\begin{equation}\label{eq:dphi}
\rd\varphi _{x}(g)=P_{g\cdot x},\qquad g\in G,\ \ x\in X.
\end{equation}

Observe that $P_{g\cdot x}:\mathfrak{g}\rightarrow
\mathfrak{x}_{g\cdot x}$ is an orthogonal projection (with respect
to the scalar product $(\cdot,\cdot)_{\mathfrak{g},g}$ on
$\mathfrak{g}$). Therefore,
\begin{equation}\label{phi-bound}
\|\rd\varphi_{x}(g)\| \leq 1
\end{equation}
and $\rd\varphi _{x}(g)^{\ast }:\mathfrak{x}_{g\cdot x}\rightarrow
\mathfrak{g}$ is an isometry. Moreover, it follows from
(\ref{decomp1}) that
\begin{equation}\label{Pr-Ad}
\rd\varphi_{x}(g)=
\Ad_{g}\circ P_{x}\circ \Ad_{g^{-1}}.
\end{equation}
Considering $\rd\varphi_{x}(\cdot)V$ as a map from $G$ to
$\mathfrak{g}$ (via the embedding $\mathfrak{x}_{x}\subset
\mathfrak{g}$), we obtain
\begin{equation}\label{d2phi}
\rd^{2}\varphi_{x}(g)(U,V)=(\ad_{U}\circ P_{g\cdot x}-P_{g\cdot
x}\circ \ad_{U})\mypp V,
\end{equation}
for any $U,V\in \mathfrak{g}$. This, together with
(\ref{phi-bound}), implies that
\begin{equation*}
\sup_{x\in X,\,g\in G}\|\rd^{2}\varphi _{x}(g)\| <\infty.
\end{equation*}

Thus, $\varphi_{x}\in C_{b}^{2}(G,X)$ uniformly in $x\in G$, and so
one can apply Theorem \ref{th:C2} provided that conditions
(\ref{beta}) and (\ref{N}) are met. Finally, if the probability
measure $Q$ is absolutely continuous with respect to a left Haar
measure on $\mathfrak{G}$, then (\ref{beta}) and (\ref{N}) can be
specified in terms of the corresponding Radon--Nikodym density. Note
that the norm used in condition (\ref{beta}) is generated in this
case by the special Riemannian structure (\ref{eq:5.15}) on $G$.

\subsection{Other examples}
In this section, we briefly discuss two further examples
illustrating possible ways of constructing cluster distributions
$\eta_x$.

\subsubsection{Manifolds of non-positive curvature.}\label{sec:5.4}

Let $X$ be a complete, path-connected manifold with non-positive
sectional curvature (\textit{Cartan--Hadamard manifold}). In this
case, for every two points $x,y\in X$ there is a unique geodesic
$g_{x,y}(t)$, $t\in[0,1]$, such that $g_{x,y}(0)=x$,
\,$g_{x,y}(1)=y$. Assume in addition that $X$ is simply connected.
It follows from the Cartan--Hadamard theorem that the exponential
map $\exp_{x}\!:T_{x}X\rightarrow X$ is a diffeomorphism for every
$x\in X$ (see, e.g., \cite{Chavel,KN2}).

Choose $x_{0}\in X$, and let
\begin{equation*}
\rd{g}_{x_{0},x}:T_{x_{0}}X\rightarrow T_{x}X
\end{equation*}
be the parallel translation along the geodesic $g_{x_0,x}$. To
deploy the construction of Section \ref{sec:droplet}, we set
$W:=T_{x_{0}}X$ and
\begin{equation*}
\varphi_{x}:=\exp_{x}\circ \,\rd g_{x_{0},x}:\,W\to X.
\end{equation*}
For a given probability measure $Q=Q_{x_0}$ on
$(T_{x_0}X)_0^\infty$, consider the corresponding translated
(push-forward) measures on $(T_{x}X)_0^\infty$,
\begin{equation}\label{eq:varpi}
Q_{x}=\rd{g}_{x_{0},x}^{\ast}Q_{x_0}\qquad x\in X.
\end{equation}
According to a general formula (\ref{eq:x0->x}), we can now define a
family of probability distributions on the space $X$ by
\begin{equation}\label{eq:eta}
\eta_{x}:=\bar{\varphi}_x^{\ast}\myp
Q_{x_0}=\exp_{x}^{\ast}Q_{x},\qquad x\in X.
\end{equation}

\begin{remark} In fact, $X$ is essentially the Euclidean space $\RR^d$ ($d=\dim X$)
with a non-constant metric which defines a family of exponential
maps $\exp_x: \RR^d\to\RR^d$ ($x\in\RR^d$). In this interpretation,
we have $W = \RR^d$ and $\varphi_x= \exp_x\! : W \to X$.
\end{remark}

\begin{remark}
Consider the diffeomorphism
\begin{equation}\label{eq:ix}
i_{x}:=\exp_{x}\circ\,{} \rd g_{x_{0},x}\circ
\exp_{x_{0}}^{-1}:\,\frakX\rightarrow \frakX,
\end{equation}
From (\ref{eq:varpi}), (\ref{eq:eta}) and (\ref{eq:ix}) it follows
that the family of distributions $\{\eta_x\}_{x\in X}$ is
translation invariant in the sense that
$\eta_{x}=i_{x}^{\ast}\eta_{x_{0}}$ \,($x\in X$).
\end{remark}

\subsubsection{Metric spaces.}\label{sec:5.5}

Let $(X,\rho)$ be a metric space, endowed with the natural topology
generated by the open balls $B^0_r(x):=\{x'\in X: \rho(x,x')<r\}$
($x\in X$, $r>0$) and equipped with a (locally finite) reference
measure $\vartheta$.

In this section, we construct an example of a family of probability
measures $\{\eta_x(\rd\bar{y})\}_{x\in X}$ on
$\mathfrak{X}=\bigsqcup_{\myp n}\mynn X^n$, based on a different
idea that avoids using any family of maps $\varphi_x$ as in Sections
\ref{sec:5.1}--\ref{sec:5.3}. To this end, note that by a
radial-angular decomposition (based on Fubini's theorem) we can
represent the $\vartheta$-volume of a (closed) ball $B_r(x):=\{x'\in
X: \rho(x,x')\le r\}$ ($x\in X$) as
\begin{equation}\label{eq:Ball_r}
\vartheta(B_r(x))=\int_0^r \left(\int_{\partial
B_s(x)}\vartheta^{\myp x}_{\rm ang}(\rd y\myp|\myp
s)\right)\vartheta^{\myp x}_{\rm rad}(\rd s)
\end{equation}
where $\partial B_r(x)=B_r(x)\setminus B^0_r(x)$ is the sphere of
radius $r$ centred at $x$, $\vartheta^{\myp x}_{\rm ang}(\rd
y\myp|\myp r)$ is the uniform ``surface'' measure on $\partial
B_r(x)$ induced by the measure $\vartheta(\rd y)$, and
$\vartheta^{\myp x}_{\rm rad}(\rd r)$ is the radial component of
$\vartheta$ as seen from $x$. According to formula
(\ref{eq:Ball_r}), the measure $\vartheta$ can be symbolically
expressed as a skew product
$$
\vartheta(\rd y)=\vartheta^{\myp x}_{\rm ang}(\rd y\myp|\myp
r)\,\vartheta^{\myp x}_{\rm rad}(\rd r)\Bigl|_{r=\rho(x,y)}.
$$

For $x\in X$ and $\bar{y}\in\frakX$, set
$\bar{\rho}(x,\bar{y}):=(\rho(x,y_i))_{y_i\in\bar{y}}\in\frakX$. As
usual, the measure $\vartheta$ can be lifted to the space
$\mathfrak{X}$ by setting
\begin{equation}\label{eq:theta-lift}
\bar{\vartheta}(\rd\bar{y}):={\textstyle\bigotimes\limits_{y_i\in\bar{y}}}
\,\vartheta(\rd{y_i}),\qquad \bar{y}\in\mathfrak{X}.
\end{equation}
Similarly, for each $x\in X$ define
\begin{align}
\label{eq:theta-rad-lift} \bar{\vartheta}^{\myp x}_{\rm rad}(\rd
\bar{r}):={}&{\textstyle\bigotimes\limits_{r_i\in\bar{r}}}
\,\vartheta^{\myp x}_{\rm rad}(\rd r_i),\\
\label{eq:theta-ang-lift} \bar{\vartheta}^{\myp x}_{\rm
ang}(\rd\bar{y}\myp|\myp\bar{r})
:={}&{\textstyle\bigotimes\limits_{y_i\in\bar{y}}}\,\vartheta^{\myp
x}_{\rm ang}(\rd y_i\myp|\myp r_i).
\end{align}

Let us now fix a point $x_0\in X$, and let
$f:\,\RR_{0}^{\infty}\rightarrow \RR_{+}$ be such that
$\int_{\mathfrak{X}} f(\bar{r})\;\bar{\vartheta}^{\myp x_0}_{\rm
rad}(\rd\bar{r})=1$. Then we can construct a family of cluster
distributions by setting, for each $x\in X$,
\begin{equation}\label{eq:eta_x-metric}
\eta_{x}(\rd\bar{y}):=f(\bar{r})\;\bar{\vartheta}^{\myp x_0}_{\rm
rad}(\rd \bar{r})\cdot\frac{\bar{\vartheta}^{\myp x}_{\rm
ang}(\rd\bar{y}\myp|\myp \bar{r})}{\bar{\vartheta}^{\myp x}_{\rm
ang}(\partial \bar{B}_{\bar{r}}(x))\myp|\myp
\bar{r})}\myp\biggl|_{\bar{r}=\bar{\rho}(x,\bar{y})},\qquad
\bar{y}\in\frakX.
\end{equation}
That is to say, under the measure $\eta_x$ a random vector $\bar{y}$
is sampled in two stages: first, a vector $\bar{r}$ of the distances
from $x$ to $\bar{y}$ is sampled with the probability density
$f(\bar{r})$ (with respect to the measure $\bar{\vartheta}^{\myp
x_0}_{\rm rad}$), and then the components $y_i$ of $\bar{y}$ are
chosen, independently of each other, with the uniform distribution
over the corresponding spheres $\partial B_{r_i}(x)$, respectively.
\begin{remark}
By definition (\ref{eq:eta_x-metric}), the measure $\eta_x$ may be
considered as a ``translation'' of the pattern measure $\eta_{x_0}$
from $x_0$ to $x$; however, this is not being done by a push-forward
of $\eta_{x_0}$ under some mapping $\varphi_x$ of the space $X$, as
prescribed by the general recipe of Section \ref{sec:droplet};
instead, we compensate the lack of such a mapping
by using the same statistics of the distances at each point $x\in X$
(prescribed by the pattern distribution $\vartheta^{\myp x_0}_{\rm
rad}$) and by taking advantage of the uniform distribution on the
corresponding spheres, which does not require any further angular
information.
\end{remark}

\begin{remark}
If there is a group $G$ of isometries of $X$ acting transitively,
then we can use the same method as for homogeneous spaces (see
Section \ref{sec:5.3}).
\end{remark}

%
%
%
%
%
%

\appendix
\section*{Appendix}

\section{On a definition of the skew-product measure $\hat{\mu}$}\label{ap:1}

We assume for simplicity that $\mu(\varGamma_{X})=1$. In order to
verify that the measure $\hat{\mu}$ is well defined by expression
(\ref{eq:mu-hat}) (which requires the internal integral in
(\ref{eq:mu-hat-int}) to be measurable as a function of
$\gamma\in\varGamma_X$), we shall construct an auxiliary measure
$\tilde\mu$ on $\varGamma_{X}\times Y^{\infty}$ and show that
$\hat{\mu}$ is its image under a certain measurable map.

Let us fix an indexation
$\mathfrak{i}=\{\mathfrak{i}_\gamma,\,\gamma\in\varGamma_X\}$ in
$\varGamma_{X}$, where $\mathfrak{i}_\gamma:\gamma\rightarrow\NN$ is
a bijection for each $\gamma \in \varGamma_{X}$. Define
\begin{equation*}
\varGamma_{X,1}:=\{(\gamma,x)\in \varGamma_{X}\times X \!: \,x\in
\gamma\}.
\end{equation*}
The indexation $\mathfrak{i}$ defines a natural bijection
\begin{equation}\label{eq:i}
\varGamma_{X,1}\ni(\gamma,x)\mapsto
(\gamma,\mathfrak{i}_\gamma(x)\in \varGamma_{X}\times \NN.
\end{equation}
Moreover, the indexation $\mathfrak{i}$ can be constructed so that
bijection (\ref{eq:i}) is measurable (see \cite{VGG}). This ensures
that the map
\begin{equation}
\varGamma_{X}\ni\gamma \mapsto
j_{k}(\gamma):=\mathfrak{i}_{\gamma}^{-1}(k)\in X
\end{equation}
is measurable for each $k\in \NN$.

Consider a family $\{\nu^\gamma,\,\gamma\in\varGamma_X\}$ of
measures on $Y^{\infty}$ defined by
\begin{equation*}
\nu^{\gamma}(\rd\bar{y}):={\textstyle\bigotimes\limits_{k\in\NN}}\,
\eta_{j_{k}(\gamma)}(\rd\bar{y}),\qquad \bar{y}\in Y^\infty.
\end{equation*}
If $A\in \mathcal{B}(Y^{\infty})$ is a cylinder set, $A=A_{1}\times
\dots\times A_{n}\times Y\times\cdots$, then
\begin{equation*}
\nu^{\gamma}(A)={\textstyle\bigotimes\limits_{k=1}^{n}}\,
\eta_{j_{k}(\gamma)}(A_{k}).
\end{equation*}
The function
$\varGamma_{X}\ni \gamma \mapsto \nu^{\gamma}(A)\in \RR$
is measurable due to the measurability of $j_{k}(\gamma)$ and
Condition~\ref{condition1}. Hence, the measure
\begin{equation*}
\tilde{\mu}(\rd\gamma\times\rd\bar{y}):=\nu^{\gamma}(\rd{y})\,\mu(\rd\gamma),\qquad
(\gamma,\bar{y})\in\varGamma_{X}\times Y^{\infty},
\end{equation*}
is well defined.

Finally, a direct calculation shows that the measure $\hat{\mu}$
defined by (\ref{eq:mu-hat}) can be represented as
$\hat{\mu}=\mathcal{I}^{\ast}\tilde{\mu}$, where
$\mathcal{I}:\varGamma_{X}\times Y^{\infty}\rightarrow
\varGamma_{X}\times Y^{\gamma}$
is a measurable map defined by
$(\gamma ,(y_{k})_{k\in \mathbb{N}})\mapsto
(\gamma,(y_{j_{k}(\gamma)})_{k\in \mathbb{N}})$.
This proves the result.

\section{Correlation functions}\label{ap:2}
For a more systematic exposition and further details, see the
classical books \cite{Ge79,Preston,Ruelle}; more recent useful
references include \cite{AKR2,KunaPhD,KKS98}.

\begin{definition}\label{def:corr}
Let $\mu$ be a probability measure on the generalised configuration
space $\varGamma_X^\sharp$, and let $\theta$ be a (locally finite)
measure on $X$. Then the correlation function
$\kappa_{\mu}^{n}:X^{n}\rightarrow \RR_{+}$ of the $n$-th order
($n\in \mathbb{N}$) of the measure $\mu$ with respect to $\theta$ is
defined by the following property: for any function $\phi\in
\mathrm{M}_{+}(X^{n})$ symmetric with respect to permutations of its
arguments, it holds
\begin{multline}
\int_{\varGamma_X^\sharp} \sum_{\{x_{1}\myn,\dots,\myp
x_{n}\}\subset \gamma}\phi (x_{1},\dots, x_{n})\,\mu(\rd\gamma)
\label{corr-funct} \\
=\frac{1}{n!}\int_{X^{n}}\phi
(x_{1},\dots,x_{n})\,\kappa_{\mu}^{n}(x_{1},\dots,x_{n})\,\theta(\rd{x}_{1})\cdots
\theta(\rd{x}_{n}).
\end{multline}
\end{definition}

\begin{remark}\label{rm:x=x} Note that possible multiple points on
the configuration $\gamma\in\varGamma_X^\sharp$ will lead
correspondingly to some coinciding points among
$\{x_1,\dots,x_n\}\subset\gamma$ on the left-hand side of formula
(\ref{corr-funct}) (cf.\ our convention on the use of set-theoretic
notations, see Section \ref{sec:2.1}).
\end{remark}
By a standard approximation argument, equation (\ref{corr-funct})
can be extended to any (symmetric) functions $\phi\in
L^1(X^n,\theta^{\otimes n})$.

\begin{condition}\label{conditionB1}
Correlation functions $\kappa_{\mu}^{m}(x_{1},\dots,x_{m})$ up to
the $n$-th order ($n\in \mathbb{N}$) of the measure $\mu$ with
respect to $\theta$ exist and are bounded.
\end{condition}

\begin{remark} Formula (\ref{corr-funct}) with $n=1$ and
$\phi(x)={\bf 1}_B(x)$ for $B\in\mathcal{B}(X)$ shows that Condition
\ref{conditionB1} automatically implies that $\mu$-a.a.\
configurations $\gamma$ are locally finite.
\end{remark}

\begin{lemma}
Assume that Condition \textup{\ref{conditionB1}} is satisfied with
some $n\in\NN$. Then $\mu \in
\mathcal{M}_{\theta}^{n}(\varGamma_{X})$.
\end{lemma}

\proof Similarly as in the proof of Lemma \ref{moments}, we obtain
(cf.\ (\ref{eq:multi}))
\begin{align}
\notag \int_{\varGamma_{X}} |\langle
f,\gamma\rangle|^{n}\,\mu(\rd\gamma)& \leq
\int_{\varGamma_{\frakX}}\left( \sum_{x\in
\gamma} |f(x)|\right)^{n}\,\mu(\rd\gamma)\\
\label{eq:multiB}& =\sum_{m=1}^{n}\int_{\varGamma_{\frakX}}
\sum_{\{x_{1},\dots,\myp x_{m}\}\subset
\gamma}\phi_{n}(x_{1},\dots,x_{m}) \,\mu(\rd\gamma),
\end{align}
where $\phi_{n}(x_{1},\dots ,x_{m})$ is a (symmetric) function given
by expression (\ref{eq:varPsi}).
%
Note that, by definition of the correlation functions (see
(\ref{corr-funct})), the integral on the right-hand side of
(\ref{eq:multiB}) is reduced to
\begin{equation}\label{eq:k==}
\frac{1}{m!}\int_{X^{m}}\phi_{n}(x_{1},\dots,x_{m})\,
\kappa_{\mu}^{m}(x_{1},\dots ,x_{m})\,\theta(\rd{x}_{1})\cdots
\theta (\rd{x}_{m}).
\end{equation}
By Condition \ref{conditionB1}, $\kappa_{\mu}^{m}\leq c_{m}$
($m=1,\dots ,n$) with some constant $c_{m}<\infty$. Hence, the
integral in (\ref{eq:k==}) is bounded by
\begin{equation}\label{4.3}
\sum_{\substack{i_{1},\dots,\myp i_{m}\geq 1\\
i_{1}+\dots +\myp i_{m}=\myp n}}\frac{c_{m}\myp n!}{i_{1}!\cdots
i_{m}!}\prod_{j=1}^{m}\int_{\mathcal{Z}}|f(x_{j})|^{i_{j}}\,
\theta(\rd{x}_{j})<\infty,
\end{equation}
since each integral in (\ref{4.3}) is finite owing to the assumption
$f\in \bigcap_{\myp 1\le q\le n}L^{q}(X,\theta)$. As a result, the
integral on the left-hand side of (\ref{eq:multiB}) is finite, and
the lemma is proved.
\endproof

\section{Integration-by-parts formula for push-forward measures}\label{ap:3}

For any Riemannian manifolds $\mathcal{W}$ and $\mathcal{Y}$, denote
by $C_{b}^{2}(\mathcal{W},\mathcal{Y})$ the space of twice
differentiable maps $\phi:\mathcal{W}\rightarrow \mathcal{Y}$ with
globally bounded derivatives $\rd\phi$, $\rd^2\phi$. In particular,
for any $\bar{w}\in\mathcal{W}$, the first derivative
$\rd\phi(\bar{w})$ is a bounded linear operator from the tangent
space $T_{\bar{w}}\mathcal{W}$ to the tangent space
$T_{\phi(\bar{w})}\mathcal{Y}$. In what follows, we fix $\phi \in
C_{b}^{2}(\mathcal{W},\mathcal{Y})$.  Note that the adjoint operator
$\rd\phi(\bar{w})^{\ast}: T_{\phi(\bar{w})}^{\ast}\mathcal{Y}\to
T_{\bar{w}}^{\ast}\mathcal{W}$ can be identified with a bounded
operator from $T_{\phi(\bar{w})}\mathcal{Y}$ to
$T_{\bar{w}}\mathcal{W}$ via the scalar products in the tangent
spaces $T_{\bar{w}}\mathcal{W}$ and $T_{\phi(\bar{w})}\mathcal{Y}$
(defined by the Riemannian structure of the manifolds $\mathcal{W}$
and $\mathcal{Y}$, respectively). Furthermore, define
$\Vect_{b}^{1}(\mathcal{W})$ as the space of differentiable vector
fields on $\mathcal{W}$ with a globally bounded first derivative.
\begin{definition}
We say that a probability measure $Q(\rd\bar{w})$ on $\mathcal{W}$
satisfies an \textit{integration-by-parts} (\textit{IBP})
\textit{formula} if for any vector field $V\in
\Vect_{b}^{1}(\mathcal{W})$ there is a function $\beta_{Q}^{V}\in
L^{1}(\mathcal{W},Q)$ (\textit{logarithmic derivative} of $Q $ in
the direction $V$) such that, for any $g\in C_{b}^{1}(\mathcal{W})$,
the following identity holds
\begin{equation}\label{IBP}
\int_{\mathcal{W}}(\nabla
g(\bar{w}),V(\bar{w}))_{T_{\bar{w}}\mathcal{W}}\:Q(\rd\bar{w})
=-\int_{\mathcal{W}}g(\bar{w})\mypp\beta_{Q}^{V}(\bar{w})\,Q(\rd
\bar{w}).
\end{equation}
Whenever it exists, the function $\beta_{Q}^{V}$ can be represented
in the form
\begin{equation}\label{LD}
\beta_{Q}^{V}(\bar{w})=(\beta
_{Q}(\bar{w}),V(\bar{w}))_{T_{\bar{w}}\mathcal{W}}+\Div
V(\bar{w}),\qquad \bar{w}\in\mathcal{W},
\end{equation}
where $\beta_{Q}$ is a vector field on $\mathcal{W}$ (called the
\textit{vector logarithmic derivative} of $Q$) satisfying
\begin{equation*}
\int_{\mathcal{W}}|\beta_{Q}(\bar{w})|_{T_{\bar{w}}\mathcal{W}}\:Q(\rd
\bar{w})<\infty.
\end{equation*}
\end{definition}

Consider the push-forward measure $\eta:=\phi^{\ast}\myp Q$ on
$\mathcal{Y}$, and denote by $\mathcal{I}_{\phi}$ the operator
acting on functions $f:\mathcal{Y}\rightarrow \RR$  by the formula
\begin{equation*}
\mathcal{I}_{\phi }f=f\circ \phi.
\end{equation*}
Because of the definition of the measure $\eta$, the operator
$\mathcal{I}_{\phi}$ is an isometry from $L^{r}(\mathcal{Y},\eta)$
to $L^{r}(\mathcal{W},Q)$, for any $r\in \lbrack 1,\infty ]$. Hence,
the adjoint operator defines an isometry between the corresponding
dual spaces,
\begin{equation*}
\mathcal{I}_{\phi}^{\ast}: L^{r}(\mathcal{W},Q)^{\prime}\rightarrow
L^{r}(\mathcal{Y},\eta)^{\prime}.
\end{equation*}
Furthermore, for any $r\in (1,\infty )$ we have the isomorphisms
$L^{r}(\mathcal{W},Q)'\cong L^{n}(\mathcal{W},Q)$ and
$L^{r}(\mathcal{Y},\eta)'\cong L^{n}(\mathcal{Y},\eta)$, where
$n=r/(r-1)$ (see, e.g., \cite[Ch.\ II, \S\myp2, p.~43]{Schaefer}).
Since $r>1$ is arbitrary, this means that
$\mathcal{I}_{\phi}^{\ast}$ can be treated as an isometry from
$L^{n}(\mathcal{W},Q)$ to $L^{n}(\mathcal{Y},\eta)$ for any $n>1$.
Moreover, repeating the arguments used in the proof of Lemma
\ref{lm:L1}, it can be shown that the same also holds for $n=1$. To
summarise, for any $n\geq 1$ the operator
\begin{equation*}
\mathcal{I}_{\phi}^{\ast}:L^{n}(\mathcal{W},Q)\rightarrow
L^{n}(\mathcal{Y},\eta)
\end{equation*}
is an isometry.

\begin{theorem}\label{IBP-proj}
Let $\phi\in C_{b}^{2}(\mathcal{W},\mathcal{Y})$ be such that the
operator
\begin{equation*}
\rd\phi(\bar{w})^{\ast }:T_{\phi (\bar{w})}\mathcal{Y}\rightarrow
T_{\bar{w}}\mathcal{W},\ \bar{w}\in \mathcal{W}
\end{equation*}
is an isometry, and suppose that the measure $Q$ satisfies the IBP
formula \textup{(\ref{IBP})}. Then the push-forward measure
$\eta=\phi^{\ast}\myp Q$ satisfies an IBP formula with the
logarithmic derivative
$\beta_{\eta}^{U}=\mathcal{I}_{\phi}^{\ast}\beta_{Q}^{V}$, where
$V=V_U$ is a vector field on $\mathcal{W}$ given by
\begin{equation*}
V(\bar{w})=\rd\phi(\bar{w})^{\ast}\, U(\phi(\bar{w})),\qquad
U\in\Vect_{b}^{1}(\mathcal{Y}).
\end{equation*}
\end{theorem}

\proof Note that $V\in \Vect_{b}^{1}(\mathcal{W})$. Applying the IBP
formula (\ref{IBP}), making the change of measure
$\eta=\phi^{\ast}\myp Q$ and taking into account that
$\mathrm{d}\phi (\bar{w})\rd\phi (\bar{w})^{\ast}$ is the identity
operator in $T_{\phi (\bar{w})}\mathcal{Y}$, we see that (\ref{IBP})
holds for $\eta$ with the corresponding logarithmic derivative
$\beta_{\eta}^{U}=\mathcal{I}_{\phi}^{\ast}\beta_{Q}^{V}$. Finally,
the $\eta$-integrability of $\beta_{\eta}^{U}$ follows by the
isometry of $\mathcal{I}^*$.
\endproof




\begin{remark}\label{rm:sqcup}
All of the above remains true in the case where
$\mathcal{W}=\bigsqcup_{i=0}^{\infty }\mathcal{W}_{i}$ and
$\mathcal{Y}=\bigsqcup_{i=0}^{\infty }\mathcal{Y}_{i}$ are countable
disjoint unions of Riemannian manifolds $(\mathcal{W}_{i})$ and
$(\mathcal{Y}_{i})$ respectively, and the mapping $\phi$ acts
component-wise, that is, $\phi:\mathcal{W}_{i}\to \mathcal{Y}_{i}$.
Although the spaces $\mathcal{W}$ and $\mathcal{Y}$ do not possess a
proper Riemannian manifold structure, all notions introduced above
(including the IBP formula (\ref{IBP})) can be understood
component-wise, and we use the analogous notations without further
explanations.
\end{remark}

\section*{Acknowledgements}
Research of L.\,Bogachev was supported in part by a Leverhulme
Research Fellowship. Partial support
from DFG Grant 436\,RUS\,113/722 and SFB701 is gratefully
acknowledged. The authors would like to thank Sergio Albeverio, Yuri
Kondratiev, Eugene Lytvynov, Stanislav Molchanov and Chris Wood for
helpful discussions.


\end{document}